# Husserl, Cantor & Hilbert

## *La Grande Crise des Fondements Mathématiques du XIX$^{ème}$ Siècle*

Arkady Nedel

## ABSTRACT

Three thinkers of the 19th century revolutionized the science of logic, mathematics, and philosophy. Edmund Husserl (1859-1938), mathematician and a disciple of Karl Weierstrass, made an immense contribution to the theory of human thought. The paper offers a complex analysis of Husserl's mathematical writings covering calculus of variations, differential geometry, and theory of numbers which laid the ground for his later phenomenological breakthrough. Georg Cantor (1845-1818), the creator of set theory, was a mathematician who changed the mathematical thinking per se. By analyzing the philosophy of set theory this paper shows how was it possible (by introducing into mathematics what philosophers call 'the subject'). Set theory happened to be the most radical answer to the crisis of foundations. David Hilbert (1862-1943), facing the same foundational crisis, came up with his axiomatic method, indeed a minimalist program whose roots can be traced back to Descartes and Cauchy. Bringing together these three key authors, the paper is the first attempt to analyze how the united efforts of philosophy and mathematics helped to dissolve the epistemological crisis of the 19th century.

# CONTENTS

*Il n'y a pas de certitude ;*
*il y a seulement des hommes certains.*

Charles Renouvier

## CHAPITRE 1. LES FONDEMENTS RETROUVES

Commençons par une prolepse. Le parcours d'Edmund Husserl, l'inventeur de la phénoménologie, peut être séquencé en deux grandes périodes : la première s'étendrait de 1887 à 1913 (de son habilitation jusqu'à la publication de son *magnum opus, Ideen I*) ; la seconde de 1913 à 1936, année de la *Krisis*, son testament philosophique, où il analyse entre autres son propre parcours intellectuel. Ces deux étapes émaillent une vie extrêmement intense et parfois tourmentée, celle d'un grand esprit toute entière consacrée à la recherche de la vérité – ou de la Vérité –, en l'existence de laquelle il n'a jamais cessé de croire. On pourrait fort aisément entreprendre la rédaction d'une biographie détaillée de Husserl (ce qui n'a d'ailleurs pas encore été fait[1]) montrant son cheminement pas à pas. Il en résulterait un récit dépourvu de faits extraordinaires, rares dans la vie d'un ascète philosophique. Si quelqu'un relevait la gageure d'écrire un roman sur Husserl, celui-ci ressemblerait à celui

[1] Excepté un portrait du philosophe fait par Malvine Husserl [Husserl 1988], cette lacune est en partie comblée par des ouvrages historiques : [Spiegelberg 1960 ; Landgrebe 1963 ; Bakker 1969 ; Avé-Lallemant 1975 ; Schuhmann 1977 ; DeBoer 1978 ; Sepp 1988]

d'un Proust ou d'un Joyce, où le temps externe serait scellé dans le temps interne et tous les événements se déploieraient à l'intérieur du héros. Bref, ce serait un roman sur la vie d'une conscience, d'un λόγος frayant avec des chimères et des spectres, des sens purs incompréhensibles aux esprits quotidiens ; une conscience qui n'éclot que dans l'étrange monde des phénomènes, qu'elle extrait de sa présence courante. Peu avant sa disparition le 27 avril 1938, Husserl dresse pour lui-même le constat d'une vie heureuse, digne d'un philosophe, et ajoute qu'il n'a à se plaindre de rien. Il s'apprête à finir son histoire comme il le doit. Un seul bémol : ses projets inachevés. Maurice Natanson nomme à juste titre Husserl le « philosophe des tâches infinies » [Natanson 1973], mais sous réserve que ces infinités aient été inclues dans son travail dès le début.

L'objectif de cette étude n'est ni historique ni génétique[1]. Nous nous bornons à une esquisse du contexte husserlien en nous attachant davantage à la logique des causes qui ont mené la pensée à la phénoménologie qu'à l'ordre chronologique. Le « contexte » signifie ici le réseau d'influences, la contiguïté des idées, aussi bien synchroniques que diachroniques, que nous essaierons de reconstruire – sans prétendre à aucun moment expliquer en détails l'itinéraire intellectuel de Husserl. En lisant les textes husserliens, nous ne leur appliquerons pas la méthode scolastique ou lecture « montante » : *littera* → *sensus* → *sententia* [Grabmann 1957] (un tel exercice exigerait une compétence philologique spécifique). Nous nous contenterons de procéder à la reconstitution des épisodes jouant un rôle majeur dans la production de la conscience pure absolue qui est *achevée* dans *Ideen I*. Notre limite temporelle sera donc l'année 1913, où la nouvelle méthode atteint sa forme la plus aboutie.

A la différence des premiers essais visant à construire une « science descriptive » – notamment des *Recherches logiques* (1901/02), où Husserl était encore sous l'empire des illusions psychologiques –, la phénoménologie tardive (ou doctrine du sujet transcendantal) n'est pas tant une description qu'une pratique. La conscience absolue, lorsqu'elle saisit les essences, ne peut à proprement parler *décrire* les choses de façon objective car elle est elle-même objectivité ultime qui ne voit que ses propres actes. L'objectivité n'a de sens que si la pensée est dirigée vers l'extérieur, vers un monde situé dans un ailleurs et restant d'une manière ou d'une autre inconnu à la conscience ; plus simplement, le concept de l'objectivité est indispensable pour décrire toute existence. Le sujet transcendantal ne voit pas l'existence et les choses comme existantes mais comme vides. L'appel husserlien à revenir aux « choses en soi » (*zu den Sachen selbst*)[2] s'explique par le projet de libérer la conscience de son nécessaire lien à une objectivité, toujours trompeuse si elle appartient à notre monde. C'est un appel à employer la conscience – à mener une vie de noéticien – de sorte à la comprendre comme la réalité ultime ; en d'autres termes, à parvenir à un niveau où la chose n'est pas telle qu'elle apparaît selon l'expérience mondaine, utile ici et maintenant, mais réduite à son état idéal. La chose en soi, c'est aussi un acte de la pensée saisi au moment du penser.

Quand nous parlons de la praxis de la conscience du sujet transcendantal, il s'agit des sens surgissant de la relation entre des essences vues phénoménologiquement. *Vues*, c'est-à-dire naturalisées par la conscience en tant que ses propres éléments constitutifs, ce qui rend cette dernière absolue. La

[1] Pour une étude génétique voir [Lavigne 2004]
[2] Pour les débats contemporains voir [Sallis 1984]

« vison d'essence », l'une des expressions-clés du Husserl tardif, est le moyen de voir les objets existants dans leur non-existence, dépouillés de leur apparition. Le vrai sens des choses n'est pas celui désigné par leur nom et délimitant la sphère de leur applicabilité, mais celui qui reste même quand la chose (πρᾶγμα) perd son objectivité. Husserl montre qu'il suffit de soustraire à un objet son corps empirique par l'abstraction de premier ordre et tout discours sur son caractère objectif change radicalement. Le nouveau corps abstrait de l'objet ne correspond plus à l'objectivité sur laquelle repose la pensée ordinaire (rappelant ici Duns Scot, chez qui « abstrahitur ab materia » est la première étape). De cette opération d'abstraction naît une autre sémantique qui n'appartient plus à la classe des choses comprenant une donation objective. L'on parle dans ce cas d'objets idéaux n'ayant qu'une vague ressemblance avec les phénomènes : les premiers n'existent qu'en mathématiques, tandis que les phénomènes sont l'idéalisation de l'idée même de l'existence.

Il est difficile d'imaginer quelle aurait la phénoménologie si Husserl n'avait pas reçu sa formation mathématique. Grâce à celle-ci a lieu sa première rencontre avec les nombres naturels, entités idéales *par excellence* ; l'une des plus belles énigmes de l'époque depuis près de cent ans, les nombres suscitaient plus de questions que de solutions. « C'est par l'analyse de [ce concept] que toute philosophie des mathématiques doit commencer » [Hua XII], dit-il en 1887 dans sa thèse d'habilitation.

Il va de soi que les entités idéales ont toujours posé problème ; dès la haute scolastique, pour ne pas remonter plus loin dans le temps, ce thème fait l'objet d'une attention assidue chez les philosophes, en particulier ceux traitant de la nature du continuum et des universaux : *sententia vacum* vs. *sententia rerum*. Pour nombre de ces auteurs, philosophie et mathématiques sont indissociables et complémentaires. Jean Buridan, Marsile d'Inghen (Nijmegen), William Ockham ou Thomas Bradwardine, les *magistri nominales*, pensent tous que le continuum est infini aussi bien par son étendue que sa profondeur ; il n'est pas d'éléments plus petits qui le constituent. Certains, tels Grégoire de Rimini ou Gauthier Burley, défendent quant à eux l'existence de ces éléments en les plaçant hors de l'*intellectus humanus*. Là surgit le schisme philosophique, dont on trouve les suites aussi bien dans la pensée de Georg Cantor que dans la phénoménologie husserlienne. Si Cantor prend nettement le parti réaliste, voyant comme Burley dans la proposition mentale une construction d'objets réels – Cantor croit que les nouvelles séries de nombres qu'il a découvertes existent *in esse reali* –, Husserl est en revanche un nominaliste convaincu. A l'instar de Marsile d'Inghen, il s'intéresse au savoir universel dont les lois sont seules légitimes ; cependant, ce savoir et ses lois ne peuvent résider que dans notre conscience. Pour les découvrir, il faut faire oublier tout ce qui n'appartient pas à cette conscience-même, tout ce qui ne procède d'elle ni ne caractérise son essence. Comme Ockham, Husserl traite les concepts universaux *qua* structures purement mentales naturalisées par le travail phénoménologique de la conscience sur soi-même. Le but du travail consiste à atteindre le stade de la conscience où toute question sur la réalité de ses actes perdrait son sens. Cette réalité doit devenir la seule.

La pensée de Husserl doit beaucoup à sa formation initiale. L'essentiel n'est pas tant que les mathématiques inculquent la rigueur mais servent au penseur de tremplin vers ses recherches. Jusqu'au bout de son raisonnement, l'esprit husserlien reste scientifique ; dans la *Krisis* il ne se lasse pas de répéter que la phénoménologie transcendantale (synonyme de la philosophie), rejetant les défaillances du passé, est une science universelle qui définit ses propres fondements avant que de traquer la vérité. En effet, la « crise » chez Husserl est au premier chef une crise des fondements dont les causes complexes datent de l'époque moderne, où le sujet du savoir prit sa forme finale. Cette crise, à laquelle Husserl offre sa solution de phénoménologue transcendantal en proposant un nouveau sujet (le « transsujet ») qui ne pense plus en termes de dichotomie subjectif-objectif, apparaît de manière plus évidente dans la mathématique de la seconde moitié du XIX$^{ème}$ siècle. La *Krisis* est le résumé d'une cinquantaine d'années de méditation sur ce problème. Walter Biemel note justement que l'univers phénoménologique créé par le philosophe aide le sujet européen à comprendre que cette crise est en fait la sienne [Hua VI : XX]; c'est lorsqu'il en prend connaissance que le sujet devient lui-même. Vie et savoir sont entrelacés.

Les phénomènes ne sont pas identiques aux entités mathématiques ; le mode de raisonnement d'un algébriste n'est pas semblable à celui d'un noéticien. C'est pourquoi il ne faut guère croire que Husserl ne fait que développer les matières apprises pendant ses études. Au contraire, son parcours est bien plus riche et complexe : ses toutes premières réflexions sur les mathématiques permettent déjà d'exclure une croyance aveugle de sa part dans cette science. Il ne faut pas non plus s'imaginer que Husserl se limite à transformer les concepts mathématiques en phénomènes, triviale démarche du « technicien » (comme il le dit dans les *Recherches logiques*), à surimposer une abstraction à une autre ; nous verrons par la suite que les procédés d'abstraction dans les deux sciences se distinguent. En revanche, le thème des derniers fondements domine la pensée husserlienne pratiquement dans tous ses écrits, variant avec le temps et acquérant des formes de plus en plus rigoureuses. Tout comme les mathématiques, la philosophie en est arrivée à être problématique pour elle-même. Si, dans la *Philosophie comme science rigoureuse* (1910/11)[1], Husserl diagnostique l'absence de méthodes et de théories logiquement fondées en philosophie, un quart de siècle plus tard il se pose encore la même question : comment une philosophie scientifique est-elle possible [Hua VI : 29] ?[2] La grande méfiance envers la philosophie, si manifeste dans la seconde moitié du XIX$^{ème}$ siècle, a déteint sur la pensée de la période suivante : son mode d'entendement et son moyen de produire des vérités seraient faibles, c'est-à-dire qu'ils ne reposent ni sur les donnés empiriques ni sur l'exactitude de la méthode. D'après Eugen Dühring et Oskar Bogan, pour ne citer qu'eux, la philosophie hors de la science est « un anachronisme ». Oswald Külpe (1862-1915), assistant de Wilhelm Wundt à Leipzig qui connaissait les travaux de Husserl par Carl Stumpf, note lui aussi qu'au XIX$^{ème}$ siècle, grâce au développement rapide des sciences naturelles, l'écart entre celles-ci et la philosophie est sans doute pour la première fois devenu si important [Külpe 1920 : 10] [3]. Dans le contexte de l'époque si la philosophie veut se réhabiliter en tant que discipline critique, comme le note Herbart Schnädelbach [Schnädelbach 1984: 131], ou d'autant plus dominer les sciences positives, elle ne peut guère demeurer dans son univers spéculatif mais, au contraire, doit établir un bon contact avec ces sciences, devenir elle-même dans une certaine mesure positive. La science (utile) est celle qui prévoit les problèmes de la nature et dirige l'esprit humain vers leurs solutions [cf. Dühring 1875: 474; Bogan 1890: 11-12]. Bref, la véritable philosophie ne renonce pas à la méthode scientifique.

[1] *Philosophie als strenge Wissenschaft.*

[2] [Wie soll eine "Philosophie", eine wissenschaftliche Erkenntnis von der Welt möglich werden…]

[3] (sur Husserl cf. p. 130-133) ; [Cattaruzza, Sinico 2004]

Voici réponse de Husserl : par la construction d'une méthode capable de nous donner un schéma conceptuel du monde permettant d'établir tous les liens causaux malgré leur nombre infini. On peut entendre ici les voix de deux auteurs déjà cités : Marsile d'Inghen, ne voyant le monde que dans son concept [Inghen 2000], et Cantor, déduisant les nouvelles infinités des nombres de la conception de l'*ensemble* quasi-transcendantal ; ce-dernier renvoie à son tour à la thèse d'Henri de Harclay qui affirme que plaise à Dieu et le monde sera infini et éternel. Les deux infinités ne se valent pas puisque Dieu les crée de façon différente ; il peut également soustraire une partie de chacune dont la structure sera infinie également car dans chacune de ses parties est contenu le tout. Pour cette raison, la comparabilité des infinis est légitime ; cette idée scolastique, rejetée par Bernard Bolzano lui-même, ne sera ressuscitée que dans la théorie des ensembles.

Quant au monde, une fois mis hors-circuit, il est restitué à la conscience par la *Krisis* ; mais c'est une toute autre conscience, radicalement distincte de ses manifestations ordinaires. La phénoménologie n'annihile donc rien. Elle éprouve la robustesse de ses propositions en gardant celles qui ne s'évanouissent pas dans le brouillard empirique.

Dans un texte de 1894 intitulé *Sur les objets intentionnels*, Husserl remarque : « Nous venons de trouver dans la mathématique un appui pour notre conception, et c'est pourtant précisément dans cette science que l'on voudra puiser des contre-arguments » [Hua XXII : 324]. Cet appui – qui sert à purifier « les impropriétés de la pensée » caractérisant les jugements portés sur le monde empirique –, est considéré en cette fin du XIX$^{ème}$ siècle comme la seule issue possible à la crise systémique. La méthode mathématique est sécurisante : toutes les opérations de base sont formelles et donc privées de toute ambiguïté. Par nature donc, elles sont le produit du travail de sens idéaux qui ne dépendent aucunement de l'intuition ; c'est un raisonnement coupé des contextes particuliers et qui s'effectue au niveau le plus élevé, abstrait ou, si l'on veut, général. Son contenu est non un résultat concret (2 + 2 = 4) mais l'opération-même ; étant transcendantale, elle ne peut être soumise aux contingences. Le caractère formel des mathématiques défini par Husserl comme le fondement de la science est le principe gnoséologique sur lequel repose tout le reste. Il est irréductible à toute expression plus simple (au début des *Principia philosophiae* Descartes va lui jusqu'à mettre en doute la démonstration mathématique). Ce principe s'oppose à la « mathématique réelle » défendu par ses adversaires ; ces derniers restent limités à la sphère de l'intuition où dominent vérités innées et arguties qui n'accompagnent jamais que des problèmes spécifiques. Idée que, pour être exact, Husserl a moins inventée lui-même qu'empruntée à Descartes : la règle VI des *Regulae* (1628/29)[1] nous enjoint à nous diriger vers les choses les plus simples, voire absolues, car ce sont elles qui renferment les vérités ultimes ou *intimam rerum veritatem* (Règle X). Bien que l'essence de la nature soit simple, cela n'implique pas qu'elle soit accessible. En revanche, le principe cartésien stipule que le simple est le général/absolu : « J'appelle l'absolu tout ce qui contient en soi la nature pure et simple, par exemple tout ce que l'on considère comme indépendant, originel, simple, général, unique, égal, similaire, droit et tout de ce genre. J'appelle aussi l'absolu le plus simple et le plus aisé »[1]. A la conception aristotélicienne du général – καθόλου, qui résulte de l'induction et n'est jamais connaissance directe des choses –, Descartes ajoute que la généralité exprime avant tout l'idée du simple, d'où le plus simple est identifié à l'absolu.

[1] Pour l'édition française, cf. [Descartes 1977]

Ainsi commence la méthode. Du *simplicissimum* et de l'évident, on progresse vers les concepts complexes ; Husserl n'hésite pas à réitérer cette approche dès ses écrits techniques. Par exemple, dans *Le problème fondamental de l'arithmétique et de l'analyse* (1889/90), il affirme qu'afin d'engager les procédures algébriques difficiles, il faut « connaître les règles générales des opérations de calcul (directes et inverses) avant que les signes algébriques ne soient introduits… » [Hua XXI : 235]. La mathématique et la logique se rencontrent dans cette simplicité servant de base pour toutes les recherches qui suivent même si son contenu, c'est-à-dire le nombre des axiomes, n'est pas fixe. En Allemagne, Husserl cite le mathématicien Carl Gauss qui, dans ses recherches, montre la clarté de la méthode inspirée des principes cartésiens[2] ; l'optique gaussienne renvoie à la théorie de la lumière que Descartes développe conjointement à sa géométrie.

La méthode de la clarté et de la simplicité maximale est considérée par Husserl comme propre aux mathématiques ; c'est en cela qu'elle se distingue de toutes les autres méthodes. Les mathématiques deviennent sa méthode : claire, certaine et univoque. Nous verrons que le rêve du philosophe sera de transposer cette démarche à la philosophie – qu'il trouve par trop métaphysique (ou scolastique), trop opaque – pour en faire une « science rigoureuse ». Rien ne nous empêche de rapprocher ces deux activités intellectuelles et, en appliquant « le caractère épistémologique de la mathématique à la philosophie <…>, nous entrons dans une sphère plus subtile et rencontrons un problème plus difficile, celui de mettre ces deux sciences sur un même plan … » [Hua XXI : 243] ; Husserl continue : « si la mathématique est une science, elle doit être capable de démontrer (*beweisen*) effectivement ses propositions » [Hua XXI : 243]. A l'époque où ces mots sont écrits, c'est bien le fait de « démontrer/*beweisen* » qui est au cœur du problème car la science mathématique se trouve alors incapable de présenter ses propres concepts de telle sorte qu'ils soient tout à fait dépourvus des ambiguïtés et résidus psychologiques. Selon Jean Cavaillès, les mathématiques se heurtent à une « double crise de rigueur [qui] provoque une remise en question des principes et des notions fondamentales de la géométrie et de l'analyse » [Cavaillès 1981: 45]. Les mathématiques ont été pénétrées par des structures de l'esprit qui ont déplacé les critères de la science rigoureuse vers la position du sujet ; il s'agit de vérités exprimées dans le langage formulaire mais fondées sur l'intuition. Léon Brunschvicg pense que les nombres imaginaires (ou complexes : *i*) [Brunschvicg 1912], qui ne résident dans aucune réalité objective, sont une des causes de cette crise ontologique en mathématiques. Sur le plan historique, ces étrangetés algébriques n'ont pas été étudiées avant le XVIII$^{\text{ème}}$ siècle où, grâce aux travaux de Gottfried Leibniz et Leonhard Euler, la science avait atteint un niveau de calcul plus avancé.

[1] [Absolutum voco, quidquid in se continet naturam puram et simplicem, de qua est quaestio, ut omne id quod consideratur quasi independens, causa, simplex, universale, unum, aequale, simile, rectum, vel alia hujusmodi ; atque idem primum voco simplicissimum et facillimum…]

[2] On voit ainsi l'influence cartésienne sur l'optique géométrique de Gauss, notamment la loi de Descartes sur la réflexion : « Tout rayon lumineux incident frappant une surface réfléchissante avec un angle $\alpha$ par rapport à la normale est réfléchi dans le plan défini par le rayon incident et la normale, au delà de celle ci, et avec un angle $\beta = \alpha$ par rapport à la normale ».

En effet, la cause de cette crise est intrinsèque à l'objet mathématique ; non seulement l'objet devient irréel (imaginaire, complexe), mais il risque aussi de devenir insaisissable par les moyens de cette science exacte. Nicolas Bourbaki remarque à juste titre que « les mathématiciens ont été persuadés qu'ils démontrent des « vérités » ou des « propositions vraies » ; une telle conviction ne peut évidemment être que d'ordre sentimental ou métaphysique, et ce n'est pas en se plaçant sur le terrain de la mathématique qu'on peut la justifier, ni même lui donner un sens qui n'en fasse pas une tautologie. L'histoire du concept de vérité en mathématique relève donc de l'histoire de la philosophie et non de celle des mathématiques » [Bourbaki 2007 : 21]. Pire encore, les concepts les plus basiques, les entités abstraites (nombre, limite, etc.) ont perdu leur crédibilité, sur laquelle reposait le savoir mathématique. La question essentielle est maintenant formulée ainsi : si le concept du nombre est inexact, s'il ne garantit plus l'évidence qui engendre la possibilité-même de construire la vérité, alors où faut-il chercher l'exactitude ?

Les algébristes se trouvent à ce moment devant un problème kantien : trouver le fondement transcendantal de leur science, en trouver la primale évidence. Leibniz est l'un des premiers à soulever le problème que des vérités aussi évidentes que $2 + 2 = 4$ s'avèrent non moins susceptibles d'être remises en question lorsqu'on réfléchit aux définitions des éléments et de leur démonstration. Leibniz met à l'épreuve non pas le nombre en tant que tel mais son existence, qui semblait auparavant découler logiquement de sa nature. Il ne s'agit pas d'un simple doute de mathématicien, lequel accompagne parfois la découverte d'une solution inattendue, mais d'un défaut ontologique. La vérité de la phrase mathématique, même la plus basique, n'est pas indubitable jusqu'au moment où ses éléments sont définis de telle manière que cette définition prouve leur existence. Le nombre entier, l'une des plus grandes préoccupations de la mathématique de l'époque, reste indéfini. Après avoir ouvert la boîte de Pandore, Leibniz s'en tient là et ne pousse pas ses réflexions plus loin. Karl Weierstraß, le père de l'Analyse du XIX$^{\text{ème}}$ siècle et *großer Lehrer* pour Husserl [Schuhmann 1977 : 7 ; de même Ahrens 1907], dont les travaux contribuèrent considérablement à conforter le point de vue algébrique, voyait la logique de la théorie des nombres dans une théorie rigoureuse des nombres réels. Dans son *Introduction à l'analyse* de 1859-1860, puis dans son cours sur le calcul intégral de 1863-1864, Weierstraß traite la question des fondements en penchant toujours pour la rigueur.

Des pas dans cette direction sont esquissés par Hermann Grassmann qui, en 1861, propose une définition de l'addition et de la multiplication des nombres entiers et démontre leurs propriétés fondamentales (associativité, distributivité, etc.). Quant à la mathématique, elle est selon Grassmann une « science de l'être particulier qui n'apparaît qu'à travers la pensée ou en tant que forme pure de pensée » [Grassmann 1894: 23][1]. Ses idées n'ont pas reçu

[1] Grassmann considère son *Ausdehnungslehre*, sur lequel il travaille à partir des années 1840, comme une tentative de fonder une nouvelle science ; cf. p. 46. Dans un autre texte où il définit la mathématique comme « la philosophie des sciences au sens strict du terme » (ce que Husserl

l'accueil mérité en raison de son modeste statut d'*Oberlehrer* au collège de Stettin : son chef-d'œuvre devait longtemps rester « scellé par sept sceaux » [Dieudonné 1979 : 5][1]. Il faut attendre la fin des années 1880 pour que Richard Dedekind et Giuseppe Peano[2] dotent l'arithmétique d'un système d'axiomes. C'est ici que le doute leibnizien se transforme en une ruée frénétique vers les fondements, quand il devient clair – en premier lieu aux auteurs de cette axiomatisation eux-mêmes – que leurs définitions ne possèdent pas de qualité transcendantale. La « question kantienne » demeure encore irrésolue.

Tous comprennent que pour faire de la mathématique une science exacte, un étalon de la rigueur scientifique, l'intuition kantienne est insuffisante. Son incomplétude réside dans son impossibilité à être exprimée par un langage symbolique. Nous sommes incapables de définir l'intuition en termes formels ou algébriques, voire de définir les limites de son applicabilité. L'intuition est toujours rattachée à un sujet, elle est inséparable non seulement de la sphère de la raison « personnelle » mais aussi de l'état subjectif *eo ipso*. Certes, l'intuition ne peut remédier à l'incertitude épistémologique pour la raison qu'elle est principalement non-formalisable. Dans la *Critique de la raison pure* (1781), Kant instrumentalise l'intuition : elle ouvre le sujet vers le monde à travers une objectivité de l'espace-temps caractérisant la nature de l'être humain. L'intuition chez Kant est un état du sujet avant qu'il ne prenne connaissance du monde, avant qu'il n'entre dans la réalité objective de l'espace de Newton qui réside dans chaque individu comme intuitionnée.

Si j'intuitionne, je ne le fais qu'en tant que moi-sujet qui visualise un objet à ma manière propre. Tout au contraire, le concept mathématique pur doit être libre d'expérience subjective de toute sorte ; il est transcendant par nature, comme la *scientia* – au sens de Scot, qui définit le mot *transcendens* comme appartenance à la vérité ultime et éternelle. La majorité des mathématiciens du XIX$^{\text{ème}}$ siècle s'adressent à l'intuition dans leurs recherches, mais comme à un point de départ ; la croyance commune que les concepts fondamentaux des mathématiques existent est pour eux bien compatible avec l'intuition. Tout travail sur une théorie abstraite commence par la reconnaissance du simple fait

---

répétera à plusieurs reprises dans ses écrits de jeunesse), Grassmann va jusqu'à accorder à celle-ci le rôle d'un *Geist* (hégélien) qui, en étendant ses propres limites, chemine vers les principes universels ; cf. [Grassmann 1878 : XXXI]

[1] Sa seule biographie disponible à ce jour : [Petsche 2006] A propos de l'influence grassmannienne sur les recherches des fondements, chez Peano en particulier, voir son article « Calcolo geometrico secondo l'*Ausdehnungslehre* di H. Grassmann, preceduto dalle operazioni della logica deduttiva » (1888) [Peano 1958 : 3-19] ; aussi [Freguglia 2011]

[2] Précisément : *Was sind und was sollen die Zahlen ?* (1888, cf. infra) de Dedekind et *Arithmetices Principia, nova methodo exposita* (1889) de Peano. Bien des historiens considèrent ces deux ouvrages, parus à très peu d'intervalle, comme les efforts convergents d'axiomatiser l'arithmétique. Philip Jourdain nous informe que Peano était au fait des recherches dedekindiennes : cf. [Jourdain 1912]. Néanmoins la technique de Peano se distingue de celle de ses collègues. Dans son article *Super Theorema de Cantor-Bernstein,* il explique sa méthode en soulignant la différence qui le sépare de Russell : « notre analyse des principes de ces sciences [arithmétique et géométrie] consiste en la réduction des postulats communs au nombre minimal des conditions nécessaires et suffisantes <…> *Nous avons l'idée du nombre donc le nombre existe* » ; cf. [Peano 1906 : 365] (Nous soulignons)

Enfin, le chemin vers l'axiomatisation parfaite sera repris par des disciples de Peano, tels que Mario Pieri, Cesare Burali-Forti ou Alessandro Padoa. De la sorte, Pieri se penche non seulement sur le programme de son maître, mais aussi sur les problèmes posés par David Hilbert (par ex., dans *Über den Zahlbegriff,* 1900) où il parle de la consistance de la théorie des nombres. Selon Pieri, la démonstration de cette consistance doit être fondée sur la logique.

que le sujet ne crée pas la vérité mais la découvre. Ce « fait » nous est donné dans notre faculté d'abstraire des choses existantes pour visualiser les choses inexistantes, c'est-à-dire transcendantes. Le concept du nombre, par exemple, que les savants traitent depuis l'Antiquité, vient d'abord – comme le dirait Descartes – sous une forme de l'esprit ; il n'existe pas pour le monde réel, mais ce dernier se fonde sur le nombre. Il y a donc dès le début deux réalités dont le lien est toujours incertain. De ce point de vue, l'intuition est le moyen de lier deux réalités sans aucune évidence formelle. Soulignons qu'il ne s'agit pas ici de l'abysse séparant le sujet qui pense du monde objectif – l'acte du *cogito* et la chose. La différence entre deux réalités se mesure par le degré avec lequel on peut les saisir dans une description exacte. A cet égard, les mathématiques sont l'art de décrire, comme l'illustre le cas de la géométrie analytique : née dans les travaux de Descartes et de Pierre de Fermat[1], cette science, ou plutôt cette méthode, résulte de la nécessité de disposer d'outils standards et précis permettant la résolution de tâches géométriques. La simplicité de cette méthode, qui apparaît avec l'invention des coordonnées, est basée sur le procédé très formel des calculs et n'exige point une grande virtuosité d'esprit.

Si le monde empirique peut être décrit tel qu'il est donné, malgré la grande diversité de ses apparences, les entités abstraites (nombre, point, etc.) ne se laissent pas décrire, étant elles-mêmes les outils de la description. Le problème de cette nature ambivalente des entités mathématiques est saisi dans toute sa complexité au milieu du XIX[ème] siècle. Augustin Cauchy, après avoir travaillé sur les nombres imaginaires et en particulier sur la notation symbolique de la trigonométrie, comprend que la science mathématique exige un modèle plus rigoureux basé sur des procédés logiques incontestables. Cette vison, partagée par d'autres auteurs – à commencer par Bolzano mais aussi Joseph-Louis Lagrange[2] –, peut être résumée ainsi :

Premièrement, tout concept doit se baser sur un autre déjà défini ;

---

[1] Les écrits de Fermat sur la géométrie analytique ne furent pas publiés de son vivant, à l'instar de ceux de Descartes, avec qui Fermat engagea un grand débat sur l'optique géométrique (notamment sur le chemin optique du rayon lumineux) suite à la publication du traité de celui-ci en 1637. Fermat suggère que la lumière se propage d'un point à un autre selon des trajectoires telles que la durée du parcours soit extrémale. On retrouve ainsi la plupart des résultats de l'optique géométrique, en particulier les lois de la réflexion sur les miroirs, les lois de la réfraction, etc. Dès 1636, Fermat entre en correspondance avec le Père Mersenne et s'enquiert dès sa première lettre des nouveautés mathématiques depuis les cinq dernières années. Pour plus de détails voir [Serfati, Descotes 2008]

[2] En 1797 Lagrange publie son livre *Fonctions analytiques* qui prétend fonder le calcul sur une base rigoureuse ou – selon l'auteur – sur la véritable métaphysique de ses principes. Bolzano et Cauchy ont tous deux pris pour modèle l'attitude lagrangienne. Leur différence consiste dans l'intérêt prééminent de Bolzano pour la méthode ; les efforts de Cauchy sont dirigés vers les résultats. Bolzano quant à lui, dans son texte *Rein analytischer Beweis des Lehrsatzes* (1817), va plus loin que Cauchy en abandonnant les « concepts inutiles » (les *infinitésimaux*, par ex.). En même temps que le livre de Lagrange paraît celui de Lazare Carnot *Réflexions sur la métaphysique du calcul infinitésimal* (1797). Ces deux ouvrages sont les premiers à poser ouvertement le problème des fondements mathématiques. Le thème du *Beweis* est la fonction continue à laquelle il donne une définition élégante. Le fait notable est ici pour nous que Bolzano se détourne de la saisie intuitive de la fonction. La coïncidence entre cette démarche et ce que l'on appelle « le programme de Cauchy », c'est-à-dire la proposition d'une Analyse algébrique reposant entièrement sur les concepts de base, est si frappante que certains auteurs avancent que Cauchy ne pouvait ignorer les travaux bolzaniens. Voir [Grattan-Guinness 1970: 51] ; [Hourya 1973]

Deuxièmement, chaque étape des théorèmes doit être démontrée ; cette démonstration doit découler de théorèmes eux-mêmes démontrées ;

Bolzano va encore plus loin : selon lui, nous devons avoir un fondement (*Grund*) sur lequel reposent nos raisonnements et conclusions. Cette exigence est formulée avec clarté aux § 512-537 de la *Wissenschaftslehre* (1820-1830) consacrés à la théorie de la démonstration. En particulier au § 525 *Erklärung des objectiven Grundes der Wahrheit,* il distingue, en suivant Aristote et en anticipant Frege, les démonstrations fondatrices (*Begründungen*) et celles qui établissent des certitudes (*Gewißmachungen*). Aucune objectivité scientifique, aucune science conceptuelle (*Begriffswissenschaft*) n'est pensable sans ces premières ;

Troisièmement, les définitions données et les théorèmes démontrés doivent avoir un caractère général afin de couvrir tous les résultats valables appartenant à un certain thème (ce point est particulièrement cher à Lagrange qui, tel Hegel, voyait dans n'importe quel résultat le moment d'un principe plus général).

Cauchy est sans aucun doute une figure-clé dans ce mouvement destiné à faire disparaître toutes les inexactitudes [Cauchy 1821 : ii-iii]. Par exemple, sa définition du concept de la limite est considérée comme révolutionnaire alors que la version classique, donnée dans l'*Encyclopédie* par d'Alembert et de la Chapelle, péchait justement par inexactitude [La Chapelle, d'Alembert : 1751-1772][1]. C'est ben toujours une question d'exactitude : si celle-ci ne vient pas de la mathématique, il faut donc la chercher ailleurs. Dans ces circonstances s'ébauche une réflexion inconnue auparavant, où la plus rigoureuse des sciences commence à s'interroger sur elle-même en se mettant en quête de la logique qui la gouverne. Le paradoxe tient à ce que cette logique n'existe pas avant qu'on ne la cherche ni ne la trouve ; autrement dit, cette logique est une construction des mathématiciens eux-mêmes.

Sitôt une évidence perdue, une autre se fait jour : chaque opération et chaque objet doivent avoir un sens logique, une valeur irréductible et irremplaçable révélés lors du processus de la construction d'une théorie (démonstration, etc.). Bolzano voit dans les mathématiques la science qui nous donne les lois (formes) universelles d'après lesquelles les choses doivent se disposer dans leur être[2]. Tel est aussi le point de vue exprimé par son contemporain Jozef Hoëné-Wroński (1776-1853), qui, dans son ouvrage consacré à la théorie de l'algorithme[3], trouve une logique même dans les nombres imaginaires comme $\sqrt{-1}$ (« monstres mentaux » selon Leibniz) : « <…> les nombres dits imaginaires sont éminemment logiques et, par conséquent, très conformes aux lois du savoir ; et cela, parce qu'ils émanent, et en toute pureté, de la faculté même qui donne des lois à l'intelligence humaine. De là vient la possibilité d'employer ces nombres sans aucune contradiction logique, dans toutes les opérations algorithmiques ; de les traiter comme des êtres privilégiés dans le domaine de notre savoir, et d'en déduire des résultats rigoureusement

---

[1] Pour la définition de Cauchy cf. [Cauchy 1821 : 19]

[2] Bolzano développe sa doctrine logique (métamathématique) dans *Contributions à une présentation plus solide des mathématiques* (1810) ; cf. [Cavaillès 1981 : 46] ; de même [Benoist 2002/3 : 289 et passim]

[3] P. Pragacz voit en Hoëné-Wroński « le pionnier de la pensée algorithmique en mathématiques » ; cf. [Pragacz 2007: 70]

conformes à la raison » [Hoëné-Wroński 1811 : 167] [1]. Il s'ensuit que dans tous nos emplois des nombres, il faut tâcher d'en extraire une pensée construite selon des lois analogues à celles de la nature. Le savoir mathématique est bon et utile s'il nous laisse voir son fondement logique, ses structures nues correspondantes à la pensée comme telle ; plus encore que Bolzano, Hoëné-Wroński insiste sur le privilège du savoir transcendantal d'accéder aux objets dans leur état le plus purifié et objectif. De plus, ces lois logiques, auxquelles nous remontons par le savoir mathématique, aiguiseront notre intelligence de sorte à la transformer en système de procédés logiques dont la subjectivité ne sera plus que temporaire : le sujet concret observe alors des lois qui ne sont pas le fruit de son imagination.

Dans les méandres d'une science véritable seulement se situent les principes absolus à partir desquels nous construirons un savoir total ; pour ce faire, il faut exploiter au maximum la force résidant dans la nature des mathématiques. La méthode doit être adaptée à son objet et si cet objet est la science fondamentale elle-même, pure et rigoureuse, alors il est nécessaire de choisir une méthode susceptible de conduire les mathématiques à prendre connaissance de leur caractère unique et de leur infaillibilité intrinsèque[2]. On sait qu'au cours du XIX^ème^ siècle ces intuitions sont confirmées dans bien des théories mathématiques (par exemple dans la théorie des ensembles, dont nous verrons les implications sur le parcours intellectuel de Husserl). Pour que la mathématique puisse s'élever vers le savoir absolu, elle doit non seulement définir et confirmer l'existence de ce dernier mais aussi en énoncer les principes. L'absolu n'est pas qu'un but ; il est avant tout pensée de la science pure manifestée en termes logiques. Si, par définition, l'absolu est dénué de contradictions et de défauts, il est logiquement parfait : la conscience qui saisit l'absolu est la même que celle qui établit des vérités logiques existant sans aucun rapport avec le monde. La connaissance de ces vérités est une forme suprême de la conscience.

Ce tournant vers la logique, qui marque les recherches de grands esprits du XIX^ème^ siècle, s'explique non pas par la position favorisée de la logique, mais plutôt par la quête d'une science universelle, fonctionnant à tous les niveaux : mathématique, logique et épistémologique en général[3]. A travers cette

---

[1] Il est curieux de constater que parmi les manuscrits de Bolzano, transmis par son disciple Robert von Zimmermann (1824-1898) à l'Académie royale de Vienne, il y ait des ébauches d'un système de mathématiques fondé sur la « nouvelle philosophie ». Il s'agit sans doute d'une philosophie qui s'appuie sur des principes autres que ceux de Kant. Cf. [Loužil 1978]. Sur Zimmermann et le cercle bolzanien voir [Winter 1933 : 127 et passim]

[2] En 1824 à Paris, il publie un ouvrage intitulé *Canons de logarithmes* (éd. Didot), où nous trouvons l'application concrète de ses méditations. Sans entrer dans des détails qui dépassent notre compétence, notons simplement qu'un des résultats les plus intéressants de cet auteur est la méthode générale de l'extension d'une fonction *f(x)* en séries infinies d'où vient sa « Loi suprême » qui sert de base au calcul différentiel (les séries wrońskiennes sont analysées avec finesse par [Deleuze 1968]. N'oublions pas qu'à la même époque, Bolzano élabore sa version du calcul infinitésimal suscitée par ses réflexions sur la géométrie euclidienne et par le prolongement des idées topologiques de Leibniz.

Quant aux *Canons*, ils s'inscrivent sur le plan philosophique dans le programme kantien, les algorithmes jouant le rôle d'un savoir transcendantal ; de même [Bushaw 1983 : 91-97]

[3] L'idée de Leibniz de créer la *lingua universalis* est perceptible chez Bolzano, Hoëné-Wroński et beaucoup d'autres. Selon Louis Couturat, en mai 1676 il médite sur le langage universel et le voit comme un calcul, comme une algèbre de la pensée ; cf. [Couturat 1901 : chap. III]. Georges Friedmann, ayant étudié cette question en détail, déplace l'invention du calcul différentiel vers les derniers mois de 1675 ; cf. [Friedmann 1962 : 99]

quête se précise l'idéal moderne de la science ; à celui des classiques (Kepler et Galilée) s'ajoute la croyance en les capacités particulières de la logique et en son statut d'ultime fondement. Non seulement la logique sert à mettre en ordre les pensées, mais aussi à retrouver leur origine. Ce faisant, la logique dépasse les mathématiques sur le plan réflexif en créant la possibilité d'examiner la « pensée exacte » du point de vue épistémologique ; bref, la logique devient le tréfonds de la connaissance.

Bolzano cherche l'axiomatique sur laquelle il convient de s'appuyer pour construire l'édifice de la science[1] ; Hoëné-Wroński cherche cet appui parmi les lois logiques les plus générales dans lesquelles se retrouve l'esprit humain[2]. Reste à saisir un contenu, un sens positif tel que ses manifestations ne modifieront jamais son statut : le sens idéal. Autrement dit, la crise peut être dépassée grâce à l'exactitude de la science, où la particularité (axiome, proposition, etc.) est en même temps un principe plus général et tend infiniment vers la vérité. C'est pourquoi Bolzano appelle son entreprise *Wissenschaftslehre* (« la théorie du savoir »[3], première critique radicale de Kant) et Hoëné-Wroński consacre son ouvrage mathématique principal à la recherche des axiomes présentant « une vérité absolue qui se détache de toute matière et que la raison peut seule poser » [Gonseth 1945 : 33][4].

Maintenant le chemin est inverse : il faut aller non pas d'une science vers la vérité mais de la vérité vers la science capable de l'exprimer en langage adéquate. La logique, après avoir purifié ses concepts et précisé leur application, sera sans aucun doute le meilleur instrument à cet effet : « [elle], je pense, apporte à la véritable pensée autant que la grammaire à l'usage correct de la langue » [Bolzano 1985: 67], d'où la question : <…> quand la logique établie non seulement les lois qui s'appliquent aux vérités spéculatives (*gedachten*), les vraies pensées comme l'on dit, mais les vérités en général ? Si cela n'est pas juste pour des propositions spéculatives mais pour les propositions en soi, alors doivent-elles reposer sur la validité de la règle logique ? »[5]. Pour autant, elle

---

Au XX[ème] c'est Gottlob Frege qui reprend le plus fidèlement ce projet leibnizien dans son *Begriffsschrift* (1879 ; *L'idéographie*) dont la tâche était de conceptualiser les mathématiques par un langage logique parfait. Pour une discussion détaillée : [Moese 1965 ; Morscher 1972]

[1]Une attention particulière à ce problème est portée au § 223 de la *Wissenschaftslehre*, dans lequel Bolzano élabore le concept des vérités de fondement (*Grundwahrheiten*) qui résident au-delà de toute expérience. Il est impensable de mettre en doute le caractère fondamental de ces vérités car il conditionne ce doute-même.

[2] Voir [Braun 2006] En français il existe un article intéressant [Phili 1996]

[3] Autre traduction possible : « théorie de la vérité » ou même « ontologie de la vérité ». D'ailleurs, Bolzano sent qu'il parvient – au moins dans la première partie de son ouvrage (*La théorie des fondements*) – à une théorie ontologique de la vérité (ou des vérités) et n'en est pas satisfait, comme l'atteste sa lettre à Zimmerman où il préconise de considérer son travail en tant que théorie critique de la vérité. Cf. [Winter 1949: 90]

[4] Et ensuite : « l'axiome est ainsi une vérité de raison, non une vérité d'expérience » ; *ibid.*

[5] *Ibid.*, p. 78. A ce propos : l'évolution de l'attitude de Husserl vis-à-vis de ces propositions mérite de s'y arrêter. A l'époque de ses études auprès de Weierstraß et même à Halle, il les considère comme des essences mythiques suspendues entre l'être et le non-être ; cf. [Husserl 1975 : 33-38] La lecture de Hermann Lotze, et en particulier son idée de validité (*Geltung* ; cf. [Lotze 1843 : 92, 97, 141-143]), qui – sans l'ombre d'un doute – contribue considérablement à la compréhension husserlienne de Platon, le fait revenir sur sa position antérieure : « la théorie de Bolzano selon laquelle les propositions sont des objets néanmoins dépourvus d'existence est assez intelligible. Elles privilégient l'existence idéale ou la validité des universaux » ; [Husserl 1994 : 202] ; de même [Fisette 2009 : 281-284]

n'est pas la seule *science*-reine de la pensée ; elle doit jouer son rôle *dans* l'épistémologie sans jamais remplacer cette dernière. « La logique [pure] doit me guider vers une théorie de la science, c'est-à-dire donner un *modus operandi*, entre autres pour ordonnancer toute la sphère de la vérité dans une *théorie de la science* adéquate pour que chaque partie élaborée et exprimée soit à sa place »[1], elle n'est pas but mais moyen ou, si l'on ose dire, le passage au vrai. Curieusement, ces recherches du début du XIX[ème] siècle montrent une similitude avec la pensée scolastique (onto-théologique), elle-même bâtie sur des axiomes quasi-donnés[2] et ne cherchant qu'à prouver, à démontrer une vérité préexistante : Dieu est la cause première et nécessaire de tout « exister »[3]. La vérité, pour Bolzano comme plus tard pour Husserl – n'est-ce pas pour cette raison que ce dernier appelle le philosophe tchèque « l'un des plus grands maîtres-logiciens de tous les temps » ?[4] –, est incontestablement objective. Il faut avant tout distinguer la vérité en soi de la vérité reconnue (*Wahrheit an sich… erkannten Wahrheit*) [5] ; la première possède la particularité (*Beschaffenheit*) que le vrai et le faux ne sont pas discernés en elle, c'est une possibilité dont la réalisation dépend de la qualité des propositions (*Sätzen*).

Bolzano souligne l'erreur de l'approche psychologique consistant à les mélanger, ainsi que la nécessité de l'effacer pour construire une vraie science de la logique. En critiquant Kant avant tout pour son épistémologie anthropologique, Bolzano semble conserver l'idée kantienne d'une objectivité ineffable : la vérité en soi (au sens propre du terme) existe mais il est impossible de l'exprimer. Tout langage, y compris celui de la logique, reste indigent devant l'ipséité de la vérité. Ce que nous pensons et exprimons dans un système de propositions (*Sätzen*) est par définition hors-soi mais le langage ne nous donne qu'une objectivité partielle car il appartient au sujet. Voici la racine de la

---

[1] *Ibid*., p. 86. Notons ici qu'au § 59 des *Ideen I* où Husserl précise la phénoménologie de la conscience transcendantale, il répète l'idée de Bolzano : « La logique formelle et toute la mathesis en général peuvent être incluses dans l'ἐποχή qui procède précisément à l'exclusion » ; Hua XIX/1, p. 113.

[2] Nous pouvons sans doute comparer ces recherches du XIX[ème] siècle avec la crise ontologique de l'époque de saint Thomas d'Aquin où la puissance ontologique de Dieu – l'acte créateur – devait être expliquée, pour ne pas dire ajustée au modèle du savoir aristotélicien. Cette crise avait été provoquée par la rencontre d'Aristote avec la foi chrétienne : comment exprimer l'essence divine, son infinité – l'axiome des scolastiques – en termes logiques et éviter toute limitation de Dieu ? Voir, par ex. : [Tresmontant 1964]

[3] Au § 25 Bolzano s'en explique : « Aus der Allwissenheit Gottes folgt zwar, dass eine jede Wahrheit, sollte sie auch von keinem anderen Wesen gekannt, ja nur gedacht werden, doch ihm, dem Allwissen, bekannt sei, und in seinem Verstande fortwärend vorgestellt werde. », p. 138 ; cf. [Bolzano 1985: 138]

[4] Pour la relation de Husserl à Bolzano voir [Sebestik 2003]

[5] Cf. *Wissenschaftslehre*, § 26. La similitude des théories logiques de ces deux penseurs a été étudiée par [Fels 1927]. Fels relève à juste titre que Husserl est passé à ses débuts du coté de l'antipsychologisme bolzanien. En 1935 Husserl confesse à Andrew Osborn (1902-1997) qu'il a eu la chance de tomber dans une librairie sur les livres de Bolzano ; cf. [Schuhmann 1977: 463] ; de même [Osborn 1934]

Notons au passage : ce concept de « vérité en soi » et l'épistémologie bolzanienne en général trouvent pour beaucoup leur origine dans l'ouvrage de J.-H. Lambert intitulé *Neues Organon oder Gedanken über die Erforschung und Berechnung des Wahren* (1764 ; en particulier vol. I, chap. II consacré à l'aléthologie). Nous verrons plus loin comment Husserl transforme la « vérité en soi » en objet noématique par son inclusion dans la sphère de la connaissance. Ce concept connaît une autre destinée chez Chalva Noutsoubidzé (1888-1969) qui critique le subjectivisme transcendantal husserlien tout en soutenant « la vérité en soi » de Bolzano, qu'il comprend comme l'objectivité privée d'humain ; cf. [Noutsoubidzé 1926] critique bolzanienne de Kant qui prend le sujet en tant que fondement transcendantal, position qui ne peut se solder que par un échec. Le transcendantal au contenu anthropologique ne définit pas tant les frontières du savoir que les possibilités du sujet-même et il ne reste à Kant qu'à démontrer ce fait. Selon Bolzano, l'être humain connaît certaines vérités mais sa connaissance n'est jamais complète. Pire encore, par cette dernière, le sujet se trouve détaché des choses comme telles existant pour lui uniquement dans son système de propositions ; celles-ci ne se confondent pas avec les énoncés linguistiques (avec la langue, sphère pratique) mais indiquent plutôt la possibilité de l'expression. Celle-ci réside dans le « en soi » (*an sich*) qui se distingue du concept kantien : chez Kant, « la chose en soi » signifie une chose principalement inconnaissable, qui le demeure en limitant notre connaissance après tout acte intentionnel. Cette limite n'appartient pas à la chose en soi même. Elle est une caractéristique fondamentale de la conscience : nous ne pouvons connaître l'« en soi » de la chose non parce qu'il existe, mais parce que l'acte du connaître s'avère inadéquat par rapport à ce « soi ». Il s'agit, soulignons-le, non pas d'une limitation ontologique, comme dans le cas de l'essence aristotélicienne, mais du principe du savoir. Kant ne tâche pas de connaître les choses comme telles, son but consiste à prendre connaissance du sujet comme il *est* et comme il *pense*, c'est-à-dire tel qu'il ne quitte pas le champ empirique même lorsqu'il applique au monde ses facultés aprioriques. C'est pourquoi Kant établit les limites du sujet anthropologique sans mettre jamais à l'épreuve sa nature ; il reste toujours avec ce Je empirique que Husserl, après Bolzano, va surmonter par l'έποχή et par la création du sujet transcendantal.

## CHAPITRE 2. L' « ANTI-KANT »[1]

Bolzano construit son propre idéal épistémologique en radicalisant, dans une certaine mesure, la critique kantienne : le sujet du savoir doit être identique à celui de l'existence. En d'autres termes, l'existence n'est plus une donnée *a priori*, elle est considérée comme l'épiphénomène de la connaissance dirigée vers l'appréhension de l'objectivité des vérités. Comprendre cette objectivité naturelle est le but ultime de toute théorie du savoir. Le sujet clôturé dans son espace anthropologique étant inapte à effectuer cette démarche, il faut donc recourir à un autre sujet capable de dépasser les restrictions kantiennes. Le sujet logique, qui n'envisage pas son existence comme condition préalable pour appréhender *les* objectivités, est la seule solution bolzanienne. Concrètement, cela signifie un refus du sujet en chair et en os qui déduit la connaissance objective de la somme de ses expériences psychiques et somatiques, qui perçoit sa vie quotidienne comme lui étant donnée à travers ses facultés perceptives. Le nouveau sujet épistémologique, s'il apparaît un jour, doit percevoir l'objectivité *en soi,* existant parallèlement à sa perception et comprendre que cette-dernière ne fait que porter atteinte aux qualités de toute objectivité.

[1] A notre connaissance Bolzano lui-même n'a jamais rédigé le texte sous un tel titre, mais il avait représenté le style de cette pensée anti-kantienne (voir l'ouvrage du chanoine tchèque Franz Příhonský (1788-1859), disciple du premier cercle de Bolzano, intitulé N*euer Anti-Kant* paru en 1857 ; cf. [Příhonský 2006]. Les titres apparaissent plus tard chez les auteurs qui, à la différence de Bolzano, critiquent le coté spéculatif de Kant et son chose en soi. Cf. par ex. : [Büchner 1884 ; Bolliger 1882]; Büchner, quant à lui, va jusqu'à faire appel aux philosophes de « s'opposer énergiquement à l'imposture kantienne (*Kant-Schwindel*) ; [Büchner 1884 : 315].

Ainsi s'effectue la réduction des actes subjectifs qui entrent peu à peu dans le système général de la connaissance composé de vérités d'ordres différents. Leur présence dans le sujet atrophie la part anthropologique de telle sorte que les limitations semblant imposées par la nature s'en trouvent fausses. Le monde empirique, avec toutes ses manifestations connues du sujet pensant, n'est qu'une plateforme pour accéder à l'objectivité en soi, distincte d'une quelconque expérience personnelle. Le sujet idéal chez Bolzano doit migrer vers la connaissance afin de se trouver maximalement indépendant des affects de l'existence ; autrement dit, il transforme sa perception basée sur ses expériences quotidiennes en conscience ; la conscience rassemble les expériences quotidiennes dans une unité dont l'analyse qui suit ne tend qu'à montrer le caractère superficiel de la perception empirique. A la différence de Kant, Bolzano insiste sur l'existence absolument objective des vérités synthétiques. Celles-ci ne résultent point du cogito : c'est le *cogito* qui les acquiert par ses réflexions. La tâche la plus ardue consiste à dénicher les vérités basiques de l'abondante végétation des choses et des impressions venant à nous à chaque instant de la vie [Bolzano 1985: 355]. Par sa nature, la logique rappelle la physique : elle est composée d'éléments primitifs entrant dans tout acte de pensée et dans tout jugement. Le nombre de leurs variations est infini ou presque. Nous prenons ces variations pour des vérités au sens propre du terme, tandis qu'elles ne font dans la plupart des cas que refléter nos propres égarements et erreurs.

La nouvelle « théorie du savoir » doit être la théorie des vérités fondamentales – des primitifs logiques – que Kant a senties mais en en donnant une description erronée, en les coupant de façon trop formelle en deux classes : analytique et synthétique[1]. Le critère kantien s'est avéré trop abstrait et par conséquent extérieur aux jugements auxquels il s'applique ; il passe dans la sphère du sujet connaissant mais non dans les vérités en tant que telles. Kant lui-même ne marque pas ce fossé ou ne le veut pas reconnaître en se bornant à postuler sa mystique *chose en soi*. De plus, Bolzano insiste sur la distinction non pas entre deux classes de jugement mais entre deux types de vérité qui conditionnent la logique de base : les énoncés comme « 7 + 5 = 12 » ou « la

[1] Dans le tome III de la *Wissenschaftslehre* (p. 240), on trouve l'analyse détaillée de cette distinction kantienne. Jean D'Alembert distingue deux types de démonstration que l'on peut utiliser dans la géométrie : directe et indirecte. Il est curieux que la première ressemble à l'idée de l'analycité chez Kant. Selon D'Alembert, la démonstration directe se déduit immédiatement du concept de la chose à laquelle nous voulons attribuer une propriété. Cette attitude, esquissée déjà dans le *Tractatus de Geometria elementari* (1716/17) et dans le *Traité de Géometrie* du leibnizien Pierre Varignon (1654-1722), a été reprise ensuite d'une manière ou d'une autre par des savants de l'époque, parmi lesquels Jean-Paul de Gua de Malves, auteur d'*Usage de l'analyse de Descartes* (1740) dont Kant pouvait avoir connaissance. Pour plus d'information cf. [Hankins 1970]

Pour la critique contemporaine de Kant, fondée largement sur l'approche bolzanienne, cf. [Quine 1980: 20]; de même [Bar-Hillel 1970]

Nous devons à Alfred Tarski une autre version du concept de l'analyticité qu'il déduit de son exemple bien connu : « la neige est blanche ». Ce jugement est vrai si la neige est blanche en réalité (au moment de l'énoncé). L'introduction de la sémantique contextuelle – ou vérité *sémantique* – modifie notablement l'idée kantienne ; cf. [Tarski 1944: 342] somme des angles du triangle est égale à deux angles droits » sont analytique car les propriétés (*Beschaffenheit*) exposées dans chaque de ces énoncé sont les vérités innées. C'est pourquoi l'analyticité bolzanienne, à la différence de Kant, se dégage de l'ontologie du jugement, elle montre les idéalités dans leur essence composée des telles ou telles qualités. En revanche, les mêmes jugements dans la théorie kantienne sont synthétiques exposant l'intuition subjective.

Le concept de l'analyticité en tant que propriété innée des éléments (incluse dans les opérations avec ces éléments[1]) est épitomé ainsi : « il est claire que la proposition '*A* possède une propriété *a*' » est analytique » [Bolzano, *Wissenschaftslehre* : § 304], si – répétons-le – cette propriété appartient à l'essence de la chose. N'oublions que le concept bolzanien du contenu est un ensemble de ses parties mais non pas des propriétés. Michael Dummett a raison de dire que la proposition analytique chez Bolzano trahit « une idée complexe contenant ses constituants… » [Dummett 1991: 29 ; cf. *Wissenschaftslehre* : § 148]. Pour Kant l'énoncé « tous les corps sont étendus » est analytique car « l'étendu » est une propriété innée de tout corps ; pour Bolzano cet énoncé peut se lire aussi bien analytique que synthétique. Synthétique, parce que notre connaissance des corps étendus est fondée sur l'expérience apostériorique. Si cette dernière fait partie de des jugements analytiques, alors il est impossible de les prendre comme le sous-œuvre de la science[2].

Selon *Wissenschaftslehre*, pour achever son idéal, le sujet tend à modifier sa subjectivité de telle manière pour la faire l'objet de la connaissance. L'idéal se réalise lorsque le sujet redécouvre son Je pensant comme l'ensemble des opérations logiques créant le champ mental du sujet. En bref, Bolzano cherche les sorties de « l'impasse kantienne » en changeant non seulement les principes critiques de Kant mais avant tout sa tâche en enlevant les limites du savoir [Bolzano, *Wissenschaftslehre*, Bd. III : 235 et passim][3]. Cela explique la position anti-aristotélicienne du bohémien envers l'infinité[4] ; pour ce dernier l'infini est par nature actuelle et elle ne se réalise qu'en tant que telle. Il n'existe pas des choses inconnaissables car leur inconnaissabilité est déjà une caractéristique de notre connaissance présente (et de notre jugement) de ces choses[5].

---

[1] L'étude relativement récente consacrée au désaccord entre Kant et Bolzano : [Berg 1999]. En effet Berg donne un nouvel éclairage à l'analyse livrée par [Příhonský 2006]. Il s'agit du premier ouvrage consacré à ce sujet. Voir aussi [Laz 1993]

[2] Toujours dans le tome III de *Wissenschaftslehre* (p. 101 et passim) Bolzano conteste la division kantienne du savoir en a priori et a posteriori.

[3] De même [Dubislav 1929]

[4] L'actualisation bolzanienne de l'infini sera poursuivie dans les travaux de Cantor dès l'année 1872 où il publie son article sur les points particuliers ; plus tard la théorie des ensembles de Cantor, basée sur des présuppositions inconcevables (échappant à l'intuition triviale), sera développée à son tour par le logicien Abraham Robinson (1918-1974), inventeur de l'analyse non standard en mathématiques renvoyant aux idées aussi bien de Leibniz et Newton que de Bolzano et de Hoëné-Wroński.

Voici un détail intrigant : en réintroduisant les infinitésimaux dans l'Analyse, en les montrant comme objets d'une analyse rigoureuse, Robinson emploie le concept, plutôt « l'axiome de l'idéalité » en l'attribuant le sens husserlien (cf. *Ideen I*).

[5] Bolzano mène sa bataille contre l'agnosticisme en particulier dans le troisième tome de son ouvrage où il précise sa critique de Kant concernant la limitation du savoir. Il faut dire que malgré une certaine sophistique de son argument, l'idée même était bouleversante compte tenu de l'époque qui a méprisé toute sorte d'optimisme gnoséologique. Pour plus de détails voir [Winter 1949]

L'idéal de la science réside donc dans son caractère universel qui ne se change pas ni dépend des opérations particulières. Si la mathématique est née comme un outil destiné à aider les hommes dans leur vie quotidienne, alors elle doit parvenir à la logique pure de laquelle la vérité se déduit de la manière la plus naturelle. Pour Bolzano l'appareil logique est simple, sa charpente ne comprend que trois éléments : sujet, copule et prédicat avec lesquels on peut construire des propositions de base [cf. Mugnai 1992] ; « Socrate est mort » exprime une idée irréductible et à l'économie parfaite. La vérité logique consiste en ce qu'elle ne peut être exprimée de façon plus économique, son objet peut cependant varier selon ce qu'elle désigne comme vrai. Autre aspect capital : la vérité dans la logique ne dépend point de sa particularité ni de son abstractivité (ou, en termes scolastiques : de *species specialissima* à *genera generalissima*) : « Socrate est mort », proposition qui concerne un personnage concret, est aussi vraie que « Dieu est bon ». Le sens logique est indépendant du statut existentiel des objets, qui peut être très précis ou bien indéfini [1]. De surcroît, Bolzano récuse les vérités ontologiques vivant par elles-mêmes hors du champ discursif, leur caractère *sui generis* n'étant pour lui qu'une illusion. Les vérités logiques surgissent au moment où nous les énonçons. « Socrate est mort » signifie la vérité non parce que cette phrase informe sur un fait historique que personne ne conteste mais parce qu'elle *signifie* ce fait en lui donnant un sens logique, et non historique. C'est pourquoi, pour fonder les mathématiques, il faut à la fois isoler les principes et décrire les procédés d'entendement logique, entreprise pour laquelle les outils logiques justement manquaient à Bolzano (à son époque la logique est encore *ancilla mathematicae*). Quel que soit le rôle de Bolzano, adversaire de Kant confronté au même problème, la critique de la raison mathématique demandait à être menée à son terme.

Il faut utiliser la logique en tant que rasoir d'Ockham, notamment scinder l'objet mathématique en deux parties : le *contenu* (nombre, point, etc.) qui remplit la conscience du mathématicien et le *sens* lui-même, c'est-à-dire une idéalité pure demeurant après toutes les abstractions d'ordre supérieur. Un néophyte ou « algébriste artisan », selon les termes de Husserl, ne discernerait pas la différence subtile qui existe entre eux. Or il s'agit en réalité d'un paramètre capital. Le contenu mathématique, même s'il traverse l'histoire sans changer en apparence, est toujours subjectif. Par exemple, non seulement le nombre 3 peut désigner des choses très différentes (trois pommes, trois mois, la Trinité chrétienne ou trois frères), mais il s'emploie dans chaque cas de manière particulière, se réalisant simultanément en général et comme individuation. Le nombre est un dessin psychologique qui ne dépasse pas le champ cognitif du sujet et par lequel nous faisons notre abstraction des objets infinis. Quand j'applique « 3 » aux pommes sur la table ce n'est pas le même « 3 » que si je l'applique aux trois couleurs du drapeau français. Dans ces deux cas se distingue non pas l'enveloppe symbolique mais la perception du nombre. En revanche, il est impossible de percevoir le sens du nombre qui reste dans tous les cas une essence intelligible ; autrement dit, le sens est ce que l'on ne peut abstraire. L'acte d'abstraire s'applique donc aux couches de perception, au

[1] Ici Bolzano suit dans une certaine mesure l'idée d'Alain de Lille qui ne reconnaît rien comme existant vraiment (*vere existit*) sauf Dieu dont l'existence est absolue, c'est-à-dire au-delà de l'idée même de l'existence. Il est notable que ce scolastique distingue clairement l'existence – catégorie indéfinie pour lui – et les pures opérations logiques qu'il demande de fonder sur des axiomes, à l'instar d'Euclide. Cf. [Migne 1844-1854 : 2]

champ subjectif retracé à mesure que l'abstraction progresse. Ainsi se constitue l'objet idéal qui n'est ni analytique ni synthétique (selon Kant) mais phénoménologique, comme le montre Husserl. C'est pourquoi quand nous parlons des objets idéaux, il ne faut pas confondre contenu et sens, qui ne se distingue pas dans les objets empiriques. Les mathématiques donnent les meilleurs exemples d'idéalités existant sans théorie qui pourrait les décrire.

L'idée que l'objet mathématique est composé de ces deux parties résulte déjà d'une abstraction – chirurgie logique, dira-t-on – permettant de le voir hors de son recouvrement par le sujet. Dès lors, qui le voit ? Réponse : celui qui ne considère pas cet objet comme *mathématique* ni même comme objet existant en soi. La chirurgie logique permet de voir dans pareil objet la structure de la connaissance, ce qui justifie la réduction des mathématiques à une série d'opérations logiques. Le but de cette réduction, appelons-le « l'établissement des vérités axiomatiques », est de libérer la pensée des propositions inutiles et de dépasser la mathématique naturelle fondée sur notre perception immédiate. Rien n'est plus naturel que la série des nombres $p_1$, $p_2$, $p_3$… $p_n$ à laquelle on peut toujours ajouter 1 ; rien n'est plus simple que les quatre opérations de base (addition, soustraction, multiplication et division) que nous apprenons à l'école primaire. La propriété fondamentale de ces nombres et opérations consiste en ce qu'ils ne se déduisent de nulle part et sont des flux hylétiques de l'intuition. En allant plus loin, on peut même dire que la série mathématique ressemble à l'idée de Dieu : toutes deux sont essentiellement intuitives et inaliénables à l'esprit humain. Si, à titre d'exercice imaginaire, nous essayions d'enlever ces deux idées (ou une seule), nous nous heurterions à un obstacle infranchissable : la possibilité-même d'exécuter une telle opération est conditionnée par ces idées. En d'autres termes, il s'agit de concepts seuls à créer la possibilité de les repenser. Il est donc des structures qui ne peuvent être exclues de la pensée ; c'est à partir d'elles que s'élabore une autoanalyse des mathématiques révélant l'objectivité qui caractérise toute la sphère de l'algèbre et de la géométrie. Précisons : les objets abstraits et leurs connexions présentent une/des hiérarchie(s) engagée(s) dans un processus de construction permanente. Sur le plan logique, les éléments de cette hiérarchie sont vides, dépourvus de sens propre ; les nombres naturels, séries, etc., ou toutes autres abstractions algébriques n'acquièrent leur sens qu'en étant en relation avec d'autres entités de même ordre. Cette emphase sur le caractère relationnel de la pensée algébrique remonte au siècle de Descartes et se distingue par deux traits : la nécessité de l'inventaire symbolique et l'exemption de toute ontologie ; celui-ci aura dans l'histoire ultérieure une importance capitale car il soulève le problème du statut de l'objet mathématique. La valeur de la preuve ontologique a été considérée irréfutable depuis l'Antiquité, elle semblait protéger chaque nouveau concept de possibles contradictions, *a fortiori* s'il s'agissait d'un concept surgi de l'intuition. Tout objet dénué de cette preuve se trouvait hors du traitement théorique. Les premières fissures de la crise des fondements apparaissent donc au milieu du XVII^ème^ siècle.

Ainsi, l'objectivité mathématique est un produit des relations ; mieux encore, elle est *cartésienne*, n'existant que mentalement. D'après les experts, la hiérarchie des relations maintient la cohérence du monde des idées. L'algébriste ne sait distinguer les idées et ce qu'elles désignent : à la différence de l'essence aristotélicienne, les entités mathématiques sont non seulement réelles *in potentia* mais elles forment leur réalité idéelle.

Que fait la logique ? Elle étudie les sens relationnels des concepts objectifs en ôtant toute la réserve psychologique qu'ils ont accumulée. C'est pourquoi l'abstraction logique diffère grandement de l'abstraction mathématique : la première reçoit son contenu de la compatibilité des éléments, la seconde de leur regroupement ; le savoir logique délimite le concept en le réduisant à sa fonction, le concept se confond alors avec les limites de son fonctionnement. La question de sa réalité n'a pas de sens. Les mathématiques effectuent la démarche inverse : le concept est reconnu comme existant avant qu'il ne soit mis en branle, ce qui témoigne du caractère non-métaphysique de cette science (« non-métaphysique » au sens cartésien, ignorant la question des conditions du savoir mathématique). Pour le mathématicien, la pensée ne précède point l'entendement spécifique concernant son champ de recherches, le *cogito ergo…* est un énoncé étranger au travail des arithméticiens, des géomètres, etc. car le doute cartésien ne saurait concerner que les choses dont l'existence n'est pas absolument révélée. Le problème de l'existence en mathématiques est remplacé par une existence idéale qui ne se déduit pas de l'acte de penser ; penser équivaut ici à manipuler des nombres. Métaphysiquement parlant, l'algébriste se laisse conduire dans des sentiers battus, préalablement définis par l'emploi de symboles respectant certaines règles. En bref, la différence entre le mathématicien professionnel et le philosophe cartésien tient à ce que le premier ne met jamais en doute le fondement ultime de son raisonnement.

## CHAPITRE 3. ONTOLOGIE OU NON-ONTOLOGIE

La crise du XIX$^{ème}$ siècle est celle du nombre [cf. Snapper 1979 : 207-216]. Dans les années quarante, déjà avec les travaux de Hermann Grassmann, esprit encyclopédique dont les mérites n'ont été vraiment reconnus qu'après sa mort, il devient évident que le nombre arithmétique ne peut se passer de conceptualisation. Il a bien un contenu mais est dépourvu d'un *sens* qui pourrait être abstrait de toutes les particularités et pris comme point de départ d'une critique quasi-kantienne. Dans quelle mesure peut-on parler de nombre pur, propriété de la conscience et non plus (seulement) de la pratique du calcul ? Ou bien, le nombre peut-il être le fondement du savoir *a priori*, d'autant plus que la conceptualisation du nombre par les moyens de la logique aboutit au contenu fonctionnel du nombre ? D'où la question : est-il suffisant pour attaquer la critique des fondations ? En 1861, Grassmann essaie de formuler sa théorie de l'axiomatisation basée sur les principes de l'induction ; il présente ses résultats dans *Lehrbuch der Mathematik für höhere Lehranstalten* (Bd. 1). Il s'agit de la tentative de réduire certains concepts de la géométrie à des symboles algébriques dont le sens est défini par leur relation réciproque (Bolzano, rappelons-le, fut un des premiers partisans de cette réduction [cf. Bolzano, *Rein analytischer Beweis...* II]). Par exemple, avec les éléments $e_1$, $e_2$, $e_3$, l'on peut engendrer l'espace linéaire de notre choix (en dépassant les restrictions des trois dimensions) ; chaque semblable produit est composé d'un nombre limité d'éléments ou, plus exactement, de leur ordre. Chaque changement du produit résulte d'un changement de l'ordre : concrètement, cela signifie que tout résultat dépend des règles par lesquelles on établit/construit un objet. L'idée de Grassmann ne manque pas d'attraits car non seulement elle propose une axiomatique, aussi imparfaite soit-elle, mais elle libère la mathématique de l'époque des propositions dont le contenu empêche de saisir son sens.

Dès la seconde moitié du XIX$^{ème}$ siècle la principale tâche des algébristes consiste à trouver un nouveau point d'appui – à l'instar du levier d'Archimède – aboutissant plus tard à une rénovation considérable de la science. Les mathématiciens n'ont pas inventé un concept transcendantal du nombre, si on entend par nombre une entité idéale qui peut n'avoir aucun lien avec le monde empirique ; en revanche, l'idée de bâtir toute la complexité du calcul sur quelques opérations produisant le sens sans aucune dépendance avec les objets n'est pas tant une solution logique aux problèmes qu'une charnière gnoséologique. En termes scolastiques, l'existence des vérités mathématiques s'est avérée mois importante que leur essence. La science exacte et la science de l'esprit, la métaphysique, convergent sur ce point dans leur commune nécessité de reconsidérer leurs fondements. Cependant, penser l'*essence* des objets idéaux ou a fortiori des opérations mathématiques, objectivité absolue par excellence, nous mène à entrevoir une autre perspective. Le sujet pensant doit s'assurer de la nature pensée de l'essence ; selon Aristote, l'essence n'a rien, excepté son intelligible nature, qui ne se laisse exister que dans l'entendement subjectif. Descartes, par exemple, n'y voit aucun problème et de préciser qu'il ne faut pas séparer l'essence et l'existence dans le cas des objets dépourvus d'ontologie, comme chez Dieu. Quelle que soit la parenté entre métaphysique et mathématiques, la science « la plus exacte » repose sur les préceptes de la première quand elle permet au sujet de déduire les vérités objectives de sa pensée sur le nombre. Soit : la méthode axiomatique actualise même l'objet idéal de telle manière que l'essence fusionne avec le contenu. Cette nouvelle approche pose une nouvelle question (aussi bien arithmétique que philosophique) : le concept du nombre peut-il exister sans la chose qu'il désigne ? Comment donc *penser* le nombre, entité pure de la conscience ? Cela devient le thème central de *Philosophie de l'Arithmétique* (1891) sur laquelle nous reviendrons plus tard ; même si à cette époque Husserl manifeste des affinités psychologistes, son but diffère aussi bien des adeptes de l'école de Franz Brentano que des algébristes (comme Ernst Schröder) : pour les premiers, les concepts mathématiques restent descriptifs, définissant les processus psychiques et servant à saisir le contenu de ce qu'ils décrivent ; pour les seconds, les concepts doivent aider à créer une science où le sens et les tâches seraient bien définis. En créant sa théorie algébrique – qui prend en considération des idées de Ch. von Sigwart et de H. Grassmann –, E. Schröder remarque que « la logique-même est une logique qui établit des *contenus* conceptuels » [Schröder 1966 : 89, italique originel][1]. D'un coté, la logique est composée d'un système d'individuations (*Einzelding*) distinguant les particularités et par cette distinction établissant les concepts (*Begriff*) – à chaque chose sa définition – qui permet à la pensée de voir chacune séparément et en

[1] Cet ouvrage a été publié en trois volumes entre 1890 et 1905, à compte d'auteur. Le volume III compte deux parties, la deuxième publiée à titre posthume, et éditée par Eugen Müller (1865-1932), disciple du mathématicien Jakob Lüroth (1844-1910). Les *Vorlesungen* constituaient une somme complète sur l'état de la logique « algébrique » (ou « symbolique » en termes d'aujourd'hui) à la fin du XIX$^{e}$ siècle.

Ici Schröder rejoint Sigwart sur sa définition du rôle de la logique qui consiste à mettre en ordre les catégories de la pensée formées au cours de la vie psychique ou, plus exactement, de la variété des représentations.

totalité et par là de se voir elle-même ; d'un autre coté, Schröder semble développer l'aspect relationnel de la logique algébrique déjà en germe chez Grassmann. Il s'agit des relations ordonnées dans le système selon certains modes où chaque élément présente justement une relation[1]. Concrètement, les relations sont prises comme base primaire pour construire des théories de toute complexité. Un détail important : Schröder insiste sur le fait que le concept/la définition logique doit avoir – se fonder sur – un contenu strict (*engeren Inhalt*) auquel « s'amarrent les sens étendus composant à leur tour les concepts d'ordres supérieurs » [Schröder 1966 : 91]. En outre, une vision philosophique n'est pas étrangère à Schröder quand (à la suite de Descartes et Leibniz) il médite sur le langage universel de la science qui pourrait décrire de façon adéquate tous les actes fondamentaux de la pensée. L'idée s'appelle « le *langage caractéristique* des concepts ou "idéographie" (*Begriffsschrift*), dont la tâche est de décrire non pas le contenu des représentations (*Inhalt der Vorstellungen*) mais… les idéalités pensées » [Schröder 1966 : 93-95]. Le *Begriffsschrift* est une allusion directe à l'ouvrage de Frege de 1879 intitulé *Idéographie* [2] qui retient l'attention de ses critiques (B. Russell, K. Gödel, A. Tarski)[3]. Ce livre marque un tournant dans la pensée de Frege[4] : ayant toujours considéré l'arithmétique en tant que domaine de la logique (avis que partagent Schröder et le jeune Husserl), il ne lui accorde pas – à la différence de la géométrie – le statut de savoir intuitif non soumis à l'axiomatisation logique ; après la parution du *Begriffsschrift*, le

[1] Un des critiques les plus intéressants de Schröder en même temps que son chantre est Platon Poretsky (1846-1907), professeur de logique à l'Université de Kazan. Dans ses travaux Poretsky (comme Schröder) s'appuie sur l'algèbre de Georges Boole (1815- 1864) et élabore la logique relationnelle ; cf. [Royce 1913; Kneale 1948].

Quant à Poretsky, sa méthode consiste à déduire de nouvelles relations à partir des relations données, nées de nouvelles combinaisons. Pour Couturat il s'agit de « la méthode des causes et des conséquences » qu'il considère comme un perfectionnement des résultats précédents, notamment de Boole, Schröder et de John Venn (1834-1923) qui d'ailleurs géométrise la logique sur le modèle leibnizien ; cf. [Venn 1980 : 76] Selon Poretsky, toute équation logique est – ou doit être – composée d'un nombre d'éléments constitutifs irréductibles. Pour passer d'une équation à sa conséquence, il suffit d'exclure les éléments qui composent cette équation même. Les idées de Poretsky précèdent l'axiomatique de Hilbert ; l'axiomatique est complète quand le contenu de chaque proposition représente une vérité prouvée. Voir [Poretsky 1899 ; 1900]

[2] Le titre original complet est : *Begriffsschrift, eine der arithmetischen nachgebildete Formelsprache des reinen Denkens* (Halle a/S: Verlag von Louis Nebert). Comme Leibniz avant lui et Schröder après, Frege réfléchissait à un langage de la pensée pure (une arithmétique idéale) qui serait aussi le fondement des mathématiques. Comme concepts de base, Frege prend les opérations irréductibles déjà présentes chez Grassmann en leur ajoutant quelques opérations de la logique de second ordre : implication, quantification universelle et existentielle… Partant de ces dernières, Frege compose neuf propositions (par exemple : [A → (B → A)]) qu'il nomme axiomes et dont la véracité est indubitable. Tous les autres procédés logiques en découlent selon certaines règles (les règles d'inférence, modus ponens, la règle de la généralisation, etc.). En effet, il s'agit de déterminer les constantes logiques (et, ou, si… alors…) de la syntaxe mentale. Ainsi Frege se donne le même objectif qu'Euclide (non atteint par Aristote) : axiomatiser – fixer dans des formes idéales – la connaissance intuitive. Pour plus d'informations voir [Korte 2010]

[3] Gödel et Tarski en particulier montrent qu'une des difficultés principales du concept de la vérité consiste dans la relation entre langage et métalangage ; la vérité est le moyen de la définir. Sur cette question cf. [Pulkkinen 2005]

[4] L'approche de Frege à Schröder reste néanmoins fort critique. Christian Thiel étudie ce problème dans [Thiel 1981 : 21-23 ; Putnam 1982]

logicien change de point de vue[1]. Quant à Schröder, dans un de ses articles consacré au problème de la variété (1892) Husserl se rappelle d'un débat entre l'inventeur de l'*Exakte Logik* et Lotze sur le contenu des concepts basiques (négation, extension, etc.) définissant la structure relationnelle des entités logiques[2].

La lecture des philosophes (tout d'abord de Descartes, Sigwart ou Lotze[3]) conduit Schröder à une réflexion sur des matières qui dépassent le cadre de la science rigoureuse. Dans sa longue introduction, qui fait davantage écho à la pensée cartésienne, il pose la question de la nature de l'abstraction, du concept et de son rapport avec le contenu. La faculté de l'esprit mathématique d'abstraire les éléments et de créer ainsi des classes logiques où des moments particuliers acquièrent de nouvelles qualités (par ex., celle de la généralité) est précisément le processus dans lequel la conscience subjective « donne sur » les vérités éternelles. Plus d'une fois Schröder évoque l'idée de l'immortalité des objets mathématiques, rejoignant en cela la plupart de ses collègues convaincus de l'objectivité absolue – pour ne pas dire la présence divine (Léopold Kronecker) – des nombres. Selon Schröder, l'abstraction est un instrument visant à isoler les éléments de telle sorte qu'ils deviennent une unité ou un isolement qui désigne une catégorie plus élevée par rapport au contenu particulier : « Par l'abstraction en général, les éléments représentés se trouvent isolés, cela donne aussi la possibilité de les reproduire dans cet isolement (*Isolirtheit*) <…> L'abstraction permet donc d'exercer (*ausüben*) les représentations singulières, la méthode qu'on applique aussi à une chose particulière » [Schröder 1966 : 82-83]. D'un tel processus d'abstraction apparaissent les concepts (*Begriffe*) qui, prévient Schröder en critiquant Frege, ne peuvent être complètement vierges de scories psychologiques car tout concept contient en lui-même l'expérience psychique menant à la formation d'un concept[4]. Si les opérations logiques, pour aussi objectives qu'elles aient pu

---

[1] Au § 8 du *Begriffsschrift* Frege parle de la « parité de contenu » (*Inhaltsgleichheit*), un des concepts clés de sa théorie logique dont la reconsidération poussera Frege plus tard de distinguer *Sinn* et *Bedeutung*. En outre, on trouve ses nouvelles idées développées en termes non techniques dans *Die Grundlagen der Arithmetik* (1884) puis dans deux volumes de *Grundgesetze der Arithmetik* (parus respectivement en 1893 et en 1903).

[2] Hua XXI, p. 411. Schröder livre la quintessence de ces débats au début de ses *Vorlesungen über die Algebra der Logik*, p. 99. Voir aussi : [Schwartz 1996 ; Hartimo 2012]

[3] Volker Peckhaus rappelle que les sources intellectuelles des *Vorlesungen* sont de même Cantor et Peirce, en soulignant à juste titre qu'un des tâches de Schröder est l'algèbre relationnelle ; cf. [Peckhaus 1990/91: 194-195; Houser 1990/91: 206-236]. L'idéologie des *Vorlesungen*, pour le dire en un mot, consiste à faire de la logique un calcul pour permettre de manier les concepts en jeu avec précision en l'émancipant des clichés de la langue naturelle.

Comme nous avons remarqué plus haut, Schröder continue les travaux de Boole qui s'appuie à son tour sur les idées de son compatriote John Wallis (1616-1703) dont l'*Arithmetica Infinitorum* (1656) poursuit les études cartésiennes sur les mathématiques. C'est dans ses deux livres intitulés *The Mathematical Analysis of Logic* (1847) et *An Investigation of The Laws of Thought* (1854) que Boole construit son système. En bref, la tâche est de décrire les opérations sur les classes (ensembles) et des énoncés en utilisant un alphabet très restreint : 1) les variables propositionnels $a_1$, $a_2$, $a_3$… ; 2) les copules logiques $\neg$, $\vee$, $\wedge$ ; 3) le signe = et 4) deux ( ). De plus, il crée aussi une algèbre binaire n'acceptant que deux valeurs numériques : 0 et 1. Boole affirme que dans le nouveau système le contenu des propositions reste une question des interprétations. Pour la discussion détaillée voir [Anellis, Houser 1991]

[4] Dans les *Vorlesungen* (Bd. I, p. 211-212) Schröder note que parfois la définition garantie (*verbürgen*) ce qu'elle définit (*des Definirten*) car dans une certaine mesure elle le crée elle-même. Une telle pensée se voit frappée d'anathème de Frege qui accuse l'auteur des *Vorlesungen* des fantaisies et des constructions purement intuitives ; cf. [Frege 1895 : 451-452]

sembler, sont mentales et donc sont effectuées par le sujet, alors la logique doit prendre en considération ce résidu psychologique. Dans le cas contraire nous risquons de demeurer dans l'illusion d'une auto-objectivité découlant non pas du *cogito* mais de la fausse hypothèse que l'on peut dépasser le subjectif en établissant un corridor logique pour les courants psychiques. Le concept, qui se réfère par nature à une chose idéale, a pour contenu – *Wesen*, explique Schröder – quelque chose de commun, « ce qui indique le savoir général d'une chose et, d'un autre coté, son contenu factuel (*faktischen Inhalt*) sur lequel se bâtit sa représentation mentale… » [Schröder 1966 : 83]. De là il s'ensuit que le concept reste en mouvement, même si cela passe inaperçu, et qu'il ne peut pas ne pas accéder à une étendue (*Umfang*). Chaque concept évolue dans son étendue considérée comme son champ de ses sens – ou champ sémantique –, des idéalités composant le contenu du concept. Un exemple qui illustre à merveille cette pensée exposée dans l'*introduction* est la matière : c'est une substance en général qui couvre une grande diversité de grandeurs telles que la masse, la température, etc. entrant dans le concept de la matière. Même si l'« étendue » est toujours à préciser (nous verrons qu'elle est un des concepts-clés de la théorie eidétique dans *Ideen I*), cette précision, c'est-à-dire la limitation de l'étendue, fait partie du déploiement conceptuel. Husserl y discerne le mécanisme de la constitution phénoménologique de la conscience ; en délimitant la sphère de son attention (sphère de l'intentionnalité), elle s'étend jusqu'à se transformer en objet complet de soi-même (cf. Hua III/1, § 16). Notons à propos que nous trouvons chez Husserl moult concepts mathématiques convertis dans le langage phénoménologique.

Le résultat des réflexions philosophiques de Schröder est le suivant : les concepts sont toujours créés par notre choix des objets de la pensée dont le nombre, même infiniment grand, est défini et donné (*Umfangsangabe*). Bref, l'étendue serait une notion idéale qui ne délimite pas tant les objets choisis que – par ce choix – la zone intentionnelle de la pensée : autrement dit, c'est une étendue du monde[1] au profit de la conscience. Sans utiliser « l'intentionnalité », terme brentanien, Schröder s'approche de cette idée. La pensée algébrique s'arrête sur l'objet idéal, dans lequel elle reconnaît toutes les singularités – éléments de la même classe – qui composent désormais son seul contenu ; ces singularités sont reconnues comme les moments idéalisants de la pensée. En voici l'exemple : « On ne peut avoir (*besitzen*) le concept idéal du cercle qu'après la définition de toutes les particularités possibles de tous les cercles et leur relation avec l'esprit, unifiées dans la conscience » [Schröder 1966 : 87]. Husserl répond à ce livre par une longue recension dont la moitié est précisément consacrée à sa partie philosophique. Sur un ton polémique, Husserl veut comprendre l'épistémologie de Schröder qui lui paraît mal fondée ; la faiblesse principale consiste en l'invention de concepts inutiles ne jouant aucun rôle majeur dans la théorie. Ceux-ci ont un caractère autoréférentiel et ne servent qu'à justifier le choix des idées concrètes d'un tel système. Quelle conclusion – demande le *Rezensent* – le lecteur peut-il tirer de ce type de concept ? Voici sa réponse : « évidemment, l'étendue pure (*Umfangsangabe*) n'est en aucun cas le moyen ou le concept à définir car toute détermination conceptuelle est une détermination de contenu (*Inhaltsbestimmung*) ; <…>

[1] Par « monde » nous comprenons ici tout objet ou tout ensemble des objets intentionnels.

[l'idéal de la logique] dont le principe est reflété dans le concept de l'étendue reste imprécis et sans objet (*gegenstandsloses*) » [Hua XXII : 16][1].

Nous avons ici davantage que la critique de la théorie schröderienne. Il s'agit d'une tentative de remettre en question le modèle de la pensée mathématique ; le problème se trouve dans le moyen-même de la production des vérités (dans lequel Henri Poincaré verra plus tard le duende des mathématiques[2]). Elles se déterminent par le libre choix de propositions pouvant ensuite être corrigées selon les résultats. Autrement dit, l'entendement en mathématique est toujours basé sur le syllogisme et non sur l'adéquation d'une hypothèse avec la vérité, qui dépend largement de la méthode choisie par le mathématicien afin de résoudre une tâche concrète.

Les mathématiques sont donc une autre vision des choses et leur caractère abstrait (dont on parle très souvent) ne correspond pas à celui de la philosophie ; le mathématicien cherche toujours la solution de problèmes précis, tandis que le philosophe cherche à poser le problème, à trouver un paradoxe dans toute évidence (pour le philosophe cette dernière n'est pas identique à la vérité), ouvrant ainsi une nouvelle perspective du savoir. Selon Bourbaki, c'est à partir des Grecs que le mot « mathématique » signifie « preuve ». Plus irréfutable est la preuve, plus elle s'approche de l'idéal mathématique. Même si la preuve ne résume pas toutes les « facettes » du travail mathématique (il y en a d'autres, par exemple : le calcul, ébauche de la preuve avec précision, description algorithmique…), c'est elle qui se situe au centre de l'attention des maîtres du XIX[ème] siècle dont les travaux ont créé le climat et peut-être le mode de pensée. Nul doute que la passion de Husserl pour l'austérité scientifique (plus tard pour « la phénoménologie comme science rigoureuse ») vient non seulement de ses années de formation auprès de Weierstraß puis de Brentano, mais aussi du rigoureux *Geist* de son époque. L'époque de la crise des mathématiques est troublante et belle (animée par la recherche intensive des fondements), et son invention la plus audacieuse est probablement la théorie des ensembles de Cantor (voir infra), fondée sur une méthode universelle et extrêmement puissante. Sur le plan général, Cantor propose une nouvelle épistémologie, restée pour beaucoup incomprise (et pour cette raison battue en brèche par Kronecker[3]), qui a définitivement transformé la pensée-même ; la révolution cantorienne (mathématique et épistémologique, sur laquelle nous reviendrons plus tard car elle permet de comprendre la genèse de la phénoménologie-même), consiste à créer de nouvelles unités opérationnelles – les ensembles (*Mengen*) – décrivant tout « objet pur ». Mieux encore : non

---

[1] Concernant ce débat voir [Schneider O'Connell 1988 : 91-125]

[2] Cf. [Poincaré 1902]

[3] Le monde mental de Kronecker s'est attardé en Grèce. Selon lui, toutes les mathématiques doivent reposer sur les nombres naturels. Hajime Tanabe (1885-1962), élève japonais de Husserl, note que la position kroneckerienne est très proche de celle de Hermann Helmholtz qui traite la théorie des nombres comme une théorie de la conscience, où le nombre naturel est un acte de la pensée (*Bewusstseinsakte*). Pour Kronecker, « les nombres naturels sont une suite régulière de signes ou de désignations (*Bezeichnungen*) servant à définir l'ordre des choses » ; [Tanabe 1915 : 99]. Le nombre n'est pas un signe pur, un véhicule vide ; au contraire, il est lié au monde par sa fonction d'ordonnancement des choses. La conviction de Kronecker concernant les nombres naturels, qui a pris chez lui un caractère quasi-religieux, est pour beaucoup à l'origine de ses attaques contre Cantor. Les défenseurs de Cantor, comme Dedekind et Hilbert, ont reconnu à ses travaux le mérite d'avoir opéré un changement de paradigme ; nous connaissons l'aphorisme de Hilbert : « Nul ne doit nous exclure du paradis que Cantor a créé ». Pour plus de détails voir [Edwards 1988; Moore 2002]

seulement la théorie cantorienne décrit un objet ou le construit, mais elle forme aussi l'état mental spécifique, le seul où cet objet peut apparaître[1]. C'est pourquoi, aux yeux de ses collègues et surtout de ses adversaires, Cantor s'avère trop philosophique, c'est-à-dire flou et imprécis ; sa méthode ne mène pas à des résultats concrets[2] (comme par exemple ce qu'on appelle « le produit de Kronecker », opération portant sur les matrices (masses de nombre) de certaine taille résultant dans leur produit commun). Rappelons que le but de fournir aux mathématiques un fondement solide n'a jamais été étranger à Kronecker ni d'ailleurs à Weierstraß[3]. Toute la différence consiste dans le « comment faire » ; pour ces deux derniers, ce fondement ne peut être trouvé que dans l'algèbre et l'arithmétique des nombres naturels. De ce point de vue, Cantor est un vrai théologien qui crée *ex nihilo* des infinités à puissance différente et joue avec elles comme avec les ombres de Dieu. Kronecker est lui un réaliste ; son nombre possède une ontologie et pour cette raison n'échappe pas au contrôle.

Les ensembles cantoriens ont une particularité stupéfiante (comme Cantor lui-même l'écrit dans une lettre à Dedekind e 1877 : « je le vois, mais je n'y crois pas »[4]) : tous les sous-ensembles d'un ensemble X (fini ou infini) possèdent une puissance supérieure à l'ensemble X. A la fin du XIXème siècle cette idée, triviale aujourd'hui, était complètement troublante, pour ne pas dire incompréhensible. Elle résume les propriétés basiques du nombre naturel : chaque nombre mesure une quantité ordonnée et se trouve ainsi en rapport avec le nombre suivant.

Nul doute que la théorie a montré le caractère illusoire de la croyance en un nombre qui semblait *ens entium*, mais s'est avéré au bout du compte le locus le plus problématique des mathématiques car interprété par chacun à sa manière.

Un autre aspect de la théorie cantorienne, moins évidente peut-être, c'est sa réponse à Descartes. En bref, le Dieu cartésien, dont l'épistémologie est particulièrement expliquée dans les *Objectiones* (chap. II), se transforme en « infinité des infinités » chez Cantor ; il s'agit de l'ensemble des nombres transfinis qui n'obéissent pas au classement des nombres naturels. Sur le plan philosophique, Cantor renforce la démarche cartésienne par l'introduction de la transcendance dans l'immanence : si Dieu est la seule cause du véritable savoir, alors il faut avant tout le connaître, le voir en acte. Husserl répétera cette démarche au § 58 des *Ideen I* où il parle de Dieu transcendant dans l'immanence. Si pour la plupart des philosophes les objets de la raison, y compris mathématiques, ont la plateforme psychique, et ce n'est que par cette dernière nous arrivons à comprendre le sens de ces objets, Husserl peu à peu mais avec certitude change la perspective. L'école brentanienne deviendra un très important épisode sur son chemin intellectuel[1]. Pour lui, les vérités ont une tendance à nous s'ouvrir parce qu'elles sont *dans* les choses mêmes. Nous prenons connaissance de ce « dans » comme notre propre phénomène psychique au lieu de voir ces vérités appartenant aux choses.

---

[1] Youri Manine souligne que Cantor a réussi à créer non seulement une belle méthode mais aussi un univers clos, autosuffisant et « <…> autoréférentiel où le grand résultat est obtenu avec un minimum de moyens » ; cf. [Manin 2002]

[2] C'est toujours Kronecker qui développe une critique déterminante à partir du moment où Cantor commence à travailler sur les fonctions à domaines de définition contenant une infinité de points, ce qui se transforme ultérieurement en théorie des ensembles dont les germes existent dans [Kronecker 1871 :294-296]. Quelques années plus tard, l'animosité de Kronecker atteint alors un degré tel qu'il barre à Cantor l'accès au *Journal de Crelle* (très prestigieux de l'époque) et affuble son ancien élève de sobriquets tels que « le charlatan scientifique », « le renégat qui pervertit la jeunesse… » ; cf. [Dauben 1977: 89]

[3] La proximité scientifique entre ces auteurs est étudiée par [Bottazzini 1981, chap. 7.1]

[4] [Ich sehe es, aber ich glaube es nicht].

## CHAPITRE 4. CASUS CANTOR

Weierstraß, directeur de la thèse doctorale de Cantor, invente une théorie selon laquelle tout nombre réel peut être exprimé par une séquence des nombres rationnels (par exemple, on peut présenter $\sqrt{2}$ par une série des nombres 1, 1,4 1,41 etc.).[2] L'élève y a vu une application géométrique : tous les nombres irrationnels peuvent se présenter comme points sur la ligne, de même que les nombres rationnels. Cette découverte de Cantor a reçu un accueil mitigé, il permettait l'existence des ensembles à des éléments infinis, de plus cela contredit à l'intuition qui dominait la plupart des esprits mathématiques de l'époque (mais en confirmant l'intuition bolzanienne[3]). Le dogme aristotélicien concernant le caractère potentiel de l'infini se trouve bouleversé. L'infinité actualisée change la conscience mathématique et avec cela la perception des objets idéals.

On peut appeler la vie scientifique de Cantor « la poursuite de l'infinité actuelle » autour de laquelle sont organisés presque tous ces travaux. Il s'agit d'un cas rare quand la révolution en mathématiques touche profondément aussi bien la théologie que la philosophie de son époque. Pourquoi ? Parce que la théorie cantorienne explose non seulement l'intuition mathématique « normale », mais aussi la pensée onto-théologique venant de saint Thomas et Suarez à travers Descartes pour laquelle le savoir de l'infini n'est pas actualisable. L'opinion de saint Thomas, fidèle à la conception aristotélicienne, est négative envers l'infini actuel mettant en doute la puissance divine. Même chez C.F. Gauss l'infini est un outil dialectique, car il distingue deux types d'infini : potentiel et actuelle, la véritable pour lui est potentielle. Il disait à un son collègue (H. Schumacher) que l'infinité au sens propre n'est jamais accomplie ; en effet, la notion de l'infinité en mathématiques est un moyen de démontrer la notion de la limite.

Dans *Ideen I* Husserl transformera l'infinité cantorienne en *épistémologie purement phénoménologique*, c'est-à-dire quand la conscience actualise son objet de telle manière qu'il perd son caractère objectif au sens d'exister dans le monde. Il semble que l'origine du concept de l'essence pure (*reine Wesen*), le concept clé qui ouvre *Ideen I,* outre Aristote, il faut chercher dans la théorie de Cantor.

Un certain pas vers l'infinité actuelle a été fait par Dedekind (cf. infra). L'on considère qu'une ligne est plus riche des points qu'une région des nombres naturels. A cela Dedekind répond que malgré la densité des points sur un segment, il est toujours possible d'y encrister un nombre infini des points irrationnels. Ainsi, tout segment a des « fentes » démontrant son caractère discontinu. Cette opération s'appelle « la coupure » (*Schnitt*), une généralisation quasi-philosophique d'un nombre. Plus tard, dans l'*Einleitung* 1906/07 et dans les *Ideen*, Husserl servira cette généralisation de la même manière dedekindienne dont il connaît bien grâce à sa formation.

---

[1] Dans ses leçons R. Ingarden confirme aussi que Husserl désapprouve la méthode psychologiste (brentanienne) qui traite les nombres comme instruments descriptifs de la psychique ; cf. [Ingarden 1974: 21]

[2] Traditionnellement l'on dessine ces nombres par une droite des réels : $1,\ -1,\ \frac{1}{2},\ 0{,}12,\ \pi,\ \sqrt{2}$ ...

[3]Pour le bilan historique cf. [Tieszen 1989]

Lorsque Dedekind ou Cantor inventent des nombres inexistants auparavant, cela ne signifie qu'une chose : ils construisent une abstraction permettant de trouver la solution la plus simple et directe à un problème mathématique précis. Dès qu'une telle généralisation, ou abstraction si l'on veut, s'avère ne pas être au service de cette solution, à la simplifier et clarifier, elle est abandonnée comme étant inutile. Ainsi fonctionne l'esprit du mathématicien.

Nombre inventé ou abstraction, que signifient-ils phénoménologiquement ? Nous verrons que Husserl se concentre sur cette question pendant les années de recherches intensives qui ont précédé *Ideen I*. Pour l'instant, tenons-nous-en à l'observation suivante : venu lui-même des entrailles de la mathématique de son époque, Husserl y emprunte non la méthode – qu'il finit par surpasser –, mais les possibilités qu'elle découvre, celles qui présentent un intérêt sur le plan philosophique. Husserl en tire une leçon capitale : l'existence d'une chose ou d'un objet – idéal – (d'un nombre, par ex.) diffère absolument de celle de tous les autres. L'objet idéal n'est jamais *donné*, c'est-à-dire donné par son existence, il n'équivaut jamais à cette dernière. C'est pourquoi, quand on traite ce type d'objets, il faut radicalement changer le concept-même de l'existence. L'existence ou la non-existence d'un objet idéal ne dépend d'aucune condition si ce n'est le travail de l'abstraction. C'est par ce travail, capable de changer l'ontologie même de l'existence, que naissent les formes de l'idéalité (formes au sens aristotélicien, transmettant seulement l'idée de l'existence).

Anticipons la suite : le concept husserlien de l'idéalité viendra du caractère limité de l'idéalité mathématique ; se bornant au champ de son application, elle n'apparaît que des opérations mathématiques-mêmes. Jean Desanti précise : « je peux l'effectuer soit comme racine d'une équation algébrique, soit comme coupure sur l'ensemble des rationnels, soit comme limite d'une suite Cauchy[1] sur l'ensemble des réels, etc. <…> toute effectuation du signifié, en tant qu'elle est explicite, appelle et mobilise le système disponible des effectuations déjà produites, et demeure à son tour, dans un tel système, comme une configuration déterminée » [Desanti 1968 : 237]. En revanche, les idéalités philosophiques surgissent *comme* saisie de l'instant de l'abstraire. L'idéalité phénoménologique est dirigée non pas vers l'objet – pour aussi abstrait qu'il soit – mais vers l'instant de conscience. Ainsi, pour l'algébriste, $\sqrt{2}$ est un nombre ; pour le phénoménologue le moment propre de sa conscience. L'idéalité mathématique existe donc dans la mesure où l'on peut la définir ; l'idéalité philosophique dans la mesure où nous pouvons saisir – phénoménologiser – l'instant du définir.

[1] Une « suite de Cauchy » est une suite de réels, de complexes, de points d'un espace métrique (un ensemble au sein duquel une notion de distance entre les éléments est définie), dont les termes se rapprochent à partir d'un certain rang :

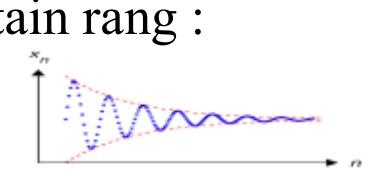

Husserl fait son premier pas vers le dépassement de la restriction mathématique dans sa recension du livre de Schröder[1], où il critique la théorie schröderienne de l'abstraction comme unité des singularités. En doutant de l'idée de Schröder, Husserl met à l'épreuve le dispositif général des mathématiciens : au lieu de créer une nouvelle espèce d'objets jamais offerts par leur simple présence, il faut se méfier de la validité épistémologique de cette espèce. Chez Schröder, comme chez les autres, il s'agit seulement d'un procédé permettant de trouver le moyen de résoudre un problème ; dès que la solution est trouvée, l'abstraction se termine. La revendication-maîtresse de Husserl – encore peu explicite dans sa critique de Schröder – consiste en ce que l'opération de l'abstraction, en remplissant une fonction cruciale en mathématique et, dans la plupart des cas, en visant une tâche concrète, ne crée jamais le sujet pensant. Celui-ci n'existe pas dans cette science, remplacé qu'il est dès le début par un opérateur (ou un artisan, comme le dira plus tard Husserl) doté d'une pensée opérationnelle. La mathématique n'a pas de sujet parce qu'elle ne met jamais en doute son existence. Le modèle mathématique se contente de décrire une espèce ou une classe de choses (soulignons le mot « modèle » car il caractérise le mieux l'essence des recherches : la création de modèles et leur application). La vérité des mathématiciens est une vérité venant de la preuve et de la démonstration[2], notamment de celle qui se divise en une série d'étapes bien précises ; la vérité philosophique surgit de l'état de la conscience où elle se met elle-même à l'épreuve.

La *Rezension* du livre de Schröder marque un tournant vers des réflexions sur la nature du savoir mathématique, systématisées en 1891 dans la *Philosophie de l'arithmétique*. La nouvelle preuve offerte la même année par Cantor[3] qu'il existe une hiérarchie des ensembles à puissance différente, c'est-à-dire que deux infinités peuvent être inégales, semble ne rien devoir au hasard. L'esprit critique exige alors un renouvellement de la pensée mathématique qui était dépourvue d'un modèle cohérent et rigoureux de ses propres fondements. La tentative de Hoëné-Wroński d'algorithmiser la science exacte en lui faisant gagner en exactitude, puis le travail de Schröder visant à apporter une base logique à l'algèbre ne firent que provoquer d'autres essais semblables, plus ou mois concluants. Comme le dit Pieri en 1900 au *Congrès International de Philosophie* de Paris, exprimant probablement là un point de vue commun : « La logique mathématique ressemble à un microscope propre à observer les plus petites différences d'idées, différences que les défauts du langage ordinaire rendent le plus souvent imperceptibles, en l'absence de quelque instrument qui les agrandisse.

Quiconque méprise les avantages d'un tel instrument, notamment dans cet ordre d'études (où souvent l'erreur résulte d'équivoques et de malentendus dans des détails en apparence insignifiants) se prive à mon avis, de propos délibéré, du plus puissant auxiliaire qu'on possède aujourd'hui pour soutenir et diriger notre esprit dans les opérations intellectuelles qui réclament une grande précision » [Pieri 1901 : 382]. Sans entrer dans les détails, notons que ces tentatives poursuivaient le même objectif : réduire voire, dans la mesure du possible, éliminer le hasard du raisonnement mathématique. Alessandro Padoa, autre intervenant au *Congrès*, nomme cet objectif « l'irréductibilité des propositions et symboles primitifs » [Padoa 1901 ; 1899] ; cela veut dire que toute proposition primitive ne peut être déduite d'une autre plus primitive et que tout symbole primitif ne peut être construit sur un autre plus primitif. Nous voyons donc la « méthode de Padoa », ainsi appelée depuis et dont la nouveauté ne sera reconnue que trente ans plus tard par Tarski [Tarski 1956 : 296-319], converger avec l'idée de Hilbert de créer un système de primitifs. D'où résulte l'animosité entre ces deux penseurs dont les idées se révèlent trop proches ; lors du *Congrès* Hilbert a ignoré de manière ostentatoire la communication de Padoa alors qu'il connaissait pertinemment son sujet.

[1] Cf. supra.

[2] On sait qu'Aristote, premier théoricien de la démonstration, la considère comme le seul véritable moyen de connaître et de construire la science ; cf. ΑΝΑΛΥΤΙΚΩΝ ΥΣΤΕΡΩΝ (*Derniers analytiques*), Γ : 72a.

[3] Il est d'usage de dater la naissance de la théorie des ensembles de son article *Sur une propriété caractéristique de tous réels algébriques* (*Über eine Eigenschaft des Inbegriffes aller reellen algebraischen Zahlen*). Le drame de cette théorie, qui est aussi celui, personnel, de Cantor, commence lorsque, malgré la réticence farouche de Kronecker et grâce à Dedekind, ce texte est publié en 1874 dans le *Journal de Crelle,* Bd. 77, 258-262. Cantor montre que les nombres réels ne sont pas dénombrables (la démonstration simplifiée de ce résultat n'apparaît que dix-sept ans plus tard avec la fameuse méthode diagonale ; cf. infra).

On croyait que la logique donnait aux mathématiques ce qui leur manquait : un choix limité qui guide vers la certitude absolue. On croyait aussi, comme Boole, Frege ou Peano, que le parfait langage descriptif pouvait produire au mieux son objet, c'est-à-dire le *faire* connaître. Les entreprises de ces deux derniers auteurs se distinguent par leur caractère rabelaisien : l'*Idéographie* de Frege, écrite sous une forme entièrement formalisée, se transforme en cryptographie. Ce texte présente une immense série de symboles « d'une effroyable complexité et fort éloignés de la pratique des mathématiciens » [Bourbaki 2007]. Pendant vingt ans, Peano rédige avec ses collègues les cinq volumes du *Formulaire de mathématiques* (*Formulario matematico*) qui vise à cataloguer en langue symbolique toutes les propositions mathématiques de base (les quatre premiers volumes sont parus en français, le dernier en latin simplifié[1]) ; il commence à travailler dès 1888 sur ce projet, dont il souligne l'importance dans sa lettre à Camille Jordan[2]. D'autres, tels Russell, Whitehead ou König[3], supposent qu'il suffirait d'atteler la pensée mathématique au collier de la logique pour la sortir de l'embarras.

[1] Cet ouvrage va de pair avec ceux de l'époque des Lumières ou de la Renaissance. Il suffit de citer *Summa* (1494) de Luca Pacioli, *Cursus seu Mundus Mathematicus* (1690) de Claude Descales ou *Elementa Matheseos universae* de Christian Wolff. Quant à la publication du *Formulario*, elle a été rendue possible grâce aux efforts des élèves de Peano (R. Bettazzi, C. Burali-Forti, G. Vailati parmi d'autres) qui, comme leur maître, croyaient en utilité ultime d'une telle entreprise. Peano avoue : « Une telle collection, extrêmement longue et difficile dans une langue ordinaire (*linguaggio comune*) profite vraiment de la natation de la logique mathématique… » ; cf. [Peano 1892 : 76]

[2] « C'est la première fois que l'on applique la logique mathématique à l'analyse d'une question de mathématiques supérieures ; et cette application est, selon moi, la chose plus importante de mon travail » (6 novembre 1894) ; cf. [Conte, Giacardi 1991 : 96]

[3] Julius König (1849-1913) présente un intérêt particulier. Hongrois d'origine, il s'installe en 1869 à Heidelberg, où il soutient l'année suivante une thèse sur les équations modulaires des fonctions elliptiques, puis pendant six mois à Berlin où il suit les cours de Weierstraß et Kronecker, qui a sa préférence. Comme Husserl avocat de la logique pure, König publie en 1914 ses *Nouveaux principes de la logique* (*Neue Grundlagen der Logik*), où il va jusqu'à affirmer que seul l'axiome d'une logique fondée sur la logique pure, science infaillible par excellence, peut être le fondement ultime du savoir comme tel. Dans cet ouvrage, l'auteur évoque quelques motifs phénoménologiques. D'après Marcel Guillaume (communication personnelle) si König n'avait commis certaines fautes, il eut pu produire la première axiomatisation de la logique brouwérienne des propositions ; cf. aussi [Guillaume 2008]

Que cela signifie-t-il techniquement ? Une chose simple en apparence : trouver le nombre total des axiomes (ou des outils déductifs) à partir desquels on peut déduire toutes les propositions logiquement vraies. Il est nécessaire que le nombre de ces axiomes soit limité. Autrement dit, l'idée que toutes les propositions vraies de l'arithmétique sont des axiomes doit nous satisfaire. Hélas, c'est là une illusion car, afin de prouver la véracité d'une seule proposition, il nous faut effectuer un nombre infini d'actions arithmétiques.

Pourtant, les efforts ultérieurs, à commencer par ceux de Hilbert, géomètre et algébriste dont Husserl a suivi les cours à Göttingen[1], de construire une axiomatique infaillible à partir de laquelle on peut déduire toute la mathématique et rendre cette déduction plus exacte encore, conduisent à la naissance d'un regard philosophique. La raison purement scientifique repose sur les axiomes, c'est-à-dire sur une structure idéale qui permet de définir les objets et d'effectuer avec eux des démonstrations. Définir, c'est-à-dire admettre l'existence de tel ou tel objet (malgré qu'Aristote récuse l'idée d'une définition impliquant l'existence du défini sans postulat ni démonstration). Il s'agit en particulier de justifier les raisonnements qui supposent un infini existant en acte, les raisonnements « transfinis » (cantoriens) en faisant l'économie de l'hypothèse de l'existence de l'infini. La science rigoureuse, coutumière du jonglage avec des matières si abstraites, commence à voir non seulement ses éléments (nombre, point, etc.) mais aussi à se voir elle-même. Ironie de l'histoire philosophique, Hilbert – qui partage le rêve leibnizien d'un langage logique idéal[2] –, en aidant Cantor à corriger l'erreur d'Aristote sur la nature de l'infini, devance Husserl dans la pratique de l'ἐποχή[3] et dans la création de son axiomatique phénoménologisée.

Si, comme nous l'avons vu, le progrès en mathématiques est obtenu par l'abstraction, alors la théorie des ensembles est une abstraction incoercible de deux aspects du nombre naturel (0, 1, 2, 3, 4…)[4] : tout nombre mesure la

[1] Husserl et Hilbert se rencontrent à Göttingen à l'automne 1901. Les deux hommes ont presque le même âge et beaucoup d'intérêts en commun. Ils se croisent à la Faculté de philosophie, l' « agora » des philosophes et des mathématiciens. Dans cette ville la philosophie a été si étroitement associée avec Lotze (qui l'a enseignée entre 1844 et 1881) que Husserl parle d'un « désert philosophique » laissé après sa mort ; cf. [Hua Dokumente, III/5:179]

Néanmoins, à Göttingen les conditions sont plus propices, comme en témoigne Malvine Husserl dans une lettre de décembre 1901 : « Ici, le climat intellectuel de l'université est fort différent de celui de Halle, et en particulier les mathématiciens Klein et Hilbert qui ont introduit Edmund dans leur milieu. Ils l'encouragent tant et si bien qu'il a récemment présenté à la Société mathématique un exposé basé sur ses anciens manuscrits… », Hua XXI, p. XIII. Il s'agit d'une double conférence des 26/11/1901 et 10/12/1901, qui suivait celle de Hilbert sur les systèmes d'axiomes. Pour l'étude détaillée de cette question voir [Schuhmann, Schuhmann 2001 : 87-123]. C'est un moment décisif pour l'un et l'autre : Hilbert ordonne la méthode de travail sur les abstractions ; Husserl, déjà auteur des *Prolégomènes à la logique pure* (cf., par ex., § 70 écrit en accord total avec l'esprit qui y régnait), continue à travailler sur la logique et sur les nouveaux concepts et se met ainsi en route vers la phénoménologie. Voir aussi [Gandt 2004 ; Plessner 1959].

[2] Non moins que Frege et Peano, pour lequel le projet leibnizien de créer une *scriptura universalis* ou *scientia generalis* est absolument fondamental, Hilbert se joint à ce projet d'autant plus volontiers après la publication du livre de [Diels 1899 ; Leibniz 1999]

[3] Nous verrons plus loin la fonction de ce concept fondamental. Notons ici que Husserl, dans son système, n'explique pas clairement son origine. Klaus Held suggère qu'il emprunte celui-ci à l'éthique hellénistique qui recommande une suspension à ceux qui cherchent l'εὐδαιμονία (le bonheur) ne venant que par l'abondance des *doxas*. Cf. [Held 1989]

[4] Le zéro est généralement considéré comme le seul nombre naturel non-positif. Dans une lettre à Ch. Hermite du 30 novembre 1895, Cantor explique que « les nombres naturels existent au

quantité et tous les nombres sont ordonnés par leur relation « $x < y$ ». De cette évidence Cantor déduit l'existence de deux types d'infinité : *cardinal* et *ordinal* (désignés par la lettre hébraïque ℵ et par la lettre grecque **ω** respectivement)[1]. Le cardinal est le nombre d'éléments d'un ensemble (par ex., le cardinal de l'ensemble {1, 2, 3, 4} est 4 ; « comme les ordinaux, les nombres cardinaux sont des concepts simples (*einfache Begriffsbildungen*) dont chacun est une vraie unité (μονάς) » [Cantor 1932: 380]. Plus généralement, la cardinalité est une propriété de l'ensemble qui caractérise celui-ci comme singularité complète.

L'ordinal est un nombre qui décrit la position numérique d'un objet (0, 1, 2, …) ; dans la théorie des ensembles, l'ordinal est encore un segment structural de l'ensemble ordonné. En termes formels, l'ensemble est un concept décrivant un objet concret avec ses particularités. Mais Cantor va plus loin : l'ensemble n'est pas seulement l'ordonnancement des éléments, y compris de ses parties (sous-ensembles), mais aussi la prise de connaissance de cela. Pour qu'un ensemble se déduise d'un autre, il ne suffit pas de le traiter comme une figure du savoir mathématique ; il faut que l'ensemble soit l'état de conscience-même dans lequel se trouve le sujet en chemin vers l'absolu. Un ensemble ne peut être vu entièrement – c'est-à-dire dans toutes les parties constituant son unité – que par la conscience dirigée non seulement sur ce qui existe maintenant mais sur ce qui peut exister (à comparer : voir les trois aspects du temps simultanément).

Phénoménologiquement parlant, dans un ensemble entier la différence entre « est » et « peut » s'annule en dépassant leur cadre ontologique. Au lieu de discerner ces deux positions, comme le fait Aristote, la raison ensembliste les néglige, comme d'ailleurs toutes les autres différences venant du domaine empirique. A vrai dire, rien dans l'ensemble cantorien ne correspond aux choses réelles ni même les simples nombres ordinaux qui ne sont que des actes mentaux dirigés vers l'extension du champ conceptuel du nombre, d'où le cardinal. Nonobstant le fait que Cantor traite les ensembles comme des êtres vivants, établissant entre eux une parenté, leur réalité n'est pas donnée mais reste toujours en donation. C'est déjà une réalité autre que celle des objets mathématiques que nous avons appris à appréhender comme entités idéales. Les ensembles ne sont pas idéaux au sens algébrique ou géométrique, ils ne sont pas – comme le dirait Husserl – des noèmes à la différence du nombre ou du point. Ces derniers sont figés (noèmes mortes), alors que l'ensemble est un espace noétique qui engendre des abstractions en tendant lui-même vers la complétude. C'est ici qu'il faut chercher la genèse de la conscience absolue élaborée par Husserl dans ses *Ideen* ; non seulement l'ensemble existe simultanément dans des temps différents, non seulement il accueille les diverses possibilités ontologiques en les transférant au plus haut niveau d'abstraction où elles se réalisent au-delà de l'intuition habituelle mais, ce faisant, l'ensemble devient lui-même mode de pensée. En pensant, il neutralise la conscience liée à l'ontologie, ce qui constituera le thème central de la phénoménologie tardive.

---

degré le plus haut de la réalité comme idées éternelles dans l'*Intellectus Divinus* » ; cité par [Dauben 1990 : 228 ; de même Décaillot 2008]

[1] Le plus petit cardinal infini est désigné par $\aleph_0$ ; le cardinal immédiatement supérieur est $\aleph_1$, etc. Pour une belle explication de cette notation voir [Verriest 1951 : 28]

Quoi qu'il en soit, la découverte de deux espèces d'infini, dont les relations dépassent le cadre de l'arithmétique classique, s'est finalement soldée par un double effet : une mathématique nouvelle et la folie de son auteur[1].

Avant d'aller plus loin il faut bien comprendre : la révolution cantorienne n'est pas sans rappeler la révolution cartésienne (la seconde n'est-elle pas la continuation de la première ?). Grâce à Cantor en mathématiques apparaît le *sujet* qui, en considérant le monde classique des nombres et des points comme un cas particulier, crée des abstractions d'ordre supérieur[2]. Ainsi, l'échelle des ensembles $\Omega_0 < \Omega_1 < \Omega_2 < \Omega_3 <\ldots$ qui croît partout signifie qu'il reste toujours un élément non-inclus dans une série d'abstractions. Ce résidu continuel engendre un temps spécifique du sujet, celui-ci y demeure sans garder aucun rapport avec le temps empirique ; Husserl reformulera cette idée pour son Je transcendantal qui n'existe que dans le temps de sa conscience. Kronecker, ex-maître de Cantor, l'accuse d'être animé d'une perversité qui menace la nature de la science en donnant le mauvais exemple aux jeunes esprits. Cette animosité s'explique par deux raisons majeures : l'une personnelle et l'autre professionnelle.

D'emblée, Kronecker est un extrémiste ; convaincu de l'existence des nombres entiers[3] uniquement, il veut fonder sur ceux-ci tout l'édifice de la mathématique. Il refuse toute opinion autre que la sienne (par exemple, celle des

---

[1] L'existence d'un deuxième type d'infini ouvre l'horizon au sens husserlien du terme permettant de parler d'« infini des infinités », monstre conceptuel qui répandait la terreur parmi les bons esprits universitaires. Cantor fait à maintes reprises allusion au caractère divin de son idée ; par exemple, dans *Mitteilung zur Lehre vom Transfiniten* (1887), il souligne que l'infini actuel (ou les nombres transfinis) existe dans un « être hors-du-monde » (*außerweltlichen Sein*), *in Deo*, dans l'état le plus parfait qui soit ; cf. *Gesammelte Abhandlungen*…, p. 378. Puis en critiquant Wundt, qui a mal compris la nature du transfini, Cantor insiste sur l'actualité de l'infini (rejoignant sur ce point Bolzano) et inscrit le nombre transfini dans la tradition scolastique et thomiste (p. 386).

Il consacre le reste de sa vie à des questions philosophiques et théologiques, cherchant entre autres la présence des nombres transfinis dans la religion. Métaphysiquement, Cantor s'oppose à certains auteurs néothomistes (G. Pecci, T. Zigliara ou J. Kleutgen) qui tendent à distinguer l'infini actuel et la connaissance divine, convaincu qu'il est que sa théorie servira à invalider cet égarement ; cf. [Dauben 1977 ; Feferman 1987]

[2] Les débats autour de la naissance du sujet mathématique sont déclenchés par un événement : en 1904 au *III^ème^ Congrès International des mathématiciens*, Ernst Zermelo présente son *Preuve. De la possibilité pour tout ensemble d'être bien ordonné* (Beweis, daβ jede Menge wohlgeordnet werden kann, publié dans le volume 59 des *Mathematische Annalen*). Dans ce texte, il démontre une idée capitale : tout ensemble peut présenter une structure de bon ordre ; cela signifie que tout sous-ensemble (non vide) a un plus petit élément. Sur le plan personnel, la communication s'avère très favorable pour Cantor qui, après la critique de König sur le problème du continu, se sentait vulnérable. La discussion initiée par l'idée de Zermelo et publiée dans le volume suivant des *Math. Annalen* concerne les principes psychologiques du Je mathématiquement pensant.

Il est fort probable que Husserl suive cette affaire : premièrement, il entretient une correspondance avec Zermelo quelques années durant ; deuxièmement, il nourrit des sentiments amicaux pour Cantor, dont il garde les travaux dans sa bibliothèque personnelle (je remercie le Dr. Thomas Vongehr, secrétaire des Archives Husserl à Louvain pour ses précisions concernant ces matériaux).

[3] On cite souvent sa phrase : « Le Bon Dieu créa les nombres, le reste est l'œuvre de l'homme » (Die ganzen Zahlen hat der liebe Gott gemacht, alles andere ist Menschenwerk) ; cf. [Weber 1893:19]. Jacqueline Boniface précise que, d'après Kronecker, à partir de ces nombres et de lettres considérées comme indéterminées, le mathématicien ouvre, c'est-à-dire construit des expressions algébriques qui constituent les phénomènes qu'il aura pour tâche de décrire », cf. [Boniface 2010]

constructivistes). A l'opposé de ses collègues Cantor, Peano ou Dedekind, qui manient des abstractions de différents ordres, Kronecker veut voir toute la mathématique devant ses yeux sans oser jouer avec le démon de l'infini actuel. Son trait le plus redoutable consiste à négliger le contenu d'une série dont la non-contradiction formelle s'avère suffisante ; cette série n'est que la visualisation de l'infini, du principe métaphysique permettant de monter dans une suite d'éléments sans arrêt et sans atteindre jamais le plus grand. – « Cette infinité actuelle, à quoi sert-elle, s'interroge Kronecker, si elle se comporte de la même manière dans chacun de ses secteurs ? »

L'intention du critique est claire : emmener la mathématique le plus loin possible dans la métaphysique, l'empêcher de s'engager sur la voie d'une mauvaise réflexion. La peur de Kronecker vient en effet de ce qu'Hugo Dingler (1881-1954), philosophe des sciences et disciple de Husserl à Göttingen, appelle « la régression infinie », qui apparaît chaque fois que nous essayons de trouver l'ultime fondement d'une science en employant son propre langage [Dingler 1931: 32]. Toujours selon Kronecker, la mathématique doit, au lieu de produire des formes pures, s'appuyer sur son efficacité ; la peur des « abstractions incontrôlables» a eu pour conséquence le refus de constituer le sujet mathématique[1]. Pour avoir travaillé sous sa direction, Cantor est familier des vues de son maître. Au début, il justifie cette position et propose un critère clair pour la démonstration, mais Cantor ne tarde pas à en percevoir la faiblesse qui limite la liberté de l'œuvre mathématique. Au fur et à mesure, l'élève se rapproche de la sensibilité mentaliste, pour ne pas dire métaphysique. Kronecker, bien évidemment, ne pardonne pas cette « trahison ». De surcroît, Cantor défend le formalisme : tout nombre existe s'il est possible de le définir de manière précise, c'est-à-dire de le distinguer des autres nombres. Exister en mathématiques signifie être défini. Par exemple, il nomme dans ses premiers travaux le symbole $\infty$ « nombre transfini » en introduisant une nouvelle série de nombres : $\infty$, $\infty + 1$, $\infty + 2$, ...[2] La définition d'un nombre équivaut à son engendrement ; la série est un espace où le nombre trouve son identité.

Pourtant, Cantor s'est avéré plus dangereux encore pour les algébristes classiques que ne le soupçonnait Kronecker. Il introduit la conscience subjective dans la mathématique, tissant ainsi des fils invisibles avec Descartes ; la pensée de ce dernier entre dans une nouvelle phase qui trouve son accomplissement dans *Ideen I*. Le bon mathématicien définit son objet mais ne le pense pas, ne le voit pas ontologiquement[3]. L'objet mathématique n'existe pas sans les opérations rigoureusement définies que l'on peut lui faire subir, cet objet est un ensemble d'opérations desquelles on déduit une structure abstraite généralisant leurs propriétés. Ainsi les nombres naturels 1, 2, 3… sont des objets contenant en eux-mêmes (analytiquement) les opérations qu'un mathématicien effectue dans son travail. Mieux encore, l'objet mathématique est une action qui résulte dans un « objet » ; c'est pourquoi le mathématicien construit les objets d'actions

---

[1] Notons que l'école néokantienne de Bade, avant tout Emil Lask et Heinrich Rickert avec qui Husserl entretient de bonnes relations, propose au lieu du sujet une *conscience impersonnelle* (das unpersönliche Bewußtsein) permettant de construire une nouvelle objectivité. Du point de vue de l'individu, les énoncés de cette conscience sont objectifs par nature ; voir particulièrement : [Rickert 1892 : chap. XVII]

[2] L'histoire précise du transfini est étudiée par [Wallace 2009]

[3] Desanti note : « la mathématique se réduit à un signifié privé de signifiant », cf. [Desanti 1968: 81]

qu'il conviendrait d'appeler généralisations[1]. Penser l'objet n'est pas dans la nature des mathématiques car l'objet pensé n'est jamais défini, il n'existe pas dans des limites fixes.

La théorie des ensembles quant à elle avance des concepts (transfini, puissance, etc.) qui transforment la *définition* en processus de pensée. Le mathématicien commence à se dégager de l'objet et devient peu à peut sujet dès lors qu'il se heurte aux difficultés de fixer les conditions initiales de son raisonnement. Pour être exact, ce tournant épistémologique est précédé des changements capitaux intervenus dans la géométrie grâce aux travaux de Carl Gauss et de Nikolaï Lobatchevski, qui donnent naissance à la géométrie non-euclidienne. En bref, les objets géométriques, qui semblaient avoir des propriétés fixes une fois pour toutes, n'en ont plus indépendamment du travail du géomètre et leur rapport avec la réalité empirique peut s'interpréter de plusieurs façons. Tout objet dans la géométrie n'est pas à décrire comme donné mais à construire, le concept du « modèle » est situationnel.

A l'instar de l'*ego* cartésien, le sujet cantorien se crée lui-même à partir de sa pensée. En revanche, à la différence de l'*ego* cartésien, le sujet chez Cantor prend connaissance non pas de son « cogito… » ni donc de son présent, mais de la possibilité du penser tout entier. La métaphysique de la théorie des ensembles se résume à ceci : Pourrait-on construire à partir des opérations mathématiques non pas des objets mais une *nouvelle conscience* ? Car c'est là le seul moyen de comprendre les nombres infinis[2]. Certes, le chemin qui mène à cette dernière passe par l'idéalisation : le sujet cartésien idéalise son existence en l'abstrayant dans sa pensée, le sujet chez Cantor idéalise son entendement dans des ensembles qui n'existent que dans l'intellectus Dei. A proprement parler, le Dieu cartésien et cantorien est le même ; la seule différence est que le premier pense quand le second compte.

Nous verrons comment les changements de la subjectivité impulsés par Descartes et Cantor aboutiront au sujet transcendantal et comment Husserl pousse cette subjectivité à son terme logique. Le sujet husserlien cependant se distingue de ses prédécesseurs par son orientation sotériologique : en mettant le monde hors circuit (*Ausschaltung* ; cf. § 31), ce sujet a pour but de créer une conscience d'ordre supérieur qui déduit Dieu d'elle-même tout en demeurant dans son autonomie transcendantale. Risquons une hypothèse : *Ideen I* élaborent l'ascétisme phénoménologique qui résulte de l'ἐποχή que le sujet accepte pour mode de vie. L'ἐποχή est non seulement la rétention des flux empiriques de la conscience, permettant d'approcher de manière plus efficace le problème de l'οὐσία aristotélicien, mais aussi le moteur d'une nouvelle doctrine sotériologique. Ici, le salut vient non pas d'une conscience qui se phénoménalise dans le divin, comme chez saint Thomas d'Aquin ou chez Descartes, mais de la possibilité d'exclure Dieu du champ intentionnel. Le paradoxe tient à ce que cette exclusion est l'acte radical de penser la transcendance ; penser celle-ci non pas comme « objet » ou réalité ultime qui, selon la démonstration cartésienne, engendre en nous l'idée d'elle-même mais comme le temps pur de la pensée. En d'autres termes, Husserl sauve le sujet en le plaçant en dehors de la distinction entre la transcendance et l'immanence. Ce faisant, il n'annule pas mécaniquement cette distinction (comme il peut sembler du § 58 des *Ideen I*) ; c'est le sujet en revanche qui commence à percevoir tout comme ses propres actes mentaux.

---

[1] Euclide dans ses *Eléments* est sans doute le premier à proposer une belle généralisation des propriétés des nombres naturels par la démonstration de leur infinité : s'il y a une série finie de nombres ($p_1$, … $p_n$), il est alors possible de prendre le dernier d'entre eux et d'y ajouter 1 ($p_n$ + 1). C'est un pur exemple de création d'un objet mathématique à partir de l'action.

[2] Malgré l'indication que les transfinis constituent un type de nombres tout à fait nouveau et que leurs lois « dépendent de l'ordre naturel des choses », ils ne sont pas des objets *donnés* auxquels on applique la raison mathématique. Cf. [Cantor 1932: 371-72]. Les transfinis constituent un nouveau type de savoir doté d'une une structure réflexive qui mène à la compréhension de l'absolu.

Quant au sujet cantorien, il est le chaînon intermédiaire entre celui de Kant et de Husserl. Chez Kant le sujet, limité par ses capacités originaires spatio-temporelles, en les appliquant au monde, obtient sa connaissance. Cantor modifie le sujet kantien : étant comprise, sa limitation-même devient le moyen de voir l'existence des espaces infinis. Le sujet est inclus dans les ensembles infinis. La valeur de la théorie des ensembles consiste avant tout en la création d'un langage capable de décrire tout objet mathématique. Mieux, à l'aide de ce langage on peut créer toute construction mathématique, et c'est pourquoi nous pouvons voir dans cette dernière l'objet proto-noématique. Sur le plan méthodologique, la théorie des ensembles est une heureuse simplification de notre vision des choses inextricablement liées entre elles dans les abstractions au niveau de leur essence et dont nous arrivons à comprendre le lien. Disons-le brièvement : le nombre transfini ne peut être *vu* que dans l'ἐποχή de Husserl.

Cantor apprend aux mathématiques à penser, Husserl quant à lui transforme cette pensée en auto-pensée. Voici un résumé des principales idées cantoriennes :

α) les ensembles sont des singularités abstraites, leur nombre est infini ; la puissance des ensembles infinis est différente, donc les infinités sont variées. La définition générale de l'ensemble est la suivante : « Par 'ensemble' (*Menge*) nous entendons toute union *M* de certains et différents objets *m* (appelés les 'éléments' *M*) qui existent dans notre intuition (*Anschauung*) ou dans nos pensées [Cantor 1932: 282] ;

β) Le « nombre cardinal » (ou la « puissance ») désigne le concept qui, grâce à notre active faculté mentale (*aktiven Denkvermögens*) dirigée sur l'ensemble *M*, résulte de l'abstraction de tous ses éléments et de leur moyen d'être ordonnés. Autrement dit, le nombre cardinal est une abstraction de second degré ;

γ) Les ensembles sont unifiables par leurs éléments, un ensemble peut être inclus dans un autre et les ensembles sont comparables par leur puissance.

C'est en 1891 que Cantor invente son célèbre procédé diagonal[1] qui permet de rendre visible la puissance de l'ensemble. Le contenu mathématique de cette invention consiste à montrer que l'ensemble de tous les sous-ensembles d'un ensemble *M* (appelé *ensemble des parties* de *M*) est strictement plus grand que *M*, même si *M* est infini. Illustrons cette idée par un exemple :

[1] Le terme plus technique : l'argument diagonal. Il s'agit de son article intitulé « Sur une question élémentaire de la théorie des multiplicités » (Über eine elementare Frage der Mannigfaltigkeitslehre) paru dans *Jahresbericht der Deutschen Mathematiker-Vereinigung*, Bd. I. Ce célèbre concept cantorien n'échappe pas aux médiations intensives de Husserl ; voir [Gauthier 2004]

Admettons que l'infinité des nombres décimaux entre 0 et 1 est la même que l'infinité des nombres naturels. Dans ce cas, tous les nombres décimaux peuvent être dénombrés dans une liste :

1 ⟶ $d_1 = 0.\ d_{11}\, d_{12}\, d_{13}\, d_{14}...$
2 ⟶ $d_2 = 0.\ d_{21}\, d_{22}\, d_{23}\, d_{24}...$
3 ⟶ $d_3 = 0.\ d_{31}\, d_{32}\, d_{33}\, d_{34}...$
4 ⟶ $d_4 = 0.\ d_{41}\, d_{42}\, d_{43}\, d_{44}...$
.
.
.
n ⟶ $d_n = 0.\ d_{n1}\, d_{n2}\, d_{n3}\, d_{n4}...$
.
.
.

Considérons maintenant le nombre décimal $x = 0.\ x_1, x_2, x_3, x_4\ldots$ où $x_1$ est n'importe quel chiffre autre que $d_{11}$ ; $x_2$ autre que $d_{22}$ ; $x_3$ autre que $d_{33}$, $x_4$ autre que $d_{44}$, etc. Si *x* est un chiffre décimal et est inférieur à 1, alors il doit figurer dans notre liste. Où est-il? Ce *x* ne peut être le premier parce que le premier chiffre de *x* diffère de celui de $d_1$ ; *x* ne peut pas non plus figurer dans la liste car ses chiffres des centièmes diffèrent de ceux de $d_2$ ; autrement dit, *x* n'est pas égal à $d_n$ car leurs $n^{\text{ièmes}}$ chiffres ne sont pas les mêmes ; ainsi, *x* ne figure pas dans la liste. En d'autres termes, nous avons un nombre décimal censé figurer dans la liste mais qui n'y figure pas. Quelle que soit la façon de classer les nombres décimaux, il en reste toujours au moins un hors du classement. Dans ces conditions, le classement de tous les nombres décimaux est impossible, d'où la conclusion : l'infinité des nombres décimaux est plus grande (plus puissante) que celle des nombres naturels.

Cette nouvelle infinité est obtenue par l'abstraction des nombres réels ; Cavaillès précise qu'un tel « processus abstrait ainsi défini ne comporte aucune limitation » [Cavaillès 1947 : 7]. Ici *z* est un nombre réel, *a* un nombre décimal de ces réels et la diagonale désigne un nouvel ensemble $\{d_{11}, d_{22}, d_{33}, \ldots\}$ dit dénombrable qui n'existe pas parmi les ensembles des nombres ordinaux. D'où cette nouvelle conception de la finitude : l'ensemble est fini si ses nombres ordinaux et cardinaux coïncident[1].

Or, pour illustrer le procédé diagonal de façon encore plus simple, ébauchons une matrice de séries infinies (E) où la série diagonale ($E_n + 1$) ne peut avoir lieu dans aucune série de type E :

[1] Pour une explication alternative, voir [Belna 2000 : 85]

$E_0 = \boxed{m}\ m\ m\ m\ m\ m\ m\ m\ m \ldots$

$E_1 = w\ \boxed{w}\ w\ w\ w\ w\ w\ w\ w \ldots$

$E_3 = m\ w\ \boxed{m}\ w\ m\ w\ m\ w\ m \ldots$

$E_4 = w\ m\ w\ \boxed{m}\ w\ m\ w\ m\ w \ldots$

$E_5 = m\ w\ m\ w\ \boxed{w}\ m\ w\ m\ w \ldots$

$E_6 = m\ w\ m\ w\ w\ \boxed{m}\ w\ w\ m \ldots$

$E_7 = w\ m\ m\ w\ m\ w\ \boxed{m}\ w\ m \ldots$

$E_8 = m\ m\ w\ m\ w\ m\ w\ \boxed{m}\ w \ldots$

$E_9 = w\ m\ w\ m\ m\ w\ w\ m\ \boxed{w} \ldots$

$E_n = \ldots\ \ldots \qquad\qquad \boxed{\ldots}$

$E_n + 1 = m\ w\ m\ m\ w\ m\ m\ m\ w \ldots$

C’est ainsi que se donne à voir la révolution métaphysique en mathématiques : les infinités sont comparables mais non égales. En les comparant comme si elles étaient des objets normaux, nous ne comparons en effet que des intuitions visualisées dans une possibilité d’ordre supérieur. A la suite de Maurice Fréchet [Fréchet 1934 : 18], pour lequel le principe-clé de la théorie des ensembles consiste en la démonstration de l’inégalité des infinités, nous constatons qu’il s’agit ici d’un procédé phénoménologique. Cette inégalité est une abstraction de la pluralité des infinis ou, plus exactement, de leur existence actuelle. Ce n’est pas l’admission de cette pluralité qui est phénoménologique – ce serait alors une autre démarche mathématique – mais sa création conscientielle.

Finalement, on peut considérer comme preuve de la nouvelle conscience le célèbre axiome du choix proposé par Zermelo en 1904. Son idée consiste à voir l’objet mathématique comme construit, c’est-à-dire comme le résultat d’une procédure subjective. Il est évident que cet axiome doit son existence à la métaphysique cantorienne qui a fait entrer le sujet en scène. Non seulement cet axiome implique le sujet pensant, mais il montre sa nécessité gnoséologique car rien excepté le sujet ne peut opérer un choix libre.

## CHAPITRE 5. LE NOMBRE COMME LIBERTE

Au début des années 1880, lorsque la théorie des ensembles entre dans sa phase dramatique – Cantor et Dedekind échangent leurs lettres surtout sur le concept de la cardinalité – Husserl entame sa carrière scientifique en faisant des relations importantes. En mars 1881 il suit à Vienne le séminaire de Leo Königsberger (1837-1921)[1], ancien élève de Weierstraß, auprès duquel en mai 1882 Husserl termine sa dissertation qui s’intitule *Contributions à la théorie du calcul des variations* (*Beiträge zur Theorie der Variationsrechnung*). Le 2 octobre 1882 Max Büdinger (1828-1902), historien et doyen de l’université, nomme Königsberger et un autre professeur de Husserl, mathématicien Emil Weyr (1848-1894), comme deux rapporteurs de sa thèse. La soutenance a lieu le 29 novembre 1882 après laquelle Weierstraß le prend en tant que son assistant à Berlin durant le semestre d'été 1883. Puis, dès le semestre 1883-84, Husserl revient à Vienne pour poursuivre ses études avec Franz Brentano. En 1887, après avoir soutenu son *Habilitationsschrift*[1], consacrée au concept du nombre (*Über den Begriff der Zahl*)[2] et saluée par Cantor en personne, le philosophe va s'installer peu de temps après à Halle où il se lie d'amitié avec Carl Stumpf.

[1]Il a écrit une thèse intitulée *De motu puncti versus duo fixa centra attracti* qu’il soutient en mai 1860 en devenant professeur au gymnase à Easter pour les trois années suivantes.

Certes, à cette époque Husserl ne peut pas encore participer à tous ces grands débats à armes égales avec ses maîtres ; il rédige ses *notes* et ses articles sur les problèmes spécifiques de l'arithmétique et de la géométrie dans lesquels son intuition philosophique, parfois assez fine, se dirige peu à peu mais avec certitude vers la critique du paradigme classique. Même s'il joue maintenant le second rôle, sa circulation entre Berlin, Vienne et Halle lui permet – lui oblige – d'être au courant de tout. Ce *tout* signifiait jadis le naufrage de la science qui repose en dernier ressort sur le *principio perennis* séparant les objets mathématiques de l'intellect humain (c'est d'ailleurs une des principales raisons de la crise des fondements car s'il existe une objectivité – disons, les postulats d'Euclide – qui ne correspond pas à une idée, c'est-à-dire demeurant incomprise par sa cause, alors cela implique la déficience conceptuelle). Soulignons une fois de plus : cette longue crise a bien décelé tout le caractère illusoire de la distinction entre l'objet et l'activité mentale. Quoi qu'ils ne soient objectifs les objets d'une science, ils ne se définissent que dans cette dernière et leur objectivité se peut voir comme telle uniquement dans les actes mentaux des scientifiques.

En 1887, l'année d'habitation de Husserl, Dedekind termine son article majeur *De ce que sont les nombres et de ce qu'ils doivent être ?*[3] Ni dans la première préface où il cite Weierstraß, Cantor et d'autres dont il estime les mérites, ni dans la préface à la deuxième édition de 1893 où il se réfère à Frege, le nom de Husserl n'y apparaît ni dans le texte lui-même[4]. L'auteur commence par un axiome philosophique : les nombres sont des libres inventions (*Schöpfungen*) de l'intellect qui servent à discerner les choses et leurs espèces[5]. A part cela Dedekind ne définit guère sa position comme entièrement nouvelle en trahissant ses sources – Schröder et Kant – qu'il réunit ainsi : « Je ne considère l'arithmétique (algèbre, Analyse) qu'une partie de la logique, de plus le concept du nombre reste pour moi complètement indépendant des représentations ou des intuitions de l'espace et du temps, il fonctionne en tant qu'une suite (*Ausfluß*) immédiate de la limite pure de pensée (*reinen Denkgesetze*) » [Dedekind 1932: 335].

Cette idée correspond pleinement au travail du mathématicien et est dans une grande mesure son duende métaphysique. On trouve son origine dans l'article de 1872 où Dedekind introduit son concept de la *coupure* qui présente un bon échantillon de l'abstraction en mathématiques[6]. Si $M$ est un ensemble des nombres rationnels, on peut le couper en deux parties (sous-ensembles) $e$ et $e_1$ où tout élément de $e$ sera inférieur à tout élément de $e_1$ et que ce dernier n'aura pas de nombre le plus grand[1]. Une telle série – continuum – des nombres est un espace à deux directions ; nous pouvons se mouvoir, à partir d'un zéro choisi, à gauche ou à droite au temps infini. La coupure donne la possibilité d'insérer dans ce continuum un *je* qui par la limite définie fait abstraction des nombres rationnels pour obtenir les nombres irrationnels incarnant l'idée de l'étendue infinie par excellence (indispensables, par exemple, pour mesurer la diagonale du carré[2]) d'où, si l'on parle en termes husserliens, la coupure est une opération noétique.

[1] Thèse d'habilitation

[2] Il existe un résumé en français de ce travail édité par [Vauthier 1983]

[3] *Was sind und was sollen die Zahlen ?* sorti au début de 1888 et republié dans [Dedekind 1932]

[4] Il est curieux que dans son ouvrage *Tonpsychologie* (Bd. II, Leipzig, Hirzel, 1890), qui a provoqué beaucoup de polémique, Stumpf cite « Herr Husserl ».

[5] A noter : à l'époque nazie et même avant cette remarque métaphysique de Dedekind sera interprétée comme le trait spécifique de l'esprit allemand (*deutsche Geist*) ; pour plus d'informations cf. [Vahlen 1923: 21-22; Segal 2003: 365]

[6] L'article s'intitule *Continuité et nombres irrationnels* (*Stetigkeit und irrationale Zahlen*).

Philosophiquement parlant, une telle coupure conceptualise des objets qui se trouvent entre *e* et $e_1$ compte tenu de ce que cet « entre » n'appartient pas forcément au *M*. Il n'est d'ailleurs pas étonnant pourquoi Dedekind fut l'un des premiers à s'intéresser sérieusement aux travaux de Cantor dont il a fait connaissance en 1874 en Suisse. Deux auteurs travaillent simultanément sur les objets qui ne viennent à l'existence que par l'établissement d'une relation entre elles. Dès lors la relation et la limite, plus exactement la limite mentale – *Denkgesetze* –, composent chez Dedekind le thème central. Pour dire plus rigoureusement, cette théorie montre l'essence de la continuité étudiée déjà par Cauchy sous la forme des fonctions continues[3]. Se posèrent donc deux problèmes :

α) conceptualiser logiquement la propriété fondamentale de la ligne droite qui ne se trouve que dans notre monde sensible. Cela signifie de définir le plus rigoureusement possible le concept de la continuité (ce qui n'avait été fait ni par Euclide ni par aucun autre des auteurs grecs) mais ce qui permet d'étudier ensuite tous les domaines de continu ;

β) construire une théorie arithmétique complète du nombre pour qu'on ne doive plus aller aux intuitions géométriques. Bref le concept du nombre doit être construit sur les principes abstraits ; la propriété du nombre s'ouvre dans un système logique, indépendant de l'arithmétique *per se*.[4]

C'est là où Dedekind s'inscrit lui-même dans le programme kantien : si le nombre en soi est un produit de l'imagination, la théorie lui fournira les concepts synthétiques grâce auxquels nous pouvons soustraire du nombre ce qui se cache derrière sa notion formelle. Non seulement Dedekind construit le concept synthétique du nombre, sur les conseils de Kant, qui accorde aux mathématiques la capacité de rendre sensible toute abstraction [Kant 1956: 287 et passim], il découvre des objets qui correspondent aux chimères arithmétiques. De plus, le nombre dedekindien est transcendant au sens strictement kantien[1], cette transcendance nous incite à démolir les bornes frontières et à entrer dans les sphères tout à fait nouvelles. De Kant à Dedekind puis à Husserl est transmise une idée suivante : le nombre signifie la liberté de créer tout objet dans la sphère des idéalités[2].

---

[1] Notons à cette occasion que la continuité de Dedekind renvoie à l'ouvrage de Guillaume Heytesbury intitulé *Regulae solvendi sophismata* (1335, Règles de résolution des sophismes) où l'auteur, entre autre, analyse le concept de la limite. Pour qu'elle soit établie, il faut appliquer ce que Heytesbury appelle « la capacité active », un agissement venant du sujet vers l'objet. La limite est inséparablement lié au concept de la mesure puisqu'en mesurant une quantité, nous la limitant. La limite c'est toujours la présence d'un sujet même s'il n'apparaît que sous la forme très abstraite (la présence d'un tel sujet ressemble à la présence de l'*ego* dans l'expression « cogito » où l'ego existe dans la conjugaison grammaticale). La limite est une autolimitation car le sujet connaît ce qu'il mesure et donc voit.

[2] Cette diagonale, étant elle-même une mesure, ne peut être décrite par aucun nombre rationnel (le fait qui a eu provoqué l'effroi mystique chez les Grecs).

[3] Bolzano lui aussi travaillait sur ce thème. Quant à savoir si Cauchy avait connaissance des travaux de son collège, cela reste discutable ; cf. supra. De même [Grabiner 1981 : 84 et passim]

[4] Dans une lettre du 27 février 1890 à son critique, l'*Oberlehrer* hambourgeois Hans Keferstein, Dedekind explique clairement son approche ; cf. [Heijenoort 1967: 99-100]. Ibidem il souligne l'apport du *Grundlagen der Arithmetik* de Frege et l'importance de la logique en générale dans la fondation d'une théorie des nombres.

D'un autre coté, ces nouvelles recherches sur la nature du nombre sont inspirées chez Dedekind, comme chez Husserl, par leur maître Weierstraß. Pour celui-ci aucune théorie des nombres ne peut pas être achevée sans la définition logique et rigoureuse des nombres réels faisant objet clé des travaux weierstrassiens. Selon Weierstraß, ces nombres sont une fraction décimale infinie (par exemple, le nombre $\pi$: 3.14159…[3]) qui se manifeste sous la forme d'un nombre irrationnel, c'est-à-dire comme une infinité non définie. La valeur de cette abstraction consiste à montrer que l'infini réside dans l'objet de façon originaire en s'actualisant au cours de notre travail avec ces nombres[4]. Ainsi se résout le problème aristotélicien : l'infini devient actuel dans la conscience du sujet qu'il la soustrait à l'objet par une théorie.

Si chez Weierstraß ces infinités ont plutôt le caractère « lourd », c'est-à-dire qu'elles sont intranscriptibles, Dedekind quant à lui allège les infinités en les coupant en deux classes : gauche et droit entre lesquelles il n'y a pas de places vides. Dans sa lettre du 2 mai 1872 à ce dernier Karl Hattendorff (1834-1882) donne une définition toute en nuance de l'entreprise dedekindienne : « de définir le nombre irrationnel indépendamment de tout appel à l'intuition géométrique » [Dugac 1976 : 36]. Cette nouvelle méthode établie des relations entre les nombres de natures différentes, ce qui permet de traiter tous les nombres comme entités abstraites dont la propriété essentielle consiste en relations entre eux. Philosophiquement Dedekind complète les réflexions cantoriennes concernant les infinités en les faisant plus compactes ; en réalité il n'y a que deux types d'infinités composant le concept du nombre réel. Créer un tel concept semblait à l'époque de trouver la solution de la crise qui est engravée encore par la théorie de Cantor[5]. Il s'agit des objets compacts et séparés déployant l'infini d'eux-mêmes, objets purement mathématiques mais qui appartiennent pleinement à la conscience du sujet. Ce dernier est lui-même le résultat d'une abstraction grâce à laquelle les aspects individuels de chacun ne sont pas pris en compte.

---

[1] Du chapitre de la *dialectique transcendantale*.

[2] C'est la conviction philosophique de Dedekind, qui lui est resté cheviller au corps toute sa vie durant, que les nombres sont créés par nous-mêmes, ce qui atteste sans équivoque notre liberté sans réserve. On se réfère, par ex., à son article de 1854 ou à sa lettre de 1888 à Weber dans laquelle il appelle les êtres humains « la race divine <…> qui possède la puissance créative pour construire des choses mentales » ; cf. [Dedekind 1932 : 428-438, 488-490]

[3] Le savant perse Ghiyāth al-Kāshī (1380-1429) fut le premier à calculer une série longue après la virgule.

[4] José Ferreirós affirme que c'est Weierstraß qui « was able to give a logically rigorous definition of the real numbers, although a rather prolix and complex one. Dedekind and Cantor will simplify the matter by taking the arithmetic of the rational numbers as given » ; cf. [Ferreirós 2007 :126]

[5] En 1882 Cantor a envoyé à Dedekind *Paradoxes de l'infini* de Bolzano. D'après Pierre Dugac, la lecture de ce livre convainc Dedekind de l'existence de l'ensemble infini ; [Dugac 1976 : 87]

La coupure introduit une nouvelle relation, mieux : elle constitue la relation comme telle, la relation *comme* une catégorie mentale qui ne peut exister au niveau des abstractions mois fortes sur lequel la relation ou la limite ne joue qu'un rôle opérationnel. Dans son article de 1887 Dedekind revient à cette problématique, le but est de définir la nature des nombres naturels par lesquels on peut construire des objets mathématiques complexes. Lorsque nous parlons des segments $e$ et $e_1$, nous admettons l'existence des nouveaux objets mais considérons-nous vraiment leur limite[1] (la coupure) comme nouvel objet ? Certes, Dedekind a une tendance à la prendre en considération mais tout le problème consiste en ce qu'une telle limite n'acquiert pas le statut de tel objet en demeurant une « chose en soi » noétique qui confirme de façon allusive la présence du *je* semblable à celui de Cantor. Comment parler autrement des opérations élémentaires sur lesquelles doit reposer la mathématique entière, celles qui doivent nous guider vers la solution cohérente du problème ? Il est peu probable que la phénoménologie transcendantale de Husserl aurait été la même si le maître n'avait réfléchi systématiquement aux chemins menant vers la solution de cette crise et s'il n'avait pris en compte l'expérience de ses pères. Même dans ses travaux tardifs, tel *Méditations cartésiennes* (1929) en expliquant la nature de la subjectivité transcendentale, Husserl fait allusion à Dedekind : « [cette dernière] n'est pas un chaos des expériences intentionnelles ; elle n'est pas non plus un chaos des types constitutifs organisés chacun à sa manière par leur relation avec des objets intentionnels. En d'autres termes : le tout (*die Allheit*), qui – pour dire transcendentalement – est pour moi comme un *ego* transcendantal composé des objets et des types d'objet imaginables (*erdenklichen*), n'est guère un chaos <…> Cela nous fait prévoir une synthèse constitutive universelle dans laquelle toutes les synthèses jouant de concert sont mises en ordre (*Weise*) d'une manière certaine… » [Hua I : 90]. L'*ego* transcendantal est en effet une unité bien ordonnée où les objectités se distinguent par ce que leurs relations à la conscience sont définies.

On voit bien que

α) la construction du sujet en mathématiques du XIX$^{\text{ème}}$ a été un processus inévitable qui concernait presque tous les grands théories de l'époque. La pensée cherche la solution de l'impasse ;

β) la naissance d'un tel sujet change une fois pour toujours le caractère de la science exacte.

La clarté n'est pas plus la propriété innée des mathématiques (selon Platon) mais un produit de l'amélioration d'une théorie qui, dans ce processus de la clarification de ses fondements, constitue le sujet dans elle-même. A strictement parler, ce sujet n'est pas *cogito* cartésien au sens qu'il n'a pas *sa* pensée justifiant son existence (c'est ici le sujet mathématique se diffère de son concept philosophique), bref il ne porte guère le caractère transcendantal ; dépourvu des soins ontologiques, ce sujet sert à éclaircir non pas le rapport de la vérité à la réalité mais la position de la vérité dans une théorie. On appelle le *vrai* non ce qui correspond au réel, au champ du vécu empirique mais ce qui permet d'édifier la théorie sans équivoque [cf. Borel 1914]. Le sujet mathématique est celui qui fait l'abstraction des données de la théorie et c'est dans ce processus la théorie se précise en s'approchant de son état de la puissance irrépréhensible. La présence du sujet est incluse dans la théorie qui ne la découvre qu'à travers ses propres limitations. Nicolas Lusin donne un bon exemple : si l'on ne possède pas une idée nette et claire de la suite illimitée des nombres entiers, il ne reste qu'à avancer une autre définition qui procure cette idée en enlevant ainsi les limites précédentes [Lusin 1930 : 35]. On peut donc dire à juste titre que ces dernières jouent le rôle capital donnant la possibilité à la théorie de surmonter ses bornes.

[1] Il est intéressant de rappeler ici qu'après Leibniz, Euler et D'Alembert, qui dans son *Encyclopédie* a travaillé beaucoup sur la nouvelle méthode de l'algèbre, certains mathématiciens (S.-F. Lacroix, par ex.) considèrent le concept de la limite comme la véritable métaphysique du calcul.

Ne peut-on y voir l'illustration du credo philosophique de Dedekind formulé dans la note 66 du *Was sind und was sollen die Zahlen ?* : « Mon monde mental, c'est-à-dire la totalité *S* de toutes les choses, capable être l'objet de ma pensée, est infini » [Dedekind 1932: 357][1]. L'expression « l'objet de ma pensée » laisse entendre une conclusion importante : le sujet considère lui-même comme système actuellement infini qui n'a pas des limites sauf, pour ainsi dire, le nombre des éléments entrant en lui. Dans la *Philosophie première* Husserl énoncera la même pensée mais dans le cadre de son programme transcendantal : « <…> *le monde* tel qu'il est en lui-même et dans sa vérité logique n'est en dernière analyse qu'une *idée située à l'infini* puisant son sens intentionnel dans l'actualité de la vie de la conscience »[2]. L'infini entre dans le sujet en tant que sa structure constituante indiquant la direction des actes de la conscience vers l'appréhension d'une généralité conceptuelle. Si l'on modifie la limite d'un objet, il se transforme en un autre objet, la chose n'est jamais la même avec la limite modifiée. L'infini n'est pas l'objet au sens strict car aucune modification de ses limites n'affecte point sa nature. L'infini ne se change ni par l'extraction de sa limite quelconque ni au cours du temps en représentant lui-même le temps mental pur. Ainsi le temps du sujet mathématique qui pense non seulement d'un élément concret, d'un nombre ou d'un point, il commence à les voir au total, tous les éléments comme une unité ; « l'objet de ma pensée » signifie cette pensée même comme objet. Les nombres naturels, étant mes propres créatures, n'ont pas d'haubans ontologiques[3], ils existent éidétiquement, sans être. C'est pourquoi le mérite philosophique de Dedekind consiste à définir clairement les procédés de la pensée mathématique pour qu'elle soit capable de créer des nouveaux objets ou de les annuler. Quant à l'existence de l'ensemble

---

[1] [Meine Gedankenwelt, d. h. die Gesamtheit *S* aller Dinge, welche Gegenstand meines Denkens sein können, ist unendlich]

Hilbert remarque à ce propos : « l'erreur classique [de Dedekind] a consisté à prendre pour point de départ le système de tous les objets. Pour aussi brillante et séduisante qu'eût été son idée de fonder le nombre fini sur l'infini, le caractère impraticable de cette voie ne fait aujourd'hui de doute pour personne… » ; cf. [Dedekind 1932 : 162]. David C. McCarty écrit que, du point de vue ensembliste, le concept du *Gedankenwelt* signifie une classe entière (non pas l'ensemble < sic ! >) ou une hiérarchie cumulative de toutes les classes. En outre, *Gedankenwelt* (ainsi que l'*eigenes Ich*, expression figurant dans le texte) se réfère à l'*ego* pensant ou à l'existence de Dedekind lui-même. McCarty a raison de considérer ces réflexions comme étrangères aux mathématiques. Remarquons que § 66 fait écho à la section 13 de *Paradoxes de l'infini* (1851) de Bolzano ; pourtant, Dedekind affirme qu'il n'avait pas connaissance de cet ouvrage au moment où il a écrit *Was sind und was sollen die Zahlen ?* Cf. [McCarty 1995: 57-59]

[2] Philosophie première, t. I, p. 350.

[3] Howard Stein considère la théorie/philosophie dedekindienne du nombre comme *anti-ontologique* libérant le concept du nombre de ses entités mondaines ; cf. [Stein 1989: 247]

infini, elle doit être démontrée logiquement pour éviter les doutes concernant son caractère non-contradictoire. Ce procédé se réalise en liberté en ne signifiant rien d'autre que la capacité humaine de s'abstraire et de généraliser nos propres pensés. Ainsi l'ensemble devient l' « objet singulier » dans la théorie dedekindienne ; Bourbaki note que c'est lui qui a introduit la notion de la « chaîne » [Bourbaki 2007 : 44] permettant de percevoir les éléments d'un ensemble comme le nombre infini des limites.

La construction du sujet faisait partir la mathématique mais sans résoudre le problème de l'origine de notre savoir. Par exemple, d'où vient notre certitude qu'on peut couper la ligne droite *ad infinitum* ? La réponse reste kantienne : à cause de notre intuition du temps qui nous apparaît infini. Si la série des nombres représente le temps pur, alors ce temps a seulement des limites mentales qu'établie notre conscience. Un autre exemple : les infinitésimaux qui se produisent dans le *temps* de la décroissance infinie ; à vrai dire, ils n'existent qu'intuitivement apparaissant toujours à un moment précis. Le moment qui trahit l'existence de l'intuition phénoménologique, c'est-à-dire du temps de la possibilité pure qui nous donne la liberté de construire toutes abstractions. En outre, le même exemple montre la limite de l'infini. Selon la règle algébrique, l'infinitésimal doit se décroitre – mais rester toujours! – au cours de ce procédé inférieur à tout nombre positif qui le limite. Paradoxalement, le sujet est constitué dans la mathématique justement aux niveaux les plus éloignés de l'expérience quotidienne, il surgit de cette pure possibilité. L'infinitésimal est simultanément fort abstrait et subjectif car l'idée de la décroissance infinie, qui est d'ailleurs parente à la coupure dedekindienne, admet l'intuition subjective ou, si l'on préfère, la liberté phénoménologique qui permet de faire de la possibilité pure l'objet du savoir. L'idée de l'infinité, fût-elle mathématique ou scolastique ou autre, implique nécessairement le sujet – le *cogito* – dont le processus de pensée se déroule dans une succession des moments.

Phénoménologiquement parlant, cette idée est innée dans l'*ego* pensant comme sa condition qui lui donne à construire lui-même ; l'infini est l'horizon noétique dans lequel le sujet transcende sa conscience individuelle en parvenant à la position d'un inter-sujet. Ainsi, si en mathématiques le sujet se cristallise bon gré mal gré par le besoin de trouver un point de vue « extérieur » par rapport à ses constructions abstraites, le sujet phénoménologique doit en revanche neutraliser toute l'extérieur, ou plus précisément, faire s'annihiler l'intériorité et l'extériorité pour priver la conscience de son architecture empirique. Enfin, ces deux sujets se trouvent liés mais leur divergence semble être prédestinée. Le conflit de Husserl avec la pensée mathématique est ici en germe ; aussi bien sa critique de Schröder que sa vision des problèmes généraux du savoir trahissent en lui la position qu'il occupe peu à peu mais avec certitude au cours de son travail sur la phénoménologie et sur son autoconstitution en tant que sujet qui ne connaît que l'espace de sa conscience. Quoi qu'il en soit, le résultat le plus intriguant pour Husserl est celui de Dedekind. Chez ce dernier le sujet est privé de toute sorte de psychologisme, il pense non pas de lui-même mais « hors-lui-même » ; l'objet de sa pensée n'est que la réalisation de l'intuition à toute moment du temps (cet antipsychologisme Frege reprendra pour son thème principal et pour sa critique de la métaphysique, y compris les premiers travaux de Husserl). Dans ce temps infini l'intuition trouve sa forme parfaite. Du surcroit, le sujet créé devient immortel car il se pense de lui-même dans le temps qui n'a pas de fin. Quant à la phénoménologie transcendantale, elle recrée ce sujet de telle façon qu'il puisse reconstitue son origine.

## CHAPITRE 6. L'AXIOMATIQUE : METHODE DE LA CONSCIENCE

La réduction de toute la mathématique à ses propriétés fondamentales où n'interviennent que les opérations élémentaires se nomme « l'axiomatisation ». Peu de temps après son arrivée à l'Université de Göttingen, Hilbert annonce un nouveau cours de géométrie. La matière de ce cours composera un des ses livres majeurs, *Fondements de la géométrie*[1], où il élabore sa méthode abstraite permettant d'axiomatiser la géométrie. En bref, l'idée est la suivante : faire abstraction de la nature des objets géométriques (le point, la droite ou le plan) et se contenter de poser entre eux des relations dont les propriétés sont explicitées par des axiomes. Peu de temps avant la seconde édition de son ouvrage[2], Hilbert explique qu' « étant engagé dans la recherche des fondements de la science, nous devons établir un système d'axiomes qui contient une description exacte et complète des relations entre les idées élémentaires de la science. Ainsi établis, les axiomes sont en même temps les définitions de ces idées élémentaires <…> »[3] La méthode a son histoire, écrite entre autres par Frege, Cantor, Dedekind et Pasch[4] ; ce dernier publie *Leçons sur la nouvelle géométrie*[5], où il propose de reconstituer la géométrie euclidienne en termes plus primitifs et plus précis en prenant pour base une liste d'axiomes. Pasch recommande dans son ouvrage de ne pas interpréter les concepts géométriques au sens physique ; tout ce que nous devons prendre en considération, ce sont les axiomes propres à la géométrie. Les réalités physique comme psychique ne jouent aucun rôle lorsque l'on travaille avec les pures constructions de l'esprit.

A la différence du jeune Husserl, Pasch est loin d'être un platonicien convaincu ; il voit les concepts de la géométrie comme construits sur des données empiriques (*empirische Material*) et les axiomes ne sont que l'extension (au sens grassmannien) de faits observables. Ensuite, les concepts géométriques constituent un groupe spécial de concepts qui ne servent qu'à décrire le monde extérieur. Mais ces concepts eux-mêmes ne font pas partie de ce monde qui demeure l'objet à décrire. Dans le monde, l'« objet » (*Körper*) signifie la présence réelle d'un objet que nous percevons par nos organes de perception. La manière dont nous le percevons, constate Pasch, est fondamentalement différente de sa véritable structure physique : la description de cet objet doit donc être basée sur un système d'axiomes infalsifiables. D'où il découle que les objets géométriques n'ont aucun rapport ni avec le temps ni avec l'existence. Or, leur existence n'étant que conceptuelle, cela fait de la géométrie « rien d'autre qu'une partie des sciences naturelles » [Pasch 1882: 3].

---

[1] Titre original : *Grundlagen der Geometrie* (1899).
[2] Parue augmentée en 1903.
[3] « Mathematische Probleme » (1901), [Hilbert 1935: 301]
[4] Moritz Pasch (1843-1930), mathématicien allemand, il suivit les cours de Kronecker et de Weierstraß à Berlin et effectua ensuite toute sa carrière universitaire à Giessen.
[5] Titre original : *Vorlesungen über neuere Geometrie* (1882). Notons au passage que l'influence de ce livre sur Peano, Hilbert et d'autres a été considérable.

C'est pourquoi il faut corriger Euclide[1], qui admet parfois des propositions évasives, et discerner les concepts primitifs entachés de l'incertitude du monde quotidien. La mathématique se fonde sur les relations pures établies par elle-même entre objets *mathématiques*. Cette approche s'appelle « formaliste » ; Jan Brouwer, intuitionniste et critique-en-chef du formalisme, exprime son principe avec subtilité : « Le formaliste soutient que la raison humaine n'a pas à sa disposition d'images exactes des lignes droites ou des nombres supérieurs à dix, par exemple <…> [Cependant] pour le formaliste, l'exactitude mathématique réside dans le développement de séries de relations et est indépendante de la signification que l'on pourrait vouloir donner à ces relations ou aux entités qu'elles relient » [Brouwer 1913: 83][2]. Autrement dit, la mathématique ne reconnaît aucun regard extérieur bien qu'elle connaisse l'existence du monde physique. L'idée de l'extension, à supposer que tous les auteurs n'aient pas lu Grassmann, a occupé une place centrale dans les esprits de sensibilité formaliste. L'axiome est le produit du procédé d'extension auquel Husserl donnera un sens plus philosophique en le considérant comme abstraction.

Tout comme Pasch, le fondateur de la phénoménologie ne nie pas le monde externe qui nous inonde des données par lesquelles nous entrons dans la sphère idéale, mais il constate que ces deux mondes ne peuvent jamais se confondre : « la géométrie pure est une science *a priori*. La possibilité qu'elle trouve une application dépend du concept empirique de l'espace et de sa capacité à correspondre ou à être soumis à son concept pur de l'espace. Cela s'avère impossible. Parmi les objets de notre expérience il n'y a pas de point au sens strict du terme, c'est-à-dire quelque chose de spatialement invisible »[3].

Or, l'idée de Pasch selon laquelle l'axiomatisation enfante l'idéalisation est étrangère à Husserl. Ces deux processus se déroulent simultanément à partir du moment où l'intuition entreprend de s'ajuster au monde sensible[4]. Un tel ajustement est naturel car l'être humain a toujours besoin de géométriser le monde alentour. « La construction théorique des idéalités pures représentant la géométrie (il s'agit en général des figures de mathématique (*Gestalten-Mathematik*) <…> des idéalités et des constructions idéales accomplies guident vers l'art de la mesure pratique, vers le monde corporel dont l'objectivation se réalise dans des sphères d'intérêt limité » [Hua VI : 35].

Voilà la conclusion de réflexions husserliennes, qui ne changent pas d'ailleurs sa position de l'époque de Göttingen. L'axiomatique n'est que l'outil pour extraire des données observables leurs propriétés essentielles et y donner une forme la plus primitive possible, c'est-à-dire non-contradictoire. « La mathématique formelle est un instrument concret pour les découvertes mathématiques ; <…> la nouvelle mathématique formelle [se distingue] des veilles méthodes quantitatives… [la première] crée une généralité et une puissance incomparables en l'art mathématique » [Hua XII : 432]. Dans *Logique formelle et logique transcendantale* (1929) Husserl revient sur ce point : « le système axiomatique défini est ainsi établi ; de chaque concept apparaissant (en fait des formes conceptuelles) l'on peut construire soit une proposition (forme propositionnelle) vraie – une suite analytique (purement déductive) de l'axiome –, soit une proposition fausse, notamment une contradiction analytique : *tertium non datur* » [Hua XVII : 100]. Cette affirmation est proche de Hilbert, qui distingue clairement la *représentation* et l'*archétype idéal* : la première est une démarche qui concrétise la chose en trouvant pour elle une place parfois éventuelle dans la réalité ; quant à l'archétype, il ne réside dans notre connaissance que sous une forme inconnaissable. Bref, l'axiomatique vise à construire des ponts entre deux mondes – naturel et idéal [1] – et à établir d'exactes corrélations entre les démarches de la pensée et la structure de l'objet. « La figure représentée est seulement une admission… » [Pasch 1882: 45], alors que l'archétype est sa condition. Néanmoins, la perception des vérités géométriques s'effectue à partir de la chose représentée vers son idéal, c'est cette première qui constitue une donnée perceptible privilégiée, apanage exclusif de la vue, et qui soutient « l'intuition d'essence ».

---

[1] Pasch critique Euclide, entre autres, pour l'incomplétude (*Unvollkommentheit*) de ses propositions et de son axiomatique, en notant que ces erreurs principales n'ont pas encore été rectifiées ; cf. [Pasch 1882 : 45]. Ulrich Majer remarque que Hilbert reprend plus tard cette critique en critiquant Pasch lui-même ; pour ce premier, Pasch tout comme Euclide demeure dans le cadre de l'axiomatique classique : « We can say that Hilbert amended the classical axiomatic approach of Euclide and Pasch with a systematic inquiry of metalogical relations among axioms and axiom systems. The principal means of this inquiry is the axiomatic method in connection with a certain type of model-theoretic consideration… » ; [Majer 2004: 104-105]. Quant à Husserl, il travaille durant sa période de Göttingen sur le problème de la complétude en développant son propre point de vue. Mais sa critique du même Pasch est moins radicale que celle de Hilbert ; Husserl en reste *grosso modo* aux fondamentaux dedekindiens. Pour les details, voir [Hartimo 2007: 281-310]

[2] Un bon résumé des idées de cet auteur on peut trouver dans [Dalen 1999]

[3] Hua, XXI, p. 296. Lothar Eley considère la conceptualisation du savoir apriorique comme une des causes capitales de la crise dans la philosophie moderne. Ajoutons que cela est directement lié au problème des fondements en mathématiques. Cf. [Eley 1962 : 26]

[4]Une autre explication est donnée par [Richir 2002 : chap. 2]

Ainsi se dégage l'axiomatique, enlevant les contradictions entre approches empirique et transcendantale. Lorsqu'il s'agit de la géométrie comme science, gardons-nous de choisir une position préalable qui dominerait ensuite règles notre moyen-même de voir les choses. Toute figure concrète (représentée) [2] surgit de formelles indifférentes à leur contenu ; « en réalité, pour que la géométrie soit véritablement déductive, il faut que le processus de déduction soit partout indépendant du sens (*Sinn*) des concepts géométriques <…> Seules les relations entre les concepts géométriques, telles qu'elles sont explicitées dans les propositions et les définitions utilisées doivent être prises en considération. Au cours d'une déduction, il est certes permis et utile de conserver à l'esprit la signification (*Bedeutung*) des concepts géométriques employés, mais ce n'est en aucune façon nécessaire » [Pasch 1882: 98]. Seules importent donc les relations pures engendrées par les axiomes. On emboîte ici le pas à Dedekind qui, comme nous l'avons vu, constitue le savoir mathématique à partir d'un semblable mécanisme relationnel en mettant de coté l'aspect du contenu. Pasch, quant à lui, s'attache à l'applicabilité des concepts géométriques à la réalité empirique. Malgré son caractère abstrait, l'axiome n'est pas un noème husserlien et le chemin de l'un à l'autre passe non seulement par la suspension du monde de l'expérience (d'Euclide), mais aussi par la modification de la conscience du sujet ; selon Pasch, l'axiomatique est le moyen de gérer le monde empirique comme système fondé sur des principes indépendants de l'individu. Autrement dit, l'axiome ne modifie guère la conscience subjective, il ne fait que livrer au sujet le savoir objectif.

[1] Ne s'agit-il pas ici d'une allusion de Pasch à l'idée cartésienne de « la Géométrie naturelle », dont Dieu seul comprend l'essence ? Si les archétypes géométriques existent, alors ils ne peuvent être les objets de notre connaissance car tout objet doit avoir un caractère représentatif.

[2] Il suffit de préciser que le caractère concret d'une figure signifie ici non pas une chose particulière, individuelle ni *a fortiori* unique ; Pasch parle d'une « copie » mentale dont le nombre est infini. La figure concrète est élaborée sur un modèle déterminant l'intégralité de la classe de ces figures. Cet engendrement est possible dans un espace doté de paramètres tels que ce modèle ne subit pas de modifications structurelles.

Pour ce faire, deux conditions sont à respecter : premièrement, définir, c'est-à-dire expliciter, toutes les propositions (par exemple géométriques), y compris celles apparaissant comme les plus évidentes : « il n'est de supposition qui ne soit énoncée », répétait Pasch et ceux qui partageaient sa méthode ; deuxièmement, il ne faut pas qu'une épreuve géométrique quelle qu'elle soit repose sur une intuition sensible, voire sur l'esthétique transcendantale (cette démarche « anti-kantienne » est déjà entreprise par Bolzano, qui note que les jugements/propositions fondées uniquement sur l'intuition offusquent leur origine). Ainsi, Pasch insiste : toute proposition doit être déduite soit d'un axiome, soit d'une proposition antérieure. Autrement dit, les vérités dites mathématiques doivent par nature constituer un système axiomatique permettant de mettre à l'épreuve toute intuition géométrique et de la convertir en une catégorie purement formelle. Pasch propose une nouvelle grammaire descriptive qui concorderait avec celle de la langue allemande ; bref, la géométrie n'est autre qu'une grammaire des intuitions d'espace[1].

Par la suite, ces réflexions se sont avérées capitales dans ce que l'on appelle « l'axiomatisation » des mathématiques et dont l'apôtre fut Hilbert, réformateur et camarade de Husserl à Göttingen. Hilbert va plus loin en disant que le jugement mathématique doit absolument avoir la propriété de non-contradiction apriorique [Contro 1976]. En outre, le mathématicien doit savoir que la vérité existe indépendamment de toute théorie et il doit voir cette vérité avant-même qu'elle ne soit formulée dans une théorie. Le mathématicien s'inspire de la vérité, imaginait Hilbert, comme un poète de sa muse.

Si Pasch répond à cette « crise de l'évidence »[2] par une tentative de créer une grammaire des axiomes en les prenant comme seules conditions possibles de toute évidence et se trouvant par là proche de Bolzano, Hilbert revient lui à Kant[3] : « le but de toute science consiste avant tout à établir un réseau de

---

[1] A la même époque (dans les années 1880) en Allemagne sont élaborées de nouvelles théories de la grammaire, tant allemande qu'indo-européenne. La théorie la plus puissante est créée par un groupe de linguistes connus sous le nom « néogrammairiens » (H. Paul, O. Behaghel, K. Brugmann et d'autres). Leur idée principale consiste à décrire la langue non comme intuition réalisée dans l'expérience mais comme un système de règles auxquelles obéit la conscience linguistique des individus.

[2] Gonseth remarque avec raison que cette crise « est ouverte depuis la découverte des géométries non euclidiennes <…> comment expliquer que le monde mathématique s'en montre si peu inquiet, et qu'on continue avec sérénité à invoquer l'évidence comme le dernier fondement du vrai » ; cf. [Gonseth 1939 : 45 ; Torretti 1978]

[3] C'est par Leonard Nelson (1882-1927), fondateur de l'école néo-friesienne, que Hilbert a été introduit à la philosophie critique. Leur rencontre a lieu à Göttingen lorsque Nelson, après avoir soutenu son doctorat à Berlin, attend son habilitation qui lui sera du reste refusée par les membres du jury composé de philosophes et de mathématiciens (rappelons qu'à Göttingen les mathématiques et la philosophie appartenaient tous les deux à la Faculté de la philosophie). Bien que Husserl y habite lui aussi au même moment, Nelson choisit de garder ses distances vis-à-vis de ce dernier. Il reproche à Husserl (et à Frege) d'être l' « adepte de préjugés dogmatiques » dont la grande erreur consiste à ne pas discerner les différents sens du « fonder » (*Begründung*) de la logique sur la psychologie. Cf. [Nelson 1973 : 167-169]. A propos : [Anderson 2005]

Kurt Grelling (1886-1942), logicien ayant travaillé avec Nelson sur la solution du paradoxe de Russell, est aussi le témoin des relations singulières entretenues par les habitants concepts (*Fachwerk von Begriffen*) dont la conceptualité nous est donnée par notre *intuition* et notre *expérience*. Idéalement, tous les phénomènes d'un domaine apparaîtront comme partie de ce réseau et tous les théorèmes, qui peuvent être dérivés des axiomes, y trouveront leur expression »[1].

L'évidence précède la connaissance que nous en avons ; autrement dit, une science exacte telle que la géométrie est engendrée par l'objectivité en-dehors de toute connaissance individuelle et donc exempte de contradictions. Mieux, cette évidence est non pas le produit de la raison mais *sa* nature, dévoilée au cours de la construction de l'exactitude – de la *non-contradiction* – géométrique. Avec Ferdinand Gonseth, admettons que « le postulat des parallèles puisse être remplacé par un postulat qui le contredise, et que ces deux postulats puissent subsister l'un à coté de l'autre, il est clair que l'un et l'autre perdraient de leur évidence… » [Gonseth 1939 : 49].

Nul doute que l'axiomatique a été inventée pour mettre fin à cette longue période critique des mathématiques par la conceptualisation de leurs fondements ; Hilbert avoue que « le but de fonder les mathématiques dans la sécurité [était] aussi le [sien] ; [il désirait] rétablir les mathématiques dans leur ancienne réputation de vérité inattaquable… » [Hilbert 1935 : 160]. Dans son texte *Sur l'infini* (1925) il décrit sa position philosophique fondamentale : les mathématiques présentent le modèle de toutes les sciences et de la pensée en général [Hilbert 1926: 171] ; ce crédo hilbertien nous montre une fois encore les préférences de son auteur. Le fait que l'axiomatique trouve ses racines philosophiques dans la pensée kantienne ne semble rien devoir au hasard, car la philosophie transcendantale exige une évidence conceptuelle structurée par deux types de jugement : analytique et synthétique. Henri Poincaré, qui ne se rapporte qu'à l'axiomatisation de la géométrie et, comme on le sait, n'a jamais renoncé à ses ambitions d'établir un lien entre l'intuition et l'évidence, insiste sur la nécessité de jugements synthétiques *a priori* [2] ; à vrai dire, cela signifie non pas *zurück zu Kant*,[3] mais *en avant avec* Kant.

Cette question possède un versant gnoséologique : « la méthode axiomatique est en effet et reste le moyen convenable et indispensable à toute recherche exacte dans tout domaine qui soit. Elle est logiquement incontestable et en même temps féconde <…> procéder axiomatiquement, c'est, dans ce sens, *avoir conscience de sa pensée* ; alors qu'avant, quand on n'avait pas la méthode axiomatique, il arrivait que l'on crût naïvement à certaines relations comme à

---

de ce paradis scientifique ; voir [Luchins, Luchins 2000]. Juste après l'échec de l'habilitation de Nelson, Hilbert reste à ses cotés ; les deux hommes apprécient de partager des débats intellectuels qui s'avèrent extrêmement féconds pour l'un comme pour l'autre. Hilbert maîtrise le projet critique de Kant ; Nelson devient partisan de la méthode axiomatique dont il trouve des parallèles dans la philosophie friesienne. Cf. [Reid 1970 : 122 et passim] des dogmes[1] <…> et plus loin : « l'axiomatisation nous contraint donc à prendre position vis-à-vis de ce difficile problème de la théorie de la connaissance » [Hilbert 1935: 161 nous soulignons]. Pour nuancer cette affirmation, Hilbert demande à ce que les objets à définir par les axiomes aient un caractère purement apriorique, que ce soient ainsi des éléments indéterminés dont seuls les axiomes montrent leurs propriétés. En d'autres termes, les axiomes ne peuvent que fixer les conditions (transcendantales) selon lesquelles nous prenons connaissance de ces objets. D'où ressort l'idée selon laquelle les axiomes ne nous disent rien des objets eux-mêmes mais seulement de la logique pour les connaître. Il s'agit donc d'une conséquence fort importante de l'approche hilbertienne : la logique du savoir portant sur des objets mathématiques est plus fondamentale que leurs propriétés formelles ; de plus, ce raisonnement annule la différence, chère à Poincaré, entre la géométrie et l'arithmétique, au sens que les deux sciences se constituent comme relations pures entre des objets idéaux (rappelons que Pasch défend une telle approche [Pasch 1930 : 2]).

[1] [Hilbert 1905]. Ce retour à Kant n'empêche pas Hilbert de noter que, pour atteindre une conscience mathématique apriorique (ou *Denkapriori*, come le disait Dingler), il faut débarrasser du « mâchefer anthropologique » kantien.

[2] Au début des années 1880, Poincaré travaille sur des problèmes complexes de la géométrie, notamment sur les groupes de transformation discontinus ; pour les détails voir [Kleiner 2007, esp. Chap. 2]. Avec Sophus Lie (1842-1899) et Félix Klein (1849-1925), dont les recherches sont dirigées dans la même direction, Poincaré est tout proche du concept des groupes abstraits, ce qui a marqué profondément la pensée mathématique.

[3] Cf. [Peckhaus 1990]

Au Congrès international de mathématiques qui s'est tenu en 1900 à Paris, Hilbert résume : « la coïncidence entre la pensée géométrique et la pensée arithmétique se révèle encore en ceci : dans les recherches arithmétiques, de même que dans les considérations géométriques, nous ne remontons pas à chaque instant la chaîne des déductions jusqu'aux axiomes » ; plus loin : « <…> la non-contradiction des axiomes géométriques est ramenée à la démonstration de la non-contradiction des axiomes de l'arithmétique » [Hilbert 1990 : 9, 15 ; de même Bernays 1922]. En 1891, après un séminaire de Hermann Wiener (1857-1939) consacré aux fondements de la géométrie, Hilbert déclare qu'il faut que des concepts comme le point, la ligne ou la platitude puissent être remplacés dans toute phrase géométrique par table, chaise, pot à bière. Par conséquent, point, ligne et platitude ne sont que « trois systèmes des chose » qui satisfont les axiomes de la géométrie[2]. Il s'ensuit que les axiomes sont des énoncés abstraits établissant un groupe d'opérations avec les objets ; si l'on dit qu'il n'existe qu'une ligne droite qui passe par deux points quelconques, il s'agit d'un axiome fixant la propriété d' « être droit » et, phénoménologiquement parlant, d'une relation entre deux idéalités. Selon le maître de Göttingen, « cet énoncé se réduit essentiellement à ce théorème d'Euclide qui veut que, dans un triangle, la somme de deux côtés soit toujours plus grande que le troisième ; il est facile de voir qu'il ne s'agit dans ce théorème que de concepts élémentaires, c'est-à-dire dérivant immédiatement d'axiomes… » [Hilbert 1990 : 19].

Une démonstration sur la surface euclidienne prouve le caractère non-contradictoire de cet axiome, mais il existe un autre niveau du problème dont personne n'a vu toute la gravité à cette époque : éviter les contradictions dans la démonstration-même. « Les axiomes et les théorèmes démontrables, c'est-à-dire les formules qui naissent au cours de cette interaction [entre eux], sont les images des pensées qui constituent la méthode habituelle de la mathématique traditionnelle, mais ce ne sont pas eux qui se révèlent être des vérités absolues » [Hilbert 1935: 180]. En revanche, ces dernières se décèlent en tant qu'évidences ou, en termes hilbertiens, *Einsichten*[1] – vues absolues qui ne peuvent contenir de quelconques contradictions. Là comme ailleurs, Hilbert opte pour une solution kantienne : créer un système de règles *a priori* capable d'écarter toute démonstration contradictoire. Ces règles doivent être strictement formelles et basées non sur l'intuition – quand bien même elle serait prouvée dans la vie quotidienne –, mais sur la totalité du savoir mathématique résumé ainsi dans sa forme la plus abstraite.

[1] Plus tard nous verrons que dans l'*Einleitung in die Logik* 1906/7 (Hua XXIV) Husserl parle dans les mêmes termes de la méthode phénoménologique par rapport au dogmatisme psychologique.

[2] Cette position de Hilbert récite un théorème de Kant formulé au chapitre III des antinomies de la raison : « La nature de la raison humaine est architectonique, c'est-à-dire qu'elle considère toutes les cognitions comme appartenant à un système possible et c'est pourquoi la raison ne tolère que des principes qui au moins ne contredisent pas ce système… » ; [Kant 1956: 479]

Précisons : la *règle* chez Hilbert est une abstraction au même titre, par exemple, que la réduction phénoménologique ayant pour but de mettre à nu l'objet et de voir son essence. Mais ni la règle ni l'axiome ne sont dirigés vers l'objet : ils le contiennent. En disant qu'il n'existe qu'une ligne droite qui passe par deux points quelconques, nous voyons tout de suite le contenu de cet axiome, exposé par sa syntaxe[2]. L'axiome est construit comme un jugement analytique dont le contenu ne signifie qu'une propriété naturelle de ses éléments. Sur ce plan, la réduction phénoménologique est un procédé subjectif dont le contenu se constitue par la réduction en cours.

Dans un séminaire de 1922/23 intitulé *Introduction à la philosophie*,[3] Husserl parle de la vision pure des choses, passant toujours par la réduction phénoménologique ; celle-ci exige de nous de sacrifier le monde avec toutes les âmes qu'il peut compter, y compris la nôtre [Hua XXXV : 72]. Ce sacrifice nous mène vers l'*ego cogito*, l'acte générateur de la conscience. C'est par la réduction que l'on crée l'évidence ou la vision à partir de laquelle se constitue la possibilité-même de la vérité absolue. Ainsi apparaît la zone spécifique que Husserl appelle « une région de l'être en soi » (*ein Seinsgebiet in sich*) [Hua XXXV : 73] – ressemblant au système des axiomes –, où tout objet est saisi dans sa forme idéalisée. La réduction phénoménologique décharge la conscience de toute propriété pour la guider hors de l'analycité. Il s'agit du contenu qui n'existe pas en tant qu'objet mais seulement comme état de la conscience. Une philosophie, pour devenir *la* philosophie, doit trouver en elle-même une évidence absolue, son commencement. C'est pourquoi l'*Einsicht* hilbertienne et l'*Evidenz* husserlien [Hua XXXV : § 10-11] sont deux signes indiquant la même direction. Il est à relever qu'un dialogue implicite entre Husserl et Hilbert s'est poursuivi bien après le séjour de ce premier à Göttingen ; par exemple, au § 20 des *Ideen I* Husserl appelle précisément les axiomes « des évidences générales [qui] expriment des faits d'expérience <…> » [Hua III/1 : 45], allusion sans équivoque à Hilbert ! Ce dernier cite à son tour *L'idée de la phénoménologie* (1907)[1] en disant que « <…> l'aspect phénoménologique de la théorie [de l'hydrodynamique] concerne au niveau le plus *immédiat* les phénomènes d'émission et d'absorption <…> »[2].

---

[1] « Als die absoluten Wahrheiten sind vielmehr die Einsichten anzusehen » ; *Ibid*.

[2] Plus tard, dans les années 1930, la syntaxe des propositions logiques devient l'un des thèmes principaux de la philosophie allemande (autrichienne) puis anglo-saxonne. Rudolph Carnap va jusqu'à circonscrire les problèmes de la philosophie aux constructions syntactiques, persuadé qu'il est que le modèle formel du langage est le seul chemin pour échapper aux ténèbres métaphysiques. « Philosophy, écrit Carnap, is to be replaced by the logic of science – that is to say, by the logical analysis of the concepts and sentences of the sciences, for the logic of science is no nothing other than the logical syntax of the language of science » [Carnap 2000: xiii]. Le projet de Carnap repose sur les mathématiques du XIXème siècle et sur les travaux hilbertiens qu'il ne cesse de citer dans son livre. En outre, Carnap répond au théorème d'incomplétude de Kurt Gödel (énoncé pour la première fois en 1931) qui a bouleversé la logique mathématique. Quant à Hilbert, il ne considère pas son axiomatique comme davantage qu'une critique des mathématiques pures. Cf. « Axiomatische Denken » (1918) [Hilbert 1935: 147-156]. De même : [Giaquinto 1983]

[3] Dans les mêmes années Hilbert publie deux textes : *Nouvelle fondation des mathématiques. Première communication* (1922) [Neubegründung der Mathematik. Erste Mitteilung] et *Les fondements logiques des mathématiques* (1923) [Die logischen Grundlagen der Mathematik], où il explique une fois encore les principes de l'axiomatique. Le premier texte décrit l'axiomatique comme une théorie de la conscience.

A l'instar de Husserl, Hilbert met le monde extérieur hors-circuit pour neutraliser son objectivité expérimentée. L'axiome n'est pas basé sur l'objectivité, il la crée en s'appuyant sur l'intuition subjective. Le monde subit l'ἐποχή ; il n'est ni la base du raisonnement, ni le critère pour démontrer ses résultats mais une des innombrables possibilités de l'interprétation des objets idéaux. La réalité pour le mathématicien, c'est la classe des énoncés mathématiques – les opérations – effectués dans des formules abstraites. Soulignons de nouveau : les axiomes hilbertiens sont les limites du monde où réside le mathématicien. Cela veut dire que pour être valable, tout énoncé mathématique ne doit pas dépasser les limites de ce monde axiomatique (dans le cas inverse l'énoncé perd sa propre objectivité). C'est pourquoi les axiomes sont non seulement les garants de la non-contradiction, mais ils fournissent aussi le sujet d'une langue *l* dans laquelle celui-ci peut s'exprimer sans contradictions. Découvert par Dedekind et Cantor, le sujet en mathématiques se transforme chez Hilbert en un honnête citoyen qui ne pense que dans les limites rigoureusement définies.

L'axiomatique est un langage qui ne se décrit pas lui-même ; plus exactement, une telle description ne relève pas de ses fonctions. Dans cet univers hilbertien au lieu de demander « cet énoncé est-il vrai ? », on demande : « cet énoncé est-il formulé en accord avec les règles de la syntaxe axiomatique ? ». Hilbert remplace ainsi la notion de vérité (toujours trop conventionnelle, comme il le dit souvent) par « le choix des axiomes » (*die Wahl der Axiome*) permettant d'élaborer le texte mathématique en évitant les contradictions. Il rappelle que « notre pensée est finie et [que] le processus du penser a de même un caractère fini » [Hilbert 1935: 187] ; il faut donc choisir parmi un nombre infini de choses celles qui correspondent à nos besoins, puis bâtir à partir d'elles un système non-contradictoire. Bref, la leçon hilbertienne comporte deux aspects :

Premièrement, toute la mathématique n'est pas donnée comme sensibilité *a priori* ; sinon, elle est à construire (c'est la différence entre Kant et Hilbert) ;

[1] Hua II, voir les leçons II et IV. Le thème de ces leçons, qui traitent des problèmes de la connaissance immédiate, sera repris au § 19 des *Ideen I.* Le texte de l'*Idée* résulte tout particulièrement des cours sur la théorie générale du savoir que Husserl donne à Göttingen en 1902/03 (semestre d'hiver). Hilbert assiste à ce séminaire. D'après un extrait de son journal, Husserl est convaincu que la critique du savoir qu'il propose est plus radicale encore que celle de tous ses prédécesseurs ; cf. [Schuhmann 1977: 74]

[2] « Begründung der elementaren Strahlungstheorie » (1912), [Hilbert 1935: 218]. Dans un manuscrit de 1913 Hilbert utilise l'expression « du point de vue phénoménologique » (*vom phänomenologischen Standpunkt*) quand il parle déjà de la théorie de l'électron ; cf. *Elektronentheorie*. *Notes de cours*, Bibliothek des Mathematischen Insitituts, Universität Göttingen (semestre d'été), p. 14.

Deuxièmement, l'intuition ne nous préserve point des erreurs lorsque nous produisons des énoncés mathématiques car la distinction entre intuition et raisonnement (énoncé) formel est radicale, notamment en ce que l'intuition n'exprime que le sens interne (l'état du sujet) alors que l'axiome expose une abstraction qui ne correspond pas toujours à notre intuition[1] ; la distance entre l'intuition et l'axiomatisation est semblable à celle qui sépare l'idée d'une très grande échelle montant jusqu'à la lune et sa construction.

On imagine souvent dans l'enfance que pour construire une échelle vers la lune, il suffirait de juxtaposer deux très longs bâtons et de les joindre par un grand nombre de petits barreaux. Cette illusion est suscitée par l'intuition du monde quotidien, sensible. Cantor rejetait l'idée que les nombres sont dépendants de l'expérience sensible de l'être humain en les localisant dans l'espace intelligible, plus exactement dans une immanence d'essence qu'il appelle « l'intrasubjectivité » [Cantor 1932: cf. §8] dotée d'une existence véritable. Il est remarquable que dans *Ideen I* Husserl reprenne ce sujet en y introduisant des concepts tels que « l'intuition d'essence » (*Wesensschauung*) et « l'intuition de l'essence des choses » (*Wesenserschauung*). Chez les deux auteurs, il s'agit d'intuitionner, c'est-à-dire de *voir* une immanence objective « hors-donnée » comme essence de la pensée-même (mathématique et phénoménologique). A l'instar de Cantor, Hilbert et d'autres chercheurs des fondements, Husserl cherche éperdument son point d'appui, des concepts naturels capables de maintenir sa doctrine.

Les nombres sont des essences existantes et, pour le mathématicien cantorien, leur nature est double : « le fondement de mes réflexions étant entièrement réaliste, mais non pas moins idéaliste, il ne fait pour moi aucun doute que ces deux types de réalité se trouvent toujours conjoints <…> » [Cantor 1932: cf. §8]. On se convaincra plus loin que Husserl – toujours dans *Ideen I* – reprend cette attitude en développant la pratique de la conscience pure. Voir les phénomènes consiste à réaliser les objets idéaux ; *réaliser* veut dire les saisir dans leur essence. Le sujet transcendantal des *Ideen I* est pour beaucoup la radicalisation du sujet en mathématiques articulé dans les travaux de Dedekind et de Cantor.

Proposons une hypothèse : le sujet découvert en mathématiques durant leur crise est *le* sujet qui se développera dans *Ideen I* en tant que sujet transcendantal. Autrement dit, la généalogie du sujet transcendantal remonte aux mathématiques, ce qui n'atténue guère la complexité phénoménologique susceptible d'être restituée dans le système husserlien. A la question : « Comment le sujet transcendantal est-il lié avec les nombres dedekindiens et les ensembles cantoriens ? », on peut apporter la réponse suivante : par les abstractions créatrices mises à nu justement au cours de la recherche des fondements.

[1] La science moderne montre bien que l'intuition peut nous tromper dès que l'on dépasse les frontières de l'expérience quotidienne. La théorie de la relativité restreinte prescrit le ralentissement du temps lorsqu'il s'approche de la vitesse de la lumière ; ou bien, selon la physique quantique, si nous essayons de localiser un objet minuscule comme un électron, nous n'arriverons pas à le localiser avec certitude. Ces effets sont incompréhensibles pour l'intuition enracinée dans le monde habituel. Un autre exemple : l'axiome du choix (cf. supra), qui nous offre de choisir un élément de l'ensemble *M* jusqu'à épuiser tous les éléments de *M*. L'intuition ordinaire s'y oppose en y voyant un effort sans fin et dépourvu de sens. Mais l'axiome du choix est une idéalisation du résultat cantorien selon lequel tout ensemble peut être ordonné en principe. Il est vrai que l'intuition d'un mathématicien, d'un philosophe ou bien d'un phénoménologue est beaucoup plus avancée, dirigée vers l'abstrait, et lui suggère un autre problème : quels sont les moyens de conceptualiser de telles intuitions ? Pour le détail des débats, voir [Torretti 1983 : chap. VI]

Sur le plan historique, le sujet transcendantal est l'abstraction créatrice à laquelle Husserl attribue le plus haut degré épistémologique. Une fois mis à jour, son rôle consiste à reconstituer la science de la conscience, à la recréer. Il ne faut jamais oublier que la théorie de la conscience pure est avancée pour la même raison que Dedekind, Cantor et d'autres ont inventé la nouvelle mathématique : mettre un terme à la crise de la philosophie dont la valeur a été mise en doute par une armée de psychologues[1]. Encore une question : est-ce que Husserl anéantit ainsi le sujet mondain, petit individu travaillant dans la mesure de ses possibilités avec les abstractions mathématiques ou philosophiques ? Pas tout à fait. Au cours du processus phénoménologique développé dans *Ideen I*, le sujet mondain se voit transformé en sujet transcendantal, comme la chrysalide en papillon.

## CHAPITRE 7. L' « ARITHMETICA REALIS » DE FREGE

En 1675, Pietro Mengoli (1626-1686), professeur de mathématiques à l'Université de Bologne dont Leibniz admirait les travaux, publie l'*Arithmetica realis*[2]. Ce texte dédié au cardinal Leopoldo de' Medici est consacré aux fondements symboliques de l'arithmétique [cf. Vacca 1915]. Bien qu'au fait des systèmes de notation précédents (celles de Descartes et de Viète), Mengoli préfère user du sien propre : pour exprimer le rapport entre deux grandeurs il utilise « ; », alors que le symbole cartésien est « à » [Mengoli 1675 : 4]. Très brièvement, le but de cet auteur consiste en l'élaboration d'un système arithmétique minimal, c'est-à-dire des règles de base visant à éviter toute ambiguïté dans les résultats. Plus de deux cents ans après la parution de l'œuvre de Mengoli, Frege revient à ce problème : comment atteindre l'objectivité dans la désignation des procédés arithmétiques, comment y voir le raisonnement logique avec le plus de clarté ?

Il est fort probable que Frege se soit mis à rédiger *Les fondements de l'arithmétique*[3] (1884) sur le conseil de Stumpf [Frege 1976: 256-257], quand il ne voyait dans ce livre qu'une propédeutique destinée à ceux qui s'intéressent aux lois générales de la pensée. Pour Frege comme pour Bolzano, les mathématiques sont la *scientia prima* ; jouant le rôle du législateur, cette dernière fournit aux autres sciences, philosophie comprise, la grammaire de la pensée. Plus précisément encore, c'est la logique innée de la raison mathématique qui « grammaticalise » les actes de pensée et les dote ainsi d'une objectivité. Dans ses travaux, Frege ne cherche rien d'autre que cette objectivité exprimée par la logique générale, celle qui gouverne les fondements de notre pensée. Le nombre existe dans une objectivité que l'on décrit par la logique : il ne résulte donc pas des actes psychiques humains – selon une idée chère aux

[1] Pour plus d'informations : [Freuler 1997]

[2] Titre complet : *Arithmetica realis, serenissimo et reverendissimo principi Leopoldo ab Etrvria cardinali Medices dicata a Petro Mengolo*.

[3] Titre originel : [Frege 1884]

psychologues –, mais se laisse présenter à la conscience qui cherche à comprendre et à fixer ses lois.

L'arithmétique est à la fois l'objet et le langage, liés l'un à l'autre par cette même logique fondamentaliste[1] ; en tant qu'objet, l'arithmétique consiste en briques élémentaires du calcul dont la version traditionnelle ne connaît que quatre opérations : l'addition, la division, la multiplication et la soustraction (cette science est ensuite élargie par l'inclusion de l'étude des nombres réels sous la forme du développement décimal illimité, ou même de concepts plus avancés, comme l'exponentiation ou la racine carrée). L'*arithmétique-objet* aux yeux de Frege présente une objectivité pure n'appartenant ni au monde empirique ni à l'*ego* pensant. C'est un objet qui n'apparaît qu'au moment où la conscience conceptualise ce qu'elle voit. Si, par exemple, nous voyons quatre arbres en les désignant par le chiffre « 4 », cela ne facilite pas notre appréhension de ces arbres mais établit un recouvrement, comme le dirait Husserl, entre les *quatre arbres vus par nous* et le *chiffre « 4 » en soi,* comme une synthèse ; « attribuer un nombre [à une chose], c'est énoncer un concept » [Fege 1884: 59][2]. Réciproquement, affirmer qu'il y a quatre arbres devant nos yeux ne signifie que le contenu de cet énoncé. Rappelons que ce thème renvoie au § 8 du *Begriffsschrift*, où Frege réfléchit sur le contenu des signes ; se référant à leur objet (ou *Sinn*), ces derniers véhiculent le même contenu nuancé selon le contexte. Du point de vue logique, ces signes sont tautologiques et leur utilisation ne peut être justifiée que contextuellement. Quant au contenu conceptuel, *Sinn*, il représente une objectivité logique tout en dépassant les limites linguistiques.

Il est inutile de chercher l'objectivité arithmétique dans une chose réelle, de même qu'il est inutile de décrire le nombre comme une réalité ; en revanche il faut s'abstraire de toute réalité à laquelle on peut appliquer les concepts arithmétiques afin de voir ces derniers dans leur essence propre. Si le *Begriffsschrift* est une tentative de restituer la logique innée et d'atteindre ainsi la non-ambiguïté des jugements de vérité, les *Fondements* tâche de clarifier la nature des propositions mathématiques, de *l'arithmétique-langage*. Pour Kant, dont l'exemple « 7 + 5 = 12 » est couramment cité, les énoncés arithmétiques sont synthétiques *a priori* ; pour Frege ils sont *a priori* analytiques. « 7 + 5 = 12 » ou toute autre égalité est un énoncé arithmétique dont nous pouvons estimer la véracité même si ses éléments sont très grands. Or, Kant se trompe en disant qu'il faut avoir recours à l'intuition afin d'obtenir le résultat d'une telle opération. L'*arithmétique-langage* obéit aux lois de la logique comme l'allemand à sa grammaire. Il s'agit donc d'une matrice analytique établie de façon apriorique et excluant la nécessité de l'intuition[3]. Comme bien d'autres critiques, Frege abandonne l'anthropologisme kantien ; les vérités mathématiques sont éternelles et ne dépendent guère de notre pensée.

Pourtant, si le nombre est une objectivité incontestable, alors il faut lui donner une définition. La difficulté principale consiste en ce que le nombre fonctionne souvent de la même façon que le prédicat : les énoncés « les tables sont quatre » et « les tables sont noires » ou « 5 + 2 = 7 » sont identiques du point de vue de leur structure logique. En tant que prédicat, le nombre semble exprimer non pas l'objectivité mais un état des choses, un contexte ou – pire encore – le regard subjectif. L'objection frégéenne est la suivante : l'usage prédicatif du nombre a peu en commun avec sa nature, voire avec son contenu conceptuel. Au § 55 nous lisons sa définition du nombre : « le nombre (n + 1) appartient à un concept *F* s'il existe un objet *a* qui tombe sous *F* tandis que le nombre *n* appartienne au concept « tomber sous *F* mais non pas *a* (*unter F fallend, aber nicht a*) » [Frege 1884 : 67]. Autrement dit, à la différence de « n + 1 », représentant une série des nombres naturels, le *n* est un méta-nombre qui signifie une nouvelle catégorie : les *concepts*[1]. Ceux-ci sont distingués des nombres d'ordre « n +1 » par être exclu de la sphère de l'expérience ou objectale. La distinction entre « n + 1 » et « n » consiste en ce qu'ils appartiennent aux catégories différentes de nombre. Au § 68 l'on peut trouver une autre définition du nombre construite en effet à partir de la première : « le nombre qui appartient au concept *F* est une extension du concept 'équipotent (*gleichzahlig*) au concept F' ». Cela veut dire que le « nombre appartenant au concept » et l' « extension du concept » sont identiques, c'est-à-dire que le nombre est un concept[2]. C'est ainsi que peut être dépassé le nombre comme prédicat en identifiant la structure logique du nombre avec celle du concept. Le nombre « n + 1 » est en partie le prédicat qui correspond à notre intuition, le défaut aux yeux de Frege, mais un défaut corrigible.

[1] Pour un nouveau regard historique de l'entreprise frégéenne cf. [Künne 2010]

[2] [… die Zahlangabe eine Aussage von einem Bregriffe enthalte]

[3] Les attaques sur l'intuitionnisme kantien commencent au § 12. En revanche, au § 89 Frege reconnaît la justesse de Kant qui appelle les vérités géométriques « synthétiques et a priori » mais rejette son idée que nos intuitions ont le caractère sensible.

Ajoutons que la tentative de Frege de classer le *conceptuel* dans une catégorie spécifique renvoie au *Wissenschaftslehre* (§ 133) de Bolzano[3], où le penseur bohémien fait la distinction entre le conceptuel et ce qu'il appelle « les jugements empirique, intuitif et perceptif ». Les jugements conceptuels, dit-il en modifiant « la méthode transcendantale » de Kant, sont tels qu'ils concernent uniquement des concepts purs. Ils sont à construire et ne dépendent jamais de l'esprit humain.

Le prédicat est un nombre appliqué, de même que la « table noire » où « noire » décrit une table. Conceptuellement, le nombre est comme la couleur qui existe dans le monde physique en-dehors de la perception du sujet. En appliquant un nombre à une chose, nous créons des classes ou des ensembles de choses et les unifions dans un concept ; dire qu' « il y a quatre chevaux dans l'écurie » signifie qu'il existe un ensemble de quatre objets nommé « chevaux dans l'écurie ».

La nature du nombre est mieux exprimée dans l'idée de la cardinalité : dans ce cas, le nombre n'est pas un prédicat mais un concept, c'est-à-dire un nom abstrait n'ayant aucun contenu concret (tables, chevaux, étoiles…), tel que « 4 », « 5 », « 12 », etc. En d'autres termes, pour Frege le nombre est un objet existant dans l'arithmétique comme dans un système spécifique de la pensée où il représente seulement l'idée de la comptabilité. La nature conceptuelle – la cardinalité – du nombre le sépare du prédicat, portant la marque indélébile de la psychologie car le prédicat vient du sujet et du contexte dans lequel il se trouve *hic et nunc*. Donner un prédicat à un objet revient donc à exprimer un avis subjectif ; le prédicat trahit la présence du sujet, sa domination et son rôle dans la fabrique du sens. Le nombre comme objet conceptuel, soutenu par Frege, fait écho aux conceptions classiques, telles la théorie des nombres nodaux [cf. Wright 1983][1] ; étant une abstraction du second degré, ils représentent un ensemble quiconque.

---

[1] C'est une approche opposée à celle de Mill, que Frege critique dans son texte (cf. §§ 9 et 23), et pour qui les nombres sont induits des faits empiriques que nous observons. L'exemple d'un très grand nombre ne correspondant à aucun fait empirique vise à ridiculiser l'idée de Mill. Dans les années soixante du XX$^{ème}$ siècle, la conception frégéenne du nombre est mise à l'épreuve par Paul Benacerraf qui souligne la non-existence des nombres comme entités objectives. Cf. [Benaceraf 1964]

[2] Matthias Schirn rejette cette identification en parlant des objets *deus ex machina* ; cf. [Schirn 1990]

[3] Pour une discussion plus détaillée cf. [Künne 2010: 759-770]

En débarrassant de tout psychologisme – démarche qui a les faveurs de Russell[2] –, Frege construit un autre schème des prédicats proche de Cantor : le prédicat devient une notion générale (cardinale) à partir de laquelle sont engendrés l'ensemble des objets. Ces notions sont égales si elles contiennent le même nombre d'éléments ou le même contenu. De ce point de vue, le nombre devient le prédicat de cette notion ; « 4 » signifie « quatre éléments d'un ensemble », « 1 » signifie que l'ensemble a un seul objet, par exemple dans l'énoncé « le satellite naturel de la terre », etc. Le nombre naturel n'est défini ni par sa référence particulière ni par une opération psychique du sujet mais comme concept logique doté de certaines propriétés. Le § 10 des *Fondements* détaille : chaque nombre possède un caractère unique ; ainsi, « 3 » est un nombre premier[3], mais « 4 » n'en est pas un.

Cependant, une difficulté se fait jour tout de suite car Frege n'arrive pas à distinguer clairement les moyens logiques. On peut supposer qu'il imagine une structure idéale comportant toutes ces propriétés du nombre, celle-ci dirigeant en outre la pensée humaine telle une loi. Sa démarche est erronée : Frege construit pour l'arithmétique un système d'ordre supérieur alors qu'il considère cette même arithmétique comme un système total de pensée. Au § 14 il affirme : « les vérités arithmétiques gouvernent (*beherrschen*) tout ce qui est dénombrable ; c'est le domaine le plus large auquel appartiennent non seulement le réel (*das Wirkliche*) et l'intuitif (*das Anschauliche*), mais aussi tout le pensable. Partant, les lois du nombre ne devraient-elles pas être liées de la façon la plus intime (*innigsten*) avec les lois de la pensée ? » [Frege 1884 : 21].

---

[1] Guillermo Rosado Haddock souligne que dans ses premiers travaux Frege distingue l'objet et le concept, ce peut en effet rendre perplexe. L'objet est un nombre particulier, par exemple « 20 », auquel nous nous référons dans un énoncé précis ; le concept se réfère à une idée de l'objet ou à une classe d'objets. Pour Frege, le nombre et le concept ont une structure logique différente. [Rosado Haddock 2006: 13]

[2] Sa lettre du 16 juin 1902 ; cf. [Frege 1980: 59]. Russell montre que l'une des règles introduite dans les *Fondements*, la compréhension non restreinte, rend la théorie de Frege contradictoire. Le paradoxe révèle une antinomie, une contradiction interne à la théorie. Frege souhaitait dans cet ouvrage fonder les mathématiques sur des bases purement logiques, tâche poursuivie par Russell lui-même (cf. les *Principia Mathematica*, 1910-1913 en collaboration avec Alfred Whitehead). Ce paradoxe est relevé pour la première fois paru dans *The Principles of Mathematics* de Russell, publié en 1903.

Cette lettre bouleverse la vie de Frege ; il suspend pour longtemps son cours sur l'arithmétique et se met à chercher la solution du paradoxe de Russell, mais en vain. Pendant cette période de désespoir, Frege pense trouver les fondements des mathématiques dans la géométrie. Juste avant la publication du second volume de *Fondements*, en 1903, Frege ajoute un appendice où il avoue : « Pour un auteur scientifique, il est peu d'infortunes plus amères que de voir l'une des fondations de son travail s'effondrer alors que celui-ci s'achève. C'est dans cette situation inconfortable que m'a mis une lettre de M. Bertrand Russell, alors que le présent volume allait paraître » ; [Wehmeier 2004 ; Kratzsch 1979 ; de même Frege 1969 ; 1973]

[3] Un nombre premier est un entier naturel qui admet exactement deux diviseurs entiers distincts et positifs, 1 et lui-même. Le plus grand nombre premier connu à ce jour est $M_{43\,112\,609} = 2^{43\,112\,609}-1$ ; il s'agit du 45[e] nombre premier de Mersenne.

Ou bien, au § 38, Frege décrit l'arithmétique comme un système symbolique dont les symboles sont irremplaçables : « elle [l'arithmétique] verrait sa fin si, à la place du nombre un, qui est toujours le même, nous introduisions pour des choses diverses des symboles similaires [mais non pas identiques] » [Frege 1884 : 49]. Le *un* désigne l'existence d'une chose concrète et donc sa quantité, d'où vient le caractère analytique de l'un : si la chose existe, elle est « un », concept généralisant l'existence d'une chose concrète. Pourtant, pour Frege, il s'agit seulement d'une généralité (*Allgemeinheit*) limitée ne couvrant qu'un objet, physique ou mental. Il explique : si l'on désigne ainsi le point dans l'espace, l'arbre dans le jardin, ou un intervalle de temps, le nombre comme tel ne nous aidera guère car l'intervalle temporel ou le point dans l'espace a sa propriété spécifique (cf. § 41), sans rapport avec celle de nombre.

L'idée consiste donc à remplacer le caractère purement numérique du nombre par un nombre symbolique ou, tout simplement, par un symbole. Ce dernier doit être le résultat d'une synthèse de différentes propriétés, non seulement des nombres naturels (chez Schröder : « le nombre naturel est une somme de uns »). Quant à Frege, quel type de synthèse construit-il ? Si le jugement synthétique kantien est subjectif ou trop anthropologique (rappelons que Hilbert critique Kant pour la même raison), car le sujet kantien intuitionne puis synthétise « 12 », le conférencier viennois cherche une synthèse faite non par le sujet mais à partir du système arithmétique-même. Ce dernier va jouer le rôle du sujet. L'intuition kantienne devient chez Frege (cf. § 48) une puissance conceptuelle : en rassemblant des choses diverses dans une unité, elle couvre largement toute apperception synthétique *a priori*. Le concept du nombre dépasse le nombre. Sur le plan pragmatique, le concept doit remplacer le nombre ; ce remplacement est un acte fondateur.

Dès le début de son traité, Frege cherche à construire l'échafaudage logique sur lequel il fonde l'arithmétique. Même après y être parvenu, l'échafaudage demeure. Voilà l'impasse méthodologique : d'un coté, la conceptualisation de l'arithmétique sert de fondement ; d'un autre coté, les concepts ne doivent pas dépasser la pensée arithmétique *per se*. Le concept du nombre résulte de sa nature, mais pour voir cela il est nécessaire d'avoir déjà une certaine intuition du concept. Bref, dans cette situation la logique s'apparente à un langage sans locuteur car Frege laisse de coté la question du sujet transcendantal qui pourrait intuitionner et bâtir le système conceptuel de l'arithmétique.

L'impasse commence au § 17 du traité, où les vérités de l'arithmétique équivalent à des théorèmes de géométrie. Cela veut dire que la logique et les axiomes géométriques sont les schèmes généraux à partir desquels on déduit graduellement des vérités particulières. Il s'agit ici d'une variation par rapport à Kant, qui traite la nature du savoir mathématique comme appartenant à la catégorie de la raison pure ; par là, d'un geste « leibnizien », Frege modifie le projet kantien en remplaçant le sujet transcendantal par une notation transcendantale[1].

[1] Cf. aussi § 15. Jean-Pierre Belna note que, tout comme Leibniz, « Frege constate l'inadéquation de la langue quotidienne pour certains buts scientifiques » ; cf. [Belna 1996 : 201]. Cependant, le projet leibnizien était trop vaste pour être réalisé et Frege se tourne vers Boole. Frege commence à s'intéresser à l'algèbre de Boole vers 1874 mais ne trouve pas dans le langage booléen un véritable instrument pour exprimer le contenu des énoncés arithmétiques

Il est intéressant de remarquer que dans la constitution de son concept du nombre (ou du nombre comme concept), Frege effectue l'opération opposée de ce que Kant appelle la « régression empirique » (*empirischen Regreß*) [1] signifiant qu'une chose intuitivement saisie peut se donner élément par élément sans aucune limite préétablie ; nul élément n'a de condition nécessaire. Autrement dit, l'objet perçu entièrement peut ne pas être reconstitué dans le monde empirique. Il ne s'agit que d'une possibilité. Le nombre chez Frege résulte d'une « de-régression » empirique, c'est-à-dire qu'avant tout il doit être saisi intuitivement (à comparer : le jugement synthétique chez Kant[2]), comme un tout dont les éléments sont reconstitués dans le concept. Il y d'ailleurs un aspect de la *Critique* soutenu par Frege : les concepts mathématiques, précise Kant, doivent être libres de tout contenu psychologique car ils résident hors de toute expérience mondaine[3] ; « aucune objectivité ne peut, bien sûr, être basée sur une impression de sensation (*Sinneseindrucke*) qui en tant que affection de notre conscience est entièrement subjective <…> cette objectivité ne peut être basée que sur la raison » [Frege 1884 : 38][4]. Frege ne fait pas partie du camp kantien ou, si l'on veut, au sens très large du terme. La trajectoire épistémologique de Kant va de « que savoir » à « comment savoir », plus exactement : comment constituer le monde du savoir transcendantal ? La question de Frege, toujours actuelle à son époque, est la suivante : comment construire la logique des énoncés arithmétiques ?[5]

C'est pourquoi la démarche cantorienne de Frege n'est pas fortuite. Au § 53 des *Fondements* il introduit « l'unicité » (*Einzigkeit*) qui représente un concept d'ordre supérieur. La tâche de l'unicité est de recueillir tous les concepts signifiant le même objet sous un seul concept singulier. Nous passons ainsi à un nouveau concept inconnu au système précédent ; comme le transfini de Cantor, l'unicité est engendrée par le système le plus faible, c'est-à-dire dont le niveau d'abstraction est inférieur au nouveau produit. Frege reconnaît lui-même ce qu'il doit au théoricien des ensembles : « Je crois que je suis d'accord avec Cantor, mais ma terminologie se distingue de la sienne. Au lieu de mon « nombre » il utilise « la puissance », tandis que son concept du nombre signifie l'ordonnancement <…> notre conception du nombre couvre tous les nombres

---

(nous nous rappelons que l'aspect de contenu Boole laisse aux interprétations) ; cf. [Sluga 1987]

[1] Cf. [Kant 1956 : Ab. IX]

[2] *Ibid*., « Einleitung », V. Néanmoins, Frege rejette l'affirmation de Kant selon laquelle tous les jugements mathématiques sont synthétiques. Cette question est spécialement étudiée dans [MacFarlane 2002 : 6 et passim]. Hao Wang remarque : « Frege thought that his reduction refuted Kant's contention that arithmetic truths are synthetic. The reduction, however, cuts both ways <…> if one believes firmly in the irreducibility of arithmetic to logic, he will conclude from Frege's or Dedekind's successful reduction that what they take to be logic contains a good deal that lies outside the domain of logic » ; cf. [Wang 1957: 80; de même Nidditch 1963: 107-109].

[3] Dummett attire l'attention sur le fait que Frege « argued that <…> psychological accounts are valueless, and must be replaced by definitions that specify the contribution made by the expression defined to the condition for the truth of a statement in which it occurs » ; cf. [Dummett 1991 : 18]

[4] Chez Kant : « Daher ist reine Vernunft diejenige, welche die Prinzipien, etwas schlechthin a priori zu erkennen, enthält » ; [Kant 1956 : 55]

[5] Des nombreuses études consacrées au traité de Frege citons : [Boolos 1987]. Hartry Field explique que la logique au sens strict du terme ne saurait fonder les objets mathématiques, cf. [Field 1984] ; Harold Hodes voit dans la réduction logiciste de Frege une incompatibilité originelle, cf. [Hodes 1984: 123]

infinis et n'exige donc aucune extension » [Frege 1884 : 97-98]. En outre, Frege affirme qu'il manque au nombre chez Cantor une précision (conceptuelle, logique) alors qu'il conserve une mystérieuse « intuition interne », allusion claire à Kant. Ce dernier, note Frege vers la fin de son traité, sous-estime la valeur des jugements analytiques, rejoignant sur ce point la critique bolzanienne. L'idée de l'analyticité consiste à établir un jugement bien défini, c'est-à-dire qu'elle est indépendante du sujet. Ce que Kant attribue à l'intuition, en étendant le champ anthropologique de la raison, est inclus dans la logique du jugement analytique, c'est-à-dire dans la vérité décrite par ce jugement. Enfin, Kant se trompe lorsqu'il postule les formes transcendantales de la sensibilité car rien ne nous empêche, insiste Frege, de conceptualiser un nombre vertigineux tel que $100^{1000}$ sans en avoir aucune sensation.

Pour récapituler : les efforts frégéens sont dirigés vers la recherche de la vérité logique dont le meilleur exemple est le nombre conceptualisé par la pensée arithmétique. Le nombre est un objet à construire et cette construction emprunte le chemin allant des notations naturelles vers les formes conceptuelles. Ce qui existe, ce que l'on peut dénombrer appartient aux concepts *in concreto* et n'englobe jamais le concept en entier. Le nombre est une complétude où sont recueillis les cas concrets et c'est à partir d'elle qu'est déduit tout énoncé. Cette tendance à modifier l'arithmétique par des moyens logiques ne fait que s'accentuer à l'époque où Frege publie ses principaux ouvrages ; certains voient même dans l'axiomatique hilbertienne une tentative de donner à la pensée mathématique une base logique, sans compter les autres recherches évoquées ci-dessus. A la fin du XIX^ème^ siècle, Husserl avec ses études sur les mathématiques entre en scène, occupée alors par deux stratégies du savoir – la logique et la psychologie – qui chacune détermine la position du chercheur par rapport à ce savoir.

## Bibliographie

Ahrens, W. 1907. « Skizzen aus dem Leben Weierstraß' », *Mathematisch-naturwissenschaftliche Blatter* 4.

Alanus ab Insulis, *Regulae de sacra theologiae*. In [Migne 1844-1854]

Anderson, L.R. 2005. «Neo-Kantianism and the Roots of Anti-Psychologism», *British Journal for the History of Philosophy* vol. 13, 2.

Andreka, H., Monk, D., Nemeti, I. 1991 (eds.) *Algebraic Logic. Colloquia Mathematica Societatis Janos Bolyai* vol. 54, éd. par (Amsterdam, Elsevier Science/North-Holland.

Anellis, I.H., Houser, N. 1991. «The Nineteenth Century Roots of Universal Algebra and Algebraic Logic: A Critical-bibliographical Guide for the Contemporary Logician». In [Andreka, H., Monk, Nemeti 1991]

Aspray, W., Kitcher, P. 1988 (eds.) *History and Philosophy of Modern mathematics*, University of Minnesota Press, Minneapolis.

Avé-Lallemant, E. 1975. *Die Nachlässe der Münchener Phänomenologen in der Bayerischen Staatsbibliothek*. Wiesbaden, Otto Harrassowitz.

Bakker, R. 1969. *De Geschiedenis van het Fenomenologisch Denken*. Nederlands, Het Spectrum.

Bar-Hillel, Y. 1970. *Aspects of Languag*e, Jerusalem, Magnes Press.

Belna, J.-P. 1996. *La Notion de Nombre chez Dedekind, Cantor, Frege*, Paris, Vrin.

Belna, J.-P. 2000. *Cantor*, Paris, Les belles lettres.

Benacerraf, P., Putnam, H. 1964 (eds.) «What Numbers Could Not Be», *Philosophy of Mathematics: Selected Readings*, Englewood Cliffs, N.J., Prentice-Hall.

Benoist, J. 2002/3. « La réécriture par Bolzano de l'esthétique transcendantale ». In *Revue de métaphysique et de morale* 35.

Berg, J. 1999. «Kant über analytische und synthetische Urteile mit Berücksichtigung der Lehren Bolzanos». In *Bernard Bolzanos geistiges Erbe für das 21. Jahrhundert,* hrsg. E. Morscher, Akademia Verlag Sankt Augustin.

Bernays, P. 1922. «Über Hilberts Gedanken zur Grundlegung der Arithmetik», *Jahresbericht der Deutschen Mathematiker-Vereinigung*, Bd. 31.

Bogan, O. 1890. *Versuch eines Beweises gegen die Lösbarkeit philosophischer Probleme*, Halle, Kaemmerer.

Bolliger, A. 1882. *Anti-Kant oder Elemente der Logik, der Physik und der Ethik*, Basel, Schneider.

Bolzano, B. 1985. *Wissenschaftslehre*, Bd. 11, Teil I., hrsg. J. Berg, Stuttgart–Bad Cannstatt, Friedrich Frommann Verlag.

Boniface, J. 2010. «Position philosophique et pratique mathématique : l'exemple de L. Kronecker», *Images des Mathématiques*, CNRS.

Boolos, G. 1987. «The Consistency of Frege's *Foundations of Arithmetic*». In [Thomson 1987]

Borel, E. 1914. «L'infini mathématique et la réalité», *Revue du Mois*, 18.

Bottazzini, U. 1981. *Il Calcolo sublime: storia dell'analisi matematica da Euler a Weierstrass*, Torino, Boringhieri.

Bourbaki, N. 2007. *Eléments d'histoire des mathématiques*, Berlin, Springer.

Braun, J. 2006. *Zarys filozofii Hoene Wrońskiego*, Warszawa.

Brouwer, L.E.J. 1913. «Intuitionism and Formalism», *Bulletin of the American Mathematical Society*, vol. 20, 2.

Brunschvicg, L. 1912. *Les étapes de la philosophie mathématique*, Paris, F. Alcan.

Büchner, L. 1884. «"Anti-Kant." *Aus Natur und Wissenschaft*». In *Studien, Kritiken und Abhandlungen*, vol. II, Leipzig, Thomas.

Bushaw, D. 1983. «Wronski's "Canons of logarithms" ». In *Mathematics Magazine*, 5 (2).

Cantor, G. 1874. «Über eine Eigenschaft des Inbegriffes aller reellen algebraischen Zahlen», *Journal de Crelle*, Bd. 77.

Cantor, G. 1932. *Gesammelte Abhandlungen mathematischen und philosophischen Inhalts*, hrsg. E. Zermelo, Berlin, Springer.

Carnap, R. 2000 (1937) *Logical Syntax of Language*, London, Routledge.

Cattaruzza, S. & Sinico, M. (éds.) 2004. *Husserl in laboratorio*, EUT.

Cauchy, A. 1821. *Cours d'analyse de l'école royale polytechnique, 1^er^ partie : analyse algébrique*, Paris.

Cavaillès, J. 1947. *Transfini et Continu*, Paris, Hermann.

Cavaillès, J. 1981. *Méthode axiomatique et formalisme : Essai sur le problème du fondement des mathématiques*, Paris, Hermann.

Cho, K. (éd.) 1984. *Phaenomenologica* 95, Dordrecht, Martinus Nijhoff.

Conte, A., Giacardi, L. 1991 (eds.) *Angelo Genocchi e i suoi interlocutori scientifici. Contributi dall'epistolario*, Torino, Deputazione Subalpina di Storia Patria.

Contro, S. 1976. «Von Pasch zu Hilbert», *Archive for History of Exact Sciences*, vol. 15, 3.

Couturat, L. 1901. *Logique de Leibniz,* Paris, Félix Alcan.

Dalen, D. van. 1999. «Luitzen Egbertus Jan Brouwer». In [James 1999]

Dauben, J.W. 1977. «Georg Cantor and Pope Leo XIII: Mathematics, Theology, and the Infinite», *Journal of the History of Ideas* 38, vol.1.

Dauben, J.W. 1990. *Georg Cantor : His Mathematics and Philosophy of the Infinite, Princeton University Press*.

DeBoer, Th. 1978. *The Development of Husserl's Thought*, The Hague, Nijhoff.

Décaillot, A.-M. 2008. *Cantor et la France. Correspondance du mathématicien allemand avec les Français à la fin du XIXe siècle*, Paris, Kimé.

Dedekind, R. 1932. *Gesammelte mathematische Werke*, hrsg. R. Fricke, E. Noether & Ö. Ore, Bd. III, Braunschweig, Vieweg & Sohn.

Deleuze, G. 1968. *Différence et répétition*, Paris, Epiméthée.

Desanti, J.T. 1968. *Les Idéalités mathématiques*, Paris, Seuil.

Descartes, R. 1977. *Règles utiles et claires pour la direction de l'esprit en la recherche de la vérité*, trad. selon le lexique cartésien et annotations conceptuelles par J.-L. Marion (avec des notes mathématiques de P. Costabel), La Haye.

Diels, H. 1899. *Über Leibniz und das Problem der Universalsprache*, Berlin, Sitzung d. Akademie.

Dieudonné, J. 1979. *The Tragedy of Grassmann*, Séminaire de Philosophie et Mathématiques ENS, éd. par l'I.R.E.M., Paris-Nord, 65.

Dingler, H. 1931. *Die Philosophie der Logik und Arithmetik*, München, Eidos.

Dubislav, W. 1929. «Über Bolzano als Kritiker Kants». In *Philosophisches Jahrbuch* Bd. 42, Hf. 3.

Dugac, P. 1976. *Richard Dedekind* et les *fondements des mathématiques*, Paris, Vrin.

Dühring, E. 1875. *Kursus der Philosophie als streng wissenschaftlicher Weltanschauung und Lebensgestaltung*, Leipzig, Koschny.

Dummett, M. 1991. *Frege: Philosophy of Mathematics*, London, Duckworth.

Edwards, H. 1988. «Kronecker's Place in History». In [Aspray, Kitcher 1988]

Eley, L. 1962. *Die Krise des* Apriori *in der transzendentalen Phänomenologie Edmund Husserls*, Den Haag, Martinus Nijhoff.

Feferman, S. 1987. «Infinity in Mathematis: Is Cantor Neessary?» In: *L'infinito nella Sienza*, Instituto dello Enciclopedia Italiana.

Feist, R. 2004 (ed.) *Husserl and Sciences*, University of Ottawa Press.

Fels, H. 1927. «Die Philosophie Bolzanos». In *Philosophisches Jahrbuch* Bd. 40.

Ferreirós, J. 2007. *Labyrinth of Thought. A History of Set Theory and Its Role in Modern Mathematics,* 2ème éd., Basel, Birkhäuser.

Field, H. 1984. «Is Mathematical Knowledge Just Logical Knowledge?», *Philosophical Review* 93.

Fisette, D. 2009. «Husserl à Halle (1886-1901)», *Philosophiques,* vol. 36, 2.

Fréchet, M. 1934. *L'arithmétique de l'infini*, Paris, Hermann.

Frege, G. 1884. *Die Grundlagen der Arithmetik. Eine logisch mathematische Untersuhung über den Begriff der Zahl*, Breslau, W. Kœbner.

Frege, G. 1895. «Kritische Beleuchtung einiger Punkte in E. Schröders *Vorlesungen über die Algebra der Logik*», *Archiv für systematische Philosophie* Bd. I.

Frege, G. 1969. *Nachgelassene Schriften*, Hamburg, Félix Meiner.

Frege, G. 1973. *Schriften zur Logik. Aus dem Nachlaß*. Mit einer Einleitung von Lothar Kreiser, Berlin, Akademie Verlag.

Frege, G. 1976. *Wissenschaftlicher Briefwechsel*, hrsg. H. Hermes, F. Kambartel, Ch. Thiel & A. Veraart, Hamburg, Felix Meiner.

Frege, G. 1980. *Briefwechsel mit D. Hilbert, E. Husserl, B. Russell*, hrsg. von G. Gabriel, F. Kambartel, C.Thiel, Hamburg, Felix Meiner.

Freguglia, P. 2011. « Geometric Calculus and Geometry Foundations in Peano ». In *Giuseppe Peano between Mathematics and Logic*, éd. par F. Skof, Springer.

Frei, G. 1985. *Der Briefwechsel David Hilbert-Felix Klein (1996–1919),* Göttingen, Vandenhoeck & Ruprecht.

Freuler, L. 1997. *La crise de la philosophie au XIXe siècle,* Paris, Vrin.

Friedmann, G. 1962. *Leibniz et Spinoza*, Paris, Gallimard.

Gandt, F. de. 2004. « Göttingen 1901 : Husserl et Hilbert ». In [Worms 2004]

Gauthier, Y. 2004. «Husserl and the Theory of Multiplicities *"Mannigfaltigkeitslehre"*». In [Feist 2004]

Giaquinto, M. 1983. «Hilbert's Philosophy of Mathematics», *British Journal for Philosophy of Science,* vol. 34.

Gonseth, F. 1939. « La méthode axiomatique », *Bulletin de la S.M.F.*, t. 67.

Gonseth, F. 1945. *La géométrie et le problème de l'espace. La doctrine préalable,* Neuchâtel, Griffon.

Grabiner, J.V. 1981. *The Origins of Cauchy's Rigorous Calculus*, New York, Dover Publications.

Grabmann, M. 1957. *Die Geschichte der scholastischen Methode*, Bd. II, Berlin, Akademie-Verlag.

Grassmann, H. 1878. *Die Ausdehnungslehre von 1844, oder die lineale Ausdehnungslehre*, Leipzig, Wiegand.

Grassmann, H. 1894. *Gesammelte mathematische und physikalische Werke*, Bd. I, T. 1, hrsg. F. Engel, Leipzig, Teubner.

Grattan-Guinness, I. 1970. *The Development of the Foundations of Mathematical Analysis from Euler to Riemann,* Cambridge, Massachusetts, MIT Press.

Guillaume, M. 2008. «Some of Julius König's Mathematical Dreams in His *New Foundations of Logic, Arithmetic, and Set Theory*». In *One Hundred Years of Intuitionism (1907-2007)*, II.

Hamacher-Hermes, A. 1991. «The Debate between Husserl and Voigt Concerning the Logic of Content and Existentional Logic», *Analecta Husserliana*, vol. 34.

Hankins, T. 1970. *Jean d'Alembert. Science and the Enlightenment*, New York, Gordon & Breach.

Hartimo, M.H. 2007. «Towards completeness: Husserl on theories of manifolds 1890-1901», *Synthese*, vol. 156.

Hartimo, M.H. 2012. «Husserl and the Algebra of Logic: Husserl's 1896 Lectures», *Axiomathes*, vol. 22, issue 1.

Heijenoort, J. van. 1967. *From Frege to Gödel*, Cambridge (MA), Harvard University Press.

Held, K. 1989. «Husserl und die Griechen», *Phänomenologische Forschung*, Bd. 22.

Hilbert, D. 1905. «Logische Principien des mathematischen Denkens», *Manuscrit de Hilbert/Notes de cours*, Bibliothek des Mathematischen Instituts, Universität Göttingen, semestre d'été, 1905 (annoté par E. Hellinger).

Hilbert, D. 1926. « Über das Unendliche », *Mathematische Annalen*, Bd. 95, Berlin, Springer.

Hilbert, D. *Gesammelte Abhandlungen*, Bd. III, Berlin, Springer, 1935

Hilbert, D. 1990. *Sur les problèmes futurs des mathématiques*, trad. par L. Laugel, Paris, J. Gabay.

Hintikka, J. 1995. *From Dedekind to Gödel. Essays on the Development of the Foundations of Mathematics*, Boston University, Kluwer Academic Publishers.

Hodes, H.T. 1984. « Logicism and the Ontological Commitments of Arithmetic », *Journal of Philosophy* 81(3).

Hoëné-Wroński, J. 1811. *Introduction à la philosophie des mathématiques*, Paris, Courcier.

Houser, N. (éd), 1990/91. «The Schröder-Peirce Correspondence», *Modern Logic* vol. 1 (2-3), 1990/91.

Hourya, B.-S. 1973. « Cauchy et Bolzano », *Revue d'Histoire des sciences*, t. 26 (2).

Husserliana (Hua): I, II, III/1, VI, XII, XVII, XXII, XIX/1, XXIV, XXI, XXXV.

Hua Dokumente III/5.

Husserl, E. 1975. *Introduction to the Logical Investigations*, The Hague.

Husserl, E. 1994. *Early Writings in the Philosophy of Logic and Mathematics*, trans. by D. Willard, Dordrecht, Kluwer.

Husserl, E. 1994. *Briefwechsel. Die Neukantianer*. In Hua Dokumente, vol. III/5, Dordrecht, Kluwer Academic Publishers.

Husserl, M. 1988. «Skizze eines Lebensbildes von E. Husserl». *Husserl Studies* 5.

Ingarden, R. 1974. *Wstęp do fenomenologii Husserla* [Introduction à la phénoménologie de Husserl], Warszawa.

Inghen, M. de. 2000. *Quaestiones super quattuor libros Sententiarum*, 2 vols. éd. par G. Wieland, M. Santos Noya, M.J.F.M. Hoenen et M. Schulze, Leiden.

Jacquette, D. 2006. *Philosophy of Logic*, Amsterdam, North Holland.James, I.M. 1999 (ed.) *History of Topology*, Amsterdam, Elsevier.

Jourdain, P.E.B. 1912. « The Development of the Theories of Mathematical logic and the Principles of Mathematics». In *Quart. Journal of Pure and Applied Mathematics*, vol. 43.

Kant, I. 1956. *Kritik der reinen Vernunft*, hrsg. R. Schmidt, Hamburg, Felix Meiner.

Kleiner, I. 2007. *A History of Abstract Algebra*, Boston-Berlin, Birkhäuser.

Kneale, K. 1948. «Boole and the Revival of Logic», *Mind*, vol. 57.

Korte, T. 2010. «Frege's *Begriffsschrift* as a *lingua characteristica*», *Synthèse*, vol. 174.

Kratzsch, I. 1979. Material zu Leben und Werken Freges aus dem Besitz Universitätsbibliothek Jena, dans Begriffsschrift. Jenaer Frege-Konferenz, Jena.

Kronecker, L. 1871. «Notiz», *Journal für die reine und angewandte Mathematik* 73.

Külpe, O. 1920. *Die Philosophie der Gegenwart in Deutschland*, Leipzig, Teubner.

Künne, W. 2009. *Die philosophische Logik Gottlob Freges. Ein Kommentar*, Frankfurt am Main, Klostermann.

Künne, W. 2010. *Die philosophische Logik Gottlob Freges. Ein Kommentar*, Frankfurt a/M., Klostermann.

La Chapelle, d'Alembert. 1751-1772. « Limite », In : *L'Encyclopédie ou Dictionnaire raisonné des sciences, des arts et des métiers*, sous la dir. de D. Diderot & D'Alembert, Paris.

Lambert, J.-H. 1764. *Neues Organon oder Gedanken über die Erforschung und Berechnung des Wahren*, Leipzig.

Landgrebe, L. 1963. *Der Weg der Phänomenologie*, Gütersloh.

Lavigne, J.-F. 2004. *Husserl et la naissance de la phénoménologie (1903-1913)*, Paris, PUF.

Laz, J. 1993. *Bolzano critique de Kant*, Paris, Vrin.

Leibniz, G.W. 1999. *Sämtliche Schriften und Briefe*, Berlin-Brandenburgischen Akademie der Wissenschaften, Akademie der Wissenschaften in Göttingen, VII *Philosophische Schriften*, Band 4, *1677-Juni 1690*, éd. par H. Schepers, M. Schneider, G. Biller, U. Franke & H. Kliege-Biller, Berlin, Akademie Verlag.

Link, G. 2004 (ed.) *One Hundred Years Of Russell's Paradox: Mathematics, Logic, Philosophy*, de Gruyter, 2004.

Lotze, H. 1843. *Logik*, Leipzig.

Loužil, J. 1978. *Bernard Bolzano*, Prague, Melantrich.

Luchins, A.S., Luchins, E.H. 2000. «Kurt Grelling: Steadfast Scholar In a Time of Madness», *Gestalt Theory*, vol. 22, 4.

Lusin, N. 1930. *Leçons sur les ensembles analytiques et leurs applications*, Paris, Gauthier-Villars.

MacFarlane, J. 2002. «Frege, Kant, and the Logic in Logicism» (Ms.)

Majer, U. 2004. «Husserl and Hilbert on Geometry». In [Fesit 2004]

McCarty, D.C. 1995. «The Mysteries of Richard Dedekind». In [Hintikka 1995]

Manin, Y. 2002. «Georg Cantor and His Heritage», http://arxiv.org/abs/math/0209244v1

Mengoli, P. 1675. *Arithmetica realis, serenissimo et reverendissimo principi Leopoldo ab Etrvria cardinali Medices dicata a Petro Mengolo*, Bologne, Benatij.

Migne, T.P. 1844-1854. *Patrologiae cursus completus*. Series latina, Parisiis, t. 210.

Moese, H. 1965. « La "Wissenschaftslehre" de Bolzano et le problème de l'unité de la science ». In *Actes du XIe Congrès International d'Histoire des Sciences*, Wroclaw-Varsovie-Cracovie, t. 2.

Moore, G.H. 2002. «Hilbert on the Infinite: The role of set theory in the evolution of Hilbert's thought, *Historia Mathematics* 29.

Morscher, E. 1972. « Vom Bolzano zu Meinong : zur Geschichte des logischen Realismus ». In *Jenseits von Sein und Nichtsein*, Beiträge zur Meinong-Forschung, hrsg. R. Haller, Graz.

Mugnai, M. 1987(1992). «Leibniz and Bolzano on the "Realm of truths"». In *Bolzano's* Wissenschaftslehre *1837-1987*, International Workshop, Firenze.

Natanson, M. 1973. *Edmund Husserl: Philosopher of Infinite Tasks*, Evanston, Northwestern University Press.

Nelson, L. 1973. «Über das sogennante Erkenntnisproblem», *Schriften zur Erkenntnistheorie* (hrsg). P. Bernays et al., *Gesammelte Schriften*, Bd. 2, Hamburg, Felix Meiner.

Nidditch, P. 1963. «Peano and the Recognition of Frege», *Mind*, LXXII (285).

Noutsoubidzé, C. (ნუცუბიძე, შ.) 1926. *Wahrheit und Erkenntnisstruktur*, Berlin, Walter de Gruyter.

Osborn, A.D. 1934. *The Philosophy of Edmund Husserl in Its Development from His mathematical Interests to His First Conception of Phenomenology in* the Logical Investigations, New York, International Press.

Padoa, A. 1899. «Ideografia delle frazioni irriducibili», *Revue des mathématiques*, 6.

Padoa, A. 1901. «Essai d'une théorie algébrique des nombres entiers, précédé d'une introduction logique à une théorie déductive quelconque »,

*Bibliothèque du Congrès International de Philosophie,* vol. III : *Logique et histoire des Sciences*, Paris, Colin.

Pasch, M. 1882. *Vorlesungen über neuere Geometrie,* Leipzig, Teubner.

Peano, G. 1906. *Super Theorema de Cantor-Bernstein*. In *Rendiconti del Circolo Matematico di Palermo,* t. 21.

Peano, G. 1958. *Opere scelte,* éd. U. Cassina, t. II, Roma, Cremonese.

Peano, G. 1892. « Sopra la raccolta di formule di matematica », *Revue des mathématiques,* t. 2.

Peckhaus, V. 1990. *Hilbertprogramm und Kritische Philosophie* : *Das Göttinger Modell interdisziplinärer Zusammenarbeit zwischen Mathematik und Philosophie,* Göttingen, Vandenhoeck & Ruprecht.

Peckhaus, V. 1990/91. « Ernst Schröder und die "pasigraphischen Systeme" von Peano und Peirce », *Modern Logic* vol. 1 (2-3).

Petsche, H.-J. 2006. *Graßmann,* Basel, Birkhäuser.

Phili, C. 1996. «La loi suprême de Hoëné-Wroński: la rencontre de la philosophie et des mathématiques». In *Paradigms and mathematics,* Siglo XXI Espana Ed., Madrid.

Pieri, M. 1901. «Sur la géométrie envisagée comme un système purement logique», *Bibliothèque du Congrès International de Philosophie,* vol. III : *Logique et histoire des Sciences,* Paris, Colin.

Plessner, H. 1959. «Husserl in Göttingen», *Rede zu Feier d. 100 Geburtstages Edmund Husserls,* Göttingen, Vandenhoeck & Ruprecht.

Poincaré, H. 1902. *La science et l'hypothèse,* Paris, Flammarion.

Poretsky, P. 1899. Sept lois fondamentales de la théorie des égalités logiques, Edition de l'Université de Kazan.

Poretsky, P. 1900. «Exposé élémentaire de la théorie des égalités logiques à deux termes *a* et *b*», *Revue de Métaphysique et de morale,* t. 8.

Pragacz, P. 2007. « Notes on the life and work of Józef Maria Hoëné-Wroński », trad. anglaise par J. Spaliński, *Wiadomości Matematyczne,* t. 43.

Příhonský, F. 2006. *Bolzano contre Kant. Le nouvel Anti-Kant,* Paris, Vrin.

Pulkkinen, J. 2005. *Thought and Logic: The debates between German-Speaking Philosophers and Symbolic logicians at the turn of the 20th century*. In Europaische Studien zur Ideen und Wissenschaftsgeschichte, Bd. 12, Peter Lang Publishing.

Putnam, H. 1982. «Peirce the Logician », In *Historia Mathematica,* vol. 9.

Quine, W.V.O. 1980. «Two Dogmas of Empiricism». In *From a Logical Point of View*, Harvard University Press.

Reid, C. 1970. *Hilbert*, Berlin, Springer-Verlag.

Richir, M. 2002. *L'institution de l'idéalité. Des schématismes phénoménologiques*, Beauvais, Mémoires des Annales de Phénoménologie.

Rickert, H. 1892. *Der Gegenstand der Erkenntnis: ein Beitrag der philosophischen Transcendenz,* Freiburg.

Rosado Haddock, G.E. 2006. *A Critical Introduction to the Philosophy of Gottlob Frege,* Ashgate, Chippenham, Wiltshire.

Rosa Massa Esteve, M. 2008. *L'algebrizatció des les matemàtiques. Pietro Mengoli (1625-1686)*, Barcelona, Societat Catalana d'Història de la Ciència i de la Tècnica.

Rowe, D., McCleary, E.J. 1989 (eds.) *History and Philosophy of Modern mathematics*, vol. I, Ideas and Their Reception, Boston, Academic Press.

Royce, J. 1913. «An Extension of the Algebra of Logic», *Journal of Philosophy* vol. 10.

Sallis, J. « The Identities of the Things Themselves », *Philosophy and Science in Phenomenological Perspective*. In [Cho 1984].

Schirn, M. 1990. «Frege's Objects of a Quite Special Kind», *Erkenntnis* 32.

Schnädelbach, H. 1984. *Philosophie in Deutschland, 1831-1933,* Frankfurt a/Main, Suhrkamp.

Schneider O'Connell, A. 1988. «Husserl and Frege on Schröder. The Shoe on the Other Foot », *Etudes phénoménologiques*, 8, 1988.

Schröder, G. 1966 [reprint] *Vorlesungen über die Algebra der Logik* (*exakte Logik*), three vols., New York, Chelsea.

Schuhmann, K. 1977. *Husserl-Chronik,* The Hague, Martinus Nijhoff.

Schuhmann, E., Schuhmann, K. 2001. «Husserls Manuskripte zu seinem Göttinger Doppelvortrag von 1901», *Husserl Studies* 17, vol. 2.

Schwartz, E. 1996. « Le jeune Husserl, lecteur de Schröder ». In *Actes du Congrès International Henri Poincaré 1994* (article paru dans *Philosophia Scientiae*, vol. I, Cahier 2).

Sebestik, J. 2003. «Husserl Reader of Bolzano». In *Husserl's* Logical Investigations *Reconsidered,* ed. by D. Fisette & J.J. Drummond, vol. 48, Dordrecht, Kluwer Academic Publishers.

Segal, S.L. 2003. *Mathematicians under the Nazi*, Princeton University Press.

Sepp, H.R. (éd.) 1988. *Edmund Husserl und die phänomenologische Bewegung*. Freiburg/München, Alber.

Serfati, M., Descotes, D. (éds.) 2008. *Mathématiciens français du XVII$^{ème}$ siècle : Descartes, Fermat, Pascal*, PUBP.

Sluga, H. 1987. «Frege against the Booleans», *Notre Dame Journal of Formal Logic* vol. 28, 1.

Snapper, E. 1979. «The Three Crises in Mathematics: Logicism, Intuitionism and Formalism». In *Mathematics Magazine* (524).

Spiegelberg, H. 1960. *The Phenomenological Movement: A Historical Introduction*, 2 vols.The Hague, Nijhoff.

Stein, H. 1989. «*Logos*, Logic, and *Logistiké* : Some Philosophical Remarks on Nineteenth-Century Transformation of mathematics». In [Rowe, McCleary 1989]

Stumpf, C. 1890. *Tonpsychologie*, Bd. II, Leipzig, Hirzel.

Tanabe, H. 1915. «Zur philosophischen Grundlegung der natürlichen Zahlen», *The Tôhoku Mathematical Journal*, éd. par T. Hayashi, vol. 7.

Tarski, A. 1944. The Semantic Conception of Truth and the Foundations of Semantics », *Philosophy and Phenomenological Research*, vol. 4 (3).

Tarski, A. 1956. *Logic, Semantics, Metamathematics*, Oxford, Clarendon Press.

Thiel, C. 1981. «A Portrait; or, How to Tell Frege from Schröder», *History and Philosophy of Logic* vol. 2 (1-2).

Thomson, J.J. 1987 (ed). *On Being and Saying: Essays in Honor of Richard Cartwright*, Cambridge, MIT Press.

Tieszen, R. 1989. *Mathematical Intuition: Phenomenology and Mathematical Knowledge*, Dordrecht, Kluwer.

Torretti, R. 1978. *Philosophy of Geometry from Riemann to Poincaré*, Dordrecht, Reidei.

Torretti, R. 1983. *Relativity and Geometry*, Oxford, Pergamon Press.

Tresmontant, C. 1964. *La métaphysique du Christianisme et la crise du treizième siècle*, Paris, PUF.

Vacca, G. 1915. «Sulle scoperte di Pietro Mengoli», *Atti dell'Accademia Nazionale dei Lincei-Rendiconti*, vol. XXIV, 5.

Vahlen, T. 1923. *Wert und Wesen der Mathematik*, Greifswald.

Vauthier, J. 1983. *Queen's Papers in Pure and Applied Mathematics*, 65, Kingston, Ontario.

Venn, J. 1980. *L'Algèbre de la logique*, 2ème éd., Paris, Blanchard.

Verriest, G. 1951. *Les nombres et les espaces*, Paris, Colin, 1951.

Voigt, A. 1892. «Was ist Logik?», *Vierteljahrsschrift für wissenschaftliche Philosophie,* Bd. XVI.

Wallace, D.F. 2009. *Die Entdeckung des Unendlichen: Georg Cantor und die Welt der Mathematik*, München, Piper Verlag.

Wang, H. 1957 (1962). «The Axiomatization of Arithmetic», *Survey of Mathematical Logic*, Peking, Science Press.

Weber, H. 1893. «Leopold Kronecker», *Jahresbericht der Deutschen Mathematiker-Vereinigung*, Bd. 2.

Wehmeier, K.F. 2004. «Russell's Paradox in Consistent Fragments of Frege's *Grundgesetze der Arithmetik*». In [Link 2004]

Winter, E. 1933. *Bernard Bolzano und sein Kreis*, Leipzig, Hegner.

Winter, E. 1949. *Leben und geistige Entwicklung des Sozialethikers und Mathematikers Bernard Bolzano 1781-1848*, Halle.

Wright, C. 1983. *Frege's Conception of Numbers as Objects*, Aberdeen University Press.

Worms, F. 2004. (ed.) *Le moment 1900 en philosophie*, Lille, Presses Univ. Septentrion

Zach, R. 2006. *Hilbert's program then and now*. In [Jacquette 2006]

Zermelo, E. 1904. «Beweis, daβ jede Menge wohlgeordnet werden kann», vol. 59, *Mathematische Annalen*.

Arkady Nedel
Paris 1 Sorbonne
Email: discover@free.fr